\documentclass[11pt,a4paper]{amsart}
\pdfoutput=1
\usepackage{amsfonts,amsmath,amscd,amssymb,amsbsy,amsthm,amstext,amsopn}
\usepackage{mathrsfs,subfigure,comment}
\usepackage{graphicx}
\usepackage{array}
\usepackage{multirow}
\usepackage[all]{xy}
\input xy
\xyoption{all}
\usepackage[colorlinks=true, pdfpagemode=none, pdfmenubar=false, pdfstartview=FitH, linkcolor=blue, citecolor=blue, urlcolor=blue, pdffitwindow=false]{hyperref}

\title{A derived approach to geometric McKay correspondence in dimension three}
\author{Sabin Cautis} 
\email{scautis@math.harvard.edu} 
\address{Department of Mathematics\\ Rice University \\ Houston, TX} 
 
\author{Timothy Logvinenko} 
\email{tlog@kth.se} 
\address{Department of Mathematics\\ KTH \\ Stockholm, Sweden}

\DeclareMathOperator{\homm}{Hom}
\DeclareMathOperator{\shhomm}{{\it Hom\rm}}

\DeclareMathOperator{\gsl}{SL}

\DeclareMathOperator{\picr}{Pic}

\DeclareMathOperator{\cl}{Cl}

\DeclareMathOperator{\spec}{Spec\;}

\DeclareMathOperator{\hilb}{Hilb}

\DeclareMathOperator{\ext}{Ext}
\DeclareMathOperator{\Ext}{Ext}
\DeclareMathOperator{\irr}{Irr}

\DeclareMathOperator{\supp}{Supp}

\DeclareMathOperator{\cohcat}{Coh}

\DeclareMathOperator{\lder}{\bf L}
\DeclareMathOperator{\rder}{\bf R}

\DeclareMathOperator{\hex}{Hex}
\DeclareMathOperator{\Except}{Exc}
\DeclareMathOperator{\sinksource}{SS}
\begin{document}

\def\bv{\mathbf{v}}
\def\kgc_{K^*_G(\mathbb{C}^n)}
\def\kgchi_{K^*_\chi(\mathbb{C}^n)}
\def\kgcf_{K_G(\mathbb{C}^n)}
\def\kgchif_{K_\chi(\mathbb{C}^n)}
\def\gpic_{G\text{-}\picr}
\def\gcl_{G\text{-}\cl}
\def\trch_{{\chi_{0}}}
\def\regring{{R}}
\def\regrep{{V_{\text{reg}}}}
\def\givrep{{V_{\text{giv}}}}
\def\lbar{{(\mathbb{Z}^n)^\vee}}
\def\genpx_{{p_X}}
\def\genpy_{{p_Y}}
\def\genpcn_{p_{\mathbb{C}^n}}
\def\gnat{gnat}
\def\twalg{{\regring \rtimes G}}
\def\L{{\mathcal{L}}}
\def\O{{\mathcal{O}}}
\def\gcd{\mbox{gcd}}
\def\lcm{\mbox{lcm}}
\def\tf{{\tilde{f}}}
\def\tD{{\tilde{D}}}

\def\mckquiv{\mbox{Q}(G)}
\def\C{{\mathbb{C}}}
\def\sF{{\mathcal{F}}}
\def\sW{{\mathcal{W}}}
\def\sL{{\mathcal{L}}}
\def\O{{\mathcal{O}}}
\def\Z{{\mathbb{Z}}}
\def\hmone{{\mathcal{W}}}

\theoremstyle{definition}
\newtheorem{defn}{Definition}[section]
\newtheorem*{defn*}{Definition}
\newtheorem{exmpl}[defn]{Example}
\newtheorem*{exmpl*}{Example}
\newtheorem{exrc}[defn]{Exercise}
\newtheorem*{exrc*}{Exercise}
\newtheorem*{chk*}{Check}
\newtheorem*{remarks*}{Remarks}
\theoremstyle{plain}
\newtheorem{theorem}{Theorem}[section]
\newtheorem*{theorem*}{Theorem}
\newtheorem{conj}[defn]{Conjecture}
\newtheorem{prps}[defn]{Proposition}
\newtheorem*{prps*}{Proposition}
\newtheorem{cor}[defn]{Corollary}
\newtheorem*{cor*}{Corollary}
\newtheorem{lemma}[defn]{Lemma}
\newtheorem*{claim*}{Claim}
\newtheorem{Specialthm}{Theorem}
\renewcommand\theSpecialthm{\Alph{Specialthm}}
\numberwithin{equation}{section}

\begin{abstract}
We propose a three dimensional generalization of the geometric McKay
correspondence described by Gonzales-Sprinberg and Verdier in
dimension two. We work it out in detail when $G$ is abelian 
and $\mathbb{C}^3/G$ has a single isolated singularity. More 
precisely, we show that the Bridgeland-King-Reid derived category 
equivalence induces a natural geometric correspondence between 
irreducible representations of $G$ and subschemes of the exceptional 
set of $G$-$\hilb (\mathbb{C}^3)$. This correspondence appears to be 
related to Reid's recipe.  
\end{abstract}

\maketitle

\section{Introduction} \label{section-intro}

The study of the McKay correspondence began with an observation by
John McKay in \cite{McKay-GraphsSingularitiesAndFiniteGroups} that
there exists a bijective correspondence between irreducible
representations of a finite subgroup $G \subseteq \gsl_2(\mathbb{C})$
and exceptional divisors of the minimal resolution $Y$ of
$\mathbb{C}^2/G$. Gonzales-Sprinberg and Verdier in
\cite{GsV-ConstructionGeometriqueDeLaCorrespondanceDeMcKay} gave a
geometric construction of this correspondence. The aim of this paper
is to give a generalization of this construction for dimension three. 
Our approach is via the derived McKay equivalence of \cite{BKR01} and
for $G$ abelian it appears to give a categorification of `Reid's recipe' from 
\cite{Kinosaki-97}.

The original construction of
\cite{GsV-ConstructionGeometriqueDeLaCorrespondanceDeMcKay} is as
follows. Denote by $K^G(\mathbb{C}^2)$ the Grothendieck ring of
$G$-equivariant coherent sheaves on $\mathbb{C}^2$ and by $K(Y)$ the
Grothendieck ring of $Y$. Let $\mathcal{M} \in \cohcat(Y \times
\mathbb{C}^2)$ be the structure sheaf of the reduced fiber product $Y
\times_{\mathbb{C}^2/G} \mathbb{C}^2$. Define the transform $\Theta:
K^G(\mathbb{C}^2) \rightarrow K(Y)$ by
\begin{align}\label{eqn-gspv-transform} \Theta (-) = \left(\pi_{Y *}
(\mathcal{M} \otimes \pi^*_{\mathbb{C}^2}(-))\right)^G \end{align}
where $\pi_Y$ and $\pi_{\mathbb{C}^2}$ are the projections from $Y
\times \mathbb{C}^2$ to $Y$ and $\mathbb{C}^2$. 

Denote by $\Except(Y)$ the set of irreducible exceptional divisors on
$Y$. Let $\rho$ be an irreducible representation of $G$. Since $\mathcal{M}$ 
is flat over $Y$ its pushforward to $Y$ is a
vector bundle and $\L_\rho := \Theta (\O_{\C^2} \otimes \rho^\vee)$ is its
$\rho$-eigensheaf. In general $\L_\rho$ is a vector bundle of
dimension $\mbox{dim}(\rho)$ and is called a \em tautological \rm or
\em GSp-V \rm sheaf. It is proven in
\cite{GsV-ConstructionGeometriqueDeLaCorrespondanceDeMcKay} that:
\begin{enumerate} 
\item \label{item-ktheory-isomorphism}
$K^G(\mathbb{C}^2) \xrightarrow{\Theta} K(Y)$ is an isomorphism of
abelian groups.  
\item \label{item-tautological-bundles}
$\{ c_1(\mathcal{L}_\rho) \;|\; \rho \in \irr(G) \mbox{ non-trivial } \}$ is the basis of $H^2(Y,\Z)$ dual to the basis $\{ E \;|\; \Except(Y)\}$ of $H_2(Y,\Z)$. This gives a $1$-to-$1$ correspondence between non-trivial $\rho$ and elements of $\Except(Y)$. Denote by $E_\rho$ the irreducible exceptional divisor dual to $c_1(\mathcal{L}_{\rho^\vee})$. 
\item \label{item-twisted-skyscraper} 
For each non-trivial $\rho \in \irr(G)$ we have $$\Theta(\mathcal{O}_{0} \otimes \rho) = [\mathcal{O}_{E_\rho}]$$ where $\mathcal{O}_0$ is the skyscraper sheaf
at the origin $(0,0) \in \mathbb{C}^2$.  
\end{enumerate}

Recall that a resolution $X' \xrightarrow{f} X$ is called crepant if
$f^* \omega_X = \omega_{X'}$ (see \cite{YPG87} for more details). In
dimension two the minimal resolution $Y \rightarrow \C^2/G$ is
the unique crepant resolution. In dimension three and above
crepant resolutions may not exist and are not necessarily unique. In 
particular $\mathbb{C}^3/G$ has a number of crepant resolutions
connected by a chain of flops.  

To figure out which crepant resolution to consider it helps to go back to the minimal resolution $Y \rightarrow \C^2/G$
and observe that $Y \cong G$-$\hilb(\mathbb{C}^2)$,
the fine module space of $G$-clusters in $\mathbb{C}^2$. A
$G$-cluster in $\mathbb{C}^n$ is a $G$-invariant subscheme $\Z
\subseteq \mathbb{C}^n$ such that $H^0(\mathcal{O}_Z)$ is
isomorphic to the regular representation $\regrep$ of $G$. The
structure sheaf $\mathcal{M}$ of the reduced fiber product $Y
\times_{\mathbb{C}^2/G} \mathbb{C}^2$ turns out to be the universal
family of $G$-clusters over $Y \times \mathbb{C}^2$. 

This suggests that for $G \subseteq \gsl_3(\mathbb{C})$ we should take $Y$ to be
$G$-$\hilb(\mathbb{C}^3)$ and $\mathcal{M}$ to be the universal
family of $G$-clusters over $Y \times \C^3$. Given a variety $X$
denote by $D(X)$ (resp. $D^G(X)$) the derived category of coherent
sheaves on $X$ (resp. $G$-equivariant coherent sheaves on $X$ if $X$
is equipped with an action of $G$). Define $\Phi: D(Y)
\rightarrow D^G(\C^3)$ to be the integral transform \begin{align}
\Phi (-) = \pi_{\mathbb{C}^3 *} \left(\mathcal{M}
\overset{\lder}{\otimes} \pi^*_{Y}(- \otimes \rho_0))\right)
\end{align} with kernel $\mathcal{M}$, where $(- \otimes \rho_0)$ is
the functor of equipping the sheaf with the trivial $G$-action. It
was proven in \cite{BKR01} that $\Phi$ is an equivalence of derived
categories and consequently that $Y$ is a crepant resolution of
$\C^3/G$. As derived equivalence implies $K$-group isomorphism this
gives a very satisfying proof of (\ref{item-ktheory-isomorphism}) in
dimension three. 

There has been some effort to generalise
(\ref{item-tautological-bundles}). One can define the vector bundles
$\L_\rho$ for $\rho \in \irr(G)$ as before. 
Since $\{ \O_{\C^3} \otimes \rho \}$ forms a basis for
$K^G(\C^3)$ the collection $\{ \L_\rho \}$ is a basis for $K(Y)$. 
Hence the Chern classes $c_1(\L_\rho)$ span $H^2(X,\Z)$. But they do 
not form a basis since there are relations. 

When $G$ is abelian the
bundles $\L_\rho$ are line bundles and $Y$ is a toric variety.
Following a conjecture by Reid, Craw showed in
\cite{Craw-AnexplicitconstructionoftheMcKaycorrespondenceforAHilbC3}
that certain $\mathcal{L}_\rho$ whose classes are redundant can be
replaced by abstract elements of $K(Y)$ in such a way that second
Chern classes of these `virtual' bundles give a basis of
$H^4(Y,\mathbb{Z})$.  This gives a $1$-to-$1$ correspondence between
$\irr(G)$ and a basis of $H^*(Y, \mathbb{Z})$. 

What dictates the way `virtual' bundles are produced is the so-called 
`Reid's recipe' (see Section 6 of \cite{Kinosaki-97}) and Section 3 of
\cite{Craw-AnexplicitconstructionoftheMcKaycorrespondenceforAHilbC3}).
It is a somewhat ad-hoc way, based on calculation of monomial bases of
each of the toric affine charts of $Y$, to assign a character of $G$
to each exceptional curve of the form $E \cap F$ for $E, F \in
\Except(Y)$ and a character or a pair of characters to each
exceptional divisor $E \in \Except(Y)$. 

We propose generalizing (\ref{item-twisted-skyscraper}).  Take the
inverse $\Psi: D^G(\C^3) \rightarrow D(Y)$ of the \cite{BKR01} derived
equivalence $\Phi$ and look at the images $\Psi(\mathcal{O}_0 \otimes
\rho) \in D(Y)$ which are supported on the exceptional locus of $Y$.
The correspondence $\rho \leftrightarrow \supp \Psi(\mathcal{O}_0
\otimes \rho)$ is what we think of as the geometric McKay 
correspondence in dimension three since in dimension two this
gives the classical McKay correspondence of
\cite{GsV-ConstructionGeometriqueDeLaCorrespondanceDeMcKay}.  

Apriori $\Psi(\mathcal{O}_0 \otimes \rho)$ will be a complex of sheaves but if $G$ is abelian we show that it is a sheaf: 

\begin{theorem}\label{theorem-images-are-pure-sheaves} Let $G \subset
\gsl_3(\mathbb{C})$ be a finite, abelian group such that
$\mathbb{C}^3/G$ has a single isolated singularity at the origin. Then, 
for any irreducible representation $\chi \in \irr(G)$, the Fourier-Mukai
transform $\Psi(\mathcal{O}_0 \otimes \chi)$ is a pure sheaf (i.e. some 
shift of a coherent sheaf).  
\end{theorem}

In Section \ref{section-inverse-transform-and-dual-family} 
we explain how to compute $\Psi(\mathcal{O}_0 \otimes
\rho)$ for any $G \subseteq \gsl_3(\mathbb{C})$. In Section
\ref{section-worked-example} we show how to explicitly calculate
the supports of $\Psi(\mathcal{O}_0 \otimes \rho)$ when $G$ is abelian
by using an explicit description of $G$-$\hilb(\mathbb{C}^3)$
(\cite{Craw-AnexplicitconstructionoftheMcKaycorrespondenceforAHilbC3})
and of the universal family over it
(\cite{Logvinenko-Natural-G-Constellation-Families}).

Based on such computational evidence, we present the following conjecture
as to the exact form of the sheaves $\Psi(\mathcal{O}_0 \otimes
\chi)$. A part of it gets proven in the course of proving Theorem
\ref{theorem-images-are-pure-sheaves}.  

\begin{conj} 
Let $G \subset \gsl_3(\mathbb{C})$ be a finite, abelian group. Then, 
for any irreducible representation $\chi \in \irr(G)$,
the Fourier-Mukai transform $\Psi(\mathcal{O}_0 \otimes \chi)$ is 
one of the following: 
\begin{enumerate} 
\item $\mathcal{L}^{-1}_{\chi} \otimes \mathcal{O}_{E_i}$ 
\item $\mathcal{L}^{-1}_{\chi} \otimes \mathcal{O}_{E_i \cap E_j}$ 
\item $\mathcal{F} [1] \quad\text{ where }\quad \supp_Y(\mathcal{F}) = 
E_{i_1} \cup \dots \cup E_{i_k}$ 
\item $\mathcal{O}_{Y}(\Except(Y)) \otimes \mathcal{O}_{\Except(Y)}[2]$
\end{enumerate} where $E_i$ are irreducible exceptional divisors, 
$\mathcal{F}$ a coherent sheaf and $\L_{\chi}$ is the tautological 
line bundle of \cite{GsV-ConstructionGeometriqueDeLaCorrespondanceDeMcKay}.
\end{conj}

We expect a similar picture to hold for non-abelian $G \subset \gsl_3(\C^3)$. In particular:

\begin{conj} 
For any finite $G \subset \gsl_3(\C)$ the images $\Psi(\mathcal{O}_0 \otimes \chi)$ are pure sheaves. 
\end{conj}

\subsection{An Example}

\begin{figure}[h] \begin{center}
\includegraphics[scale=0.25]{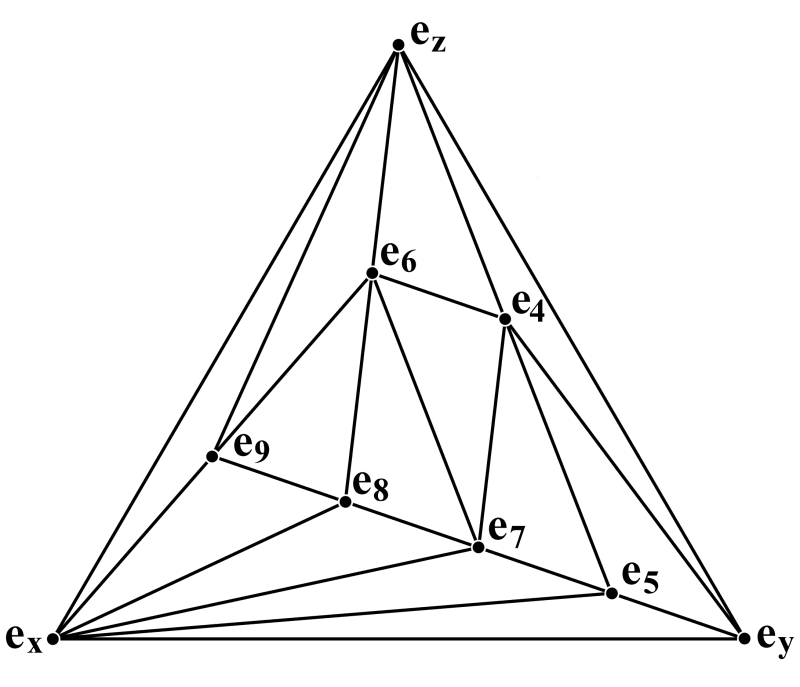} \end{center}
\caption{\label{figure-03} Toric fan of $G$-$\hilb(\mathbb{C}^3)$}
\end{figure}

Let $G = \frac{1}{13}(1,5,7) \subset \gsl_3(\C)$. For the definition 
see Section \ref{section-worked-example} which works out
this example in detail. Then $G$-$\hilb(\mathbb{C}^3)$ is
defined by the toric fan in Figure $\ref{figure-03}$. 

Computing the images $\Psi(\mathcal{O}_0 \otimes \chi)$ for each $\chi
\in \irr(G)$ we summarise the results in the following table.

\begin{footnotesize} \begin{align}
\label{table-13-1-5-7-supports-of-the-images}
\setlength{\extrarowheight}{0.2cm} \begin{array}{|r|r|r|r|} \hline
\chi & \supp H^{-2}(\Psi(\mathcal{O}_0 \otimes \chi)) & \supp
H^{-1}(\Psi(\mathcal{O}_0 \otimes \chi)) & \supp
H^{-0}(\Psi(\mathcal{O}_0 \otimes \chi)) \\ \hline \chi_0 & E_4 \cup
E_5 \cup E_6 \cup E_7 \cup E_8 \cup E_9 & 0 & 0 \\ \chi_1 & 0 & 0 &
E_4 \\ \chi_2 & 0 & 0 & E_7 \\ \chi_3 & 0 & E_7 & 0 \\ \chi_4 & 0 & 0
& E_6 \\ \chi_5 & 0 & 0 & E_9 \\ \chi_6 & 0 & E_4 \cup E_6 \cup E_9 &
0 \\ \chi_7 & 0 & 0 & E_5 \\ \chi_8 & 0 & E_4 \cup E_5 & 0 \\ \chi_9 &
0 & 0 & E_6 \cap E_7 \\ \chi_{10} & 0 & 0 & E_8 \\ \chi_{11} & 0 & E_6
\cup E_8 & 0 \\ \chi_{12} & 0 & E_5 \cup E_7 \cup E_8 \cup E_9 & 0 \\
\hline \end{array} \end{align} \end{footnotesize} 

\subsection{Relation to Reid's recipe}

From Proposition 9.3 of \cite{Craw-Ishii-02} we know that if $\chi$ marks
an irreducible exceptional divisor $E$ by Reid's recipe then
$\Psi(\O_0 \otimes \chi) \cong \L^{-1}_{\chi} \otimes \O_E$. 

\begin{figure}[h] \begin{center}
\includegraphics[scale=0.25]{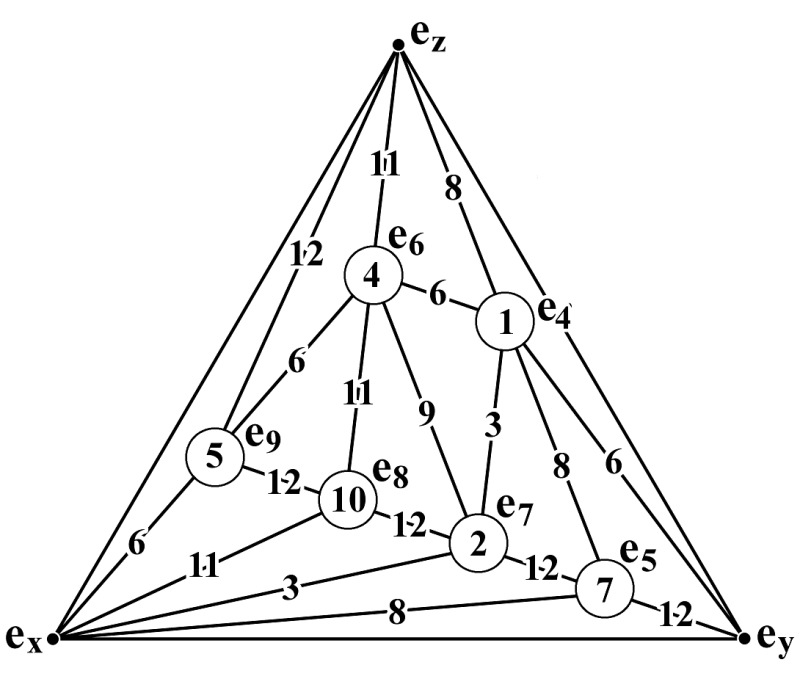} \end{center}
\caption{\label{figure-28} Reid's recipe marking for $G =
\frac{1}{13}(1,5,7)$.} \end{figure}

However it is worth noting that in all the examples we computed the
supports of $\Psi(\mathcal{O}_0 \otimes \chi)$ appear to be more
strongly related to Reid's recipe. As an example, we give in Figure
\ref{figure-28} the markings according to Reid's recipe for $G =
\frac{1}{13}(1,5,7)$ and invite the reader to compare it to the data
in Table \ref{table-13-1-5-7-supports-of-the-images}. For instance,
$\Psi(\mathcal{O}_0 \otimes \chi_{12})$ is supported on $E_5 \cup E_7
\cup E_8 \cup E_9$ -- precisely those divisors which contain
two of the curves marked by $\chi_{12}$ by Reid's recipe.

We hope to prove the following:

\begin{conj}
Let $G \subset \gsl_3(\mathbb{C})$ be a finite, abelian group.
In Reid's recipe, \begin{enumerate} \item if $\chi$ marks a
divisor then $\Psi(\mathcal{O}_0 \otimes \chi)$ is supported in degree
$0$ on that divisor \item if $\chi$ marks several irreducible curves
then $\Psi(\mathcal{O}_0 \otimes \chi)$ is supported in degree $-1$ on
the union of all divisors containing two or more of these curves \item
if $\chi$ marks a single irreducible curve then $\Psi(\mathcal{O}_0
\otimes \chi)$ is supported in degree $0$ on that curve
\end{enumerate} \end{conj}

\subsection{Further Remarks}

\begin{enumerate}

\item \label{item-doesnt-hold-for-gcons} Theorem \ref{theorem-images-are-pure-sheaves} does not
necessarily hold if you replace $Y = G$-$\hilb(\C^3)$ by 
another crepant resolution of $\mathbb{C}^3/G$. For example,
let $G = \frac{1}{11}(1,3,7)$ and let $Y$ be the fine 
module space $M_\theta$ of $\theta$-stable $G$-constellations for 
$\theta = (13, -4, 5, -4, -4, 1, -4, 9, -4, -4, -4)$ (see
\cite{Craw-Ishii-02} for more detail). Set $\mathcal{M}$ to be 
the universal family on $Y \times \mathbb{C}^3$. 
Compute the toric fan of $M_\theta$ as described in
\cite{Craw-Maclagan-Thomas-05-I}, Section 6 and compute $\mathcal{M}$ 
as described in \cite{Logvinenko-DerivedMcKayCorrespondenceViaPureSheafTransforms}, Section 4.5. 
We can then compute $\Psi(\mathcal{O}_0 \otimes \chi_0)$
just like in Section \ref{section-worked-example} and it turns out
to have non-zero cohomologies both in degree $-2$ and in degree $-1$. 

If one examines the proof of Theorem \ref{theorem-images-are-pure-sheaves} then one would find that Propositions \ref{prps-just-one-out-of-three} and \ref{prps-gcluster-no-three-arrows} are the basic results which depend on the fact we are using (the universal family over) $G$-$\hilb(\C^3)$ instead of some other crepant resolution. 

\item We prove Theorem \ref{theorem-images-are-pure-sheaves} by 
doing some very explicit computations with the representation
$\mckquiv_{\mathcal{M}}$ of the McKay quiver $\mckquiv$
associated to the universal family $\mathcal{M}$. It would be nice 
to obtain a more conceptual proof of this result. 

For example, Alexei Bondal suggested the following argument why Theorem
\ref{theorem-images-are-pure-sheaves} might hold. The
idea is to use his result which states that any exceptional object on
a Del-Pezzo surface is a sheaf. On the other hand, there is the
following result of Seidel and Thomas. If $i: X \rightarrow Y$ 
is an inclusion of a divisor satisfying $i^* \omega_Y = \O_X$ then 
$i_* \sF$ is spherical if and only if $\sF$ is exceptional.

In our case we take $Y = G$-$\hilb(\C^3)$ and pretend that the union $X \subset Y$ of its
exceptional divisors is a Del-Pezzo surface. Now $\O_0 \otimes \chi$ is a spherical object in $D^G(\C^3)$ so
$\Psi(\O_0 \otimes \chi)$ is a spherical object
in $Y$. If one could prove directly that $\Psi(\O_0 \otimes \chi)$ is 
the pushforward of an object in $D(X)$ then this object should be 
a sheaf by Bondal's result. At present, such a direct proof is 
beyond us. Note that in light of the remark 
\ref{item-doesnt-hold-for-gcons} it would need to use 
the fact that $\mathcal{M}$ is the family of $G$-clusters and 
not of arbitrary $G$-constellations.
 
\item Finally, we see in Section \ref{section-inverse-transform-and-dual-family} that $\Psi: D^G(\C^3) \rightarrow D(Y)$, the inverse of \cite{BKR01} transform $\Phi$, is the transform defined by the dual of the universal family $\mathcal{M}$. We could have instead taken the transform defined by $\mathcal{M}$ itself, as was done in \cite{GsV-ConstructionGeometriqueDeLaCorrespondanceDeMcKay}. This choice is only for convenience. In both cases the transforms of $\O_0 \otimes \chi$ are (shifted) sheaves. However, using our choice the structure sheaves of irreducible divisors appear in degree $0$, while $\Psi(\mathcal{O}_0 \otimes \rho_0)$ appears in degree $-2$ and not vice versa. Moreover, this choice matches up better with the marking given by Reid's recipe. 

\end{enumerate}

The paper is organized as follows. In section
\ref{section-inverse-transform-and-dual-family} we describe the
transform $\Psi(\O_0 \otimes \rho)$ as a cube complex. In section
\ref{section-cohomology-of-skew-commutative-cubes} we compute the
cohomology of a general cube complex. In section
\ref{section-the-mckay-quiver} we introduce the associated McKay
quiver representation and relate it to the cube complex from section
\ref{section-inverse-transform-and-dual-family}. We also prove some
properties which it satisfies. In section \ref{section-main-results}
we use these properties to prove Theorem \ref{theorem-images-are-pure-sheaves}. Finally, in section
\ref{section-worked-example} we explicitly work out $\Psi(\O_0 \otimes
\rho)$ when $G = \frac{1}{13}(1,5,7)$.

\bf Acknowledgments: \rm We would like to thank Alexei Bondal, Torsten
Ekedahl and Alastair Craw for useful discussions and suggestions. This
paper was conceived during the authors' stay at Institut
Mittag-Leffler (Sweden) during its "Moduli Spaces" program in
2006-2007. The authors would like to thank the organizers of this
program as well as everyone else at the institute for their
hospitality. The first author would also like to thank the mathematics
department at Rice University. The second author did most of his work
on this paper during his stay at KTH (Sweden) and would like to thank
them for their support. 

\section{The inverse transform and the dual family}
\label{section-inverse-transform-and-dual-family}

\subsection{Notation}

In this section, let $G$ be an arbitrary finite subgroup of
$\gsl_n(\mathbb{C})$ and $Y$ a smooth $n$-dimensional separable scheme
of finite type over $\mathbb{C}$. We equip $Y$ with the trivial
$G$-action. Then $G$ acts naturally on $Y \times \mathbb{C}^n$ and we
can consider the bounded derived category of $G$-equivariant coherent
sheaves $D^G(Y \times \mathbb{C}^n)$. 

We denote by $\givrep$ the representation of $G$ induced by its
inclusion into $\gsl_n(\mathbb{C})$ and by $\regring$ the symmetric
algebra $S(\givrep^\vee)$. We identify $\mathbb{C}^n$ with the affine
$G$-scheme $\spec \regring$. We also call a $G$-equivariant sheaf a
$G$-sheaf for short (cf. \cite{BKR01}, Section 4).  

There is an equivalence between the category of quasi-coherent
$G$-sheaves on $\mathbb{C}^n$ and of $\twalg$-modules which is induced
by the functor $\Gamma$ of taking global sections. The relative
version of this is an equivalence between the category of
quasi-coherent $G$-sheaves on $Y \times \mathbb{C}^n$ and of
quasi-coherent sheaves of $(\twalg) \otimes_\mathbb{C}
\mathcal{O}_Y$-modules induced by the pushdown functor $\pi_{Y*}$ to
$Y$. 

For an object $\sF \in D^G(Y \times \mathbb{C}^n)$ we denote by $\Phi_{\sF}$
the integral transform $D(Y) \rightarrow D^G(\mathbb{C}^n)$ defined by
\begin{align} 
\Phi_{\sF} (-) = \pi_{\mathbb{C}^3 *} \left(\sF \overset{\lder}{\otimes} \pi^*_{Y}(- \otimes \rho_0))\right)
\end{align} 
where $(- \otimes \rho_0)$ is the functor of equipping the
sheaf with the trivial $G$-action. And we denote by $\Psi_{\sF}$ the
integral transform in the opposite direction $D^G(\mathbb{C}^n)
\rightarrow D(Y)$ defined by 
\begin{align} 
\label{eqn-transform-psi}
\Psi_{\sF} (-) = \pi_{Y *} \left(\sF \overset{\lder}{\otimes} \pi^*_{\mathbb{C}^n}(-)\right)^G 
\end{align}

In brief, a $\gnat$-family is a family on $Y \times \C^n$, flat
over $Y$, of certain finite-length sheaves on $\C^n$ called
$G$-constellations, each of which is supported on a single $G$-orbit. 
Sometimes we will abuse notation by identifying a $\gnat$-family with its
pushforward to $Y$ via $\pi_{Y*}$ as a $(\twalg) \otimes_{\C}
\mathcal{O}_Y$-module. The basic example to keep in mind is the
pushforward to $Y = G$-$\hilb(\C^n)$ of the universal family of
$G$-clusters. For the precise definition see
\cite{Logvinenko-DerivedMcKayCorrespondenceViaPureSheafTransforms},
Section 3.2 and \cite{Logvinenko-Natural-G-Constellation-Families}.  

Suppose $\mathcal{F}$ is a $\gnat$-family on $Y \times \C$. It was shown in
\cite{Logvinenko-DerivedMcKayCorrespondenceViaPureSheafTransforms},
Lemma 4 that the left adjoint to $\Phi_\mathcal{F}$ is the integral
transform $\Psi_{\mathcal{F}^\vee[n]}$.  In this section, we show
that $\mathcal{F}^\vee[n]$ is the $\gnat$-family $\tilde{\mathcal{F}}$
dual to $\mathcal{F}$. Notice that $\tilde{\mathcal{F}}$ is a
sheaf (i.e. a complex concentrated in degree zero).  We then use this
to compute $\Psi_{\mathcal{F}^\vee[n]}(\mathcal{O}_0 \otimes \rho)$.
The application we have in mind is the case $n=3$ with $Y = G$-$\hilb
(\mathbb{C}^3)$ and $\mathcal{F}$ the universal family of $G$-clusters
on $Y \times \mathbb{C}^3$.  

\subsection{The dual family}

If $W$ is a $\twalg$-module, then its dual $W^\vee =
\homm_\mathbb{C}(W,\mathbb{C})$ is naturally a $\twalg$-module via
\begin{align} \label{eqn-natural-twalg-module-structure-on-dual} g
\cdot \alpha(\bullet) = \alpha(g^{-1}\cdot \bullet) \quad \text{ and }
\quad m \cdot \alpha(\bullet) = \alpha(m\cdot \bullet) \end{align}
where $\alpha \in \homm_\mathbb{C}(W,\mathbb{C}), g \in G$ and $m \in
\regring$. 

Given a finite-length $G$-sheaf $\mathcal{V}$ on $\C^n$ we can define
its {\tt dual} $\mathcal{V}^\vee$ to be the $G$-sheaf corresponding to
the dual of the $\twalg$-module $\Gamma(\mathcal{V})$. Similarly, 
if $\mathcal{F} \in \cohcat^G(Y \times \C^n)$ is a
$\gnat$-family then the corresponding family of
$\twalg$-modules is $\pi_{Y *} \mathcal{F}$. We define the {\tt dual
$\gnat$-family} $\tilde{\mathcal{F}}$ to be the $\gnat$-family
corresponding to $(\pi_* \sF)^\vee = \shhomm_{Y}(\pi_{Y*} \mathcal{F}, \O_Y)$.

\begin{prps} \label{prps-Fdual} If $\mathcal{F}$ is a $\gnat$-family
on $Y \times \mathbb{C}^n$ then $\mathcal{F}^\vee[n] \cong
\tilde{\mathcal{F}}$.  \end{prps} 

\begin{proof} 
We would like to use Grothendieck-Serre duality for the morphism 
$\pi_Y: Y \times \mathbb{C}^n \rightarrow Y$. But $\pi_Y$ is
manifestly non-proper. However, using Deligne's approach via
compactification and Pro-categories described in 
\cite{DeligneCohomologieASupportPropre} we can still obtain the
duality for the full subcategory of $D^G(Y \times \mathbb{C}^n)$
consisting of objects whose support is proper over $Y$. Specifically,
we get a natural isomorphism
\begin{align} \label{eqn-relative-grothendieck-duality-with-compact-supports}
\pi_{Y *} \rder \shhomm_{Y \times \mathbb{C}^n}(A, \pi^*_Y B \times 
\omega_{Y \times \mathbb{C}^n / Y}[n]) \xrightarrow{\sim}
\rder \shhomm_Y(\pi_{Y *} A, B)
\end{align}
for any $A \in D^G(Y \times \mathbb{C}^n)$ whose support is proper
over $Y$ and any $B \in D^G(Y)$. Observe that as $\pi_{Y}$ is both
flat and affine, $\pi^*_Y$ and $\pi_{Y *}$ are both exact. 

Setting $A = \mathcal{F}$ and $B = \mathcal{O}_Y$ yields
\begin{align}
\pi_{Y*}\mathcal{F}^\vee[n] \xrightarrow{\sim} \rder \shhomm_Y(\pi_{Y
*} \mathcal{F}, \mathcal{O}_Y).
\end{align}
since $\omega_{Y \times \mathbb{C}^n /Y}$ is trivial. Because the right-hand side is by definition $\pi_{Y *} \tilde{\mathcal{F}}$ the claim follows.
\end{proof}

For the skeptical reader unconvinced of the validity of \eqref{eqn-relative-grothendieck-duality-with-compact-supports} we give an alternative proof of Proposition \ref{prps-Fdual} in the Appendix. It computes $\mathcal{F}^\vee$ by resolving $\mathcal{F}$ by locally free sheaves. To do this we give a natural locally free resolution for finite-length sheaves on $\mathbb{C}^n$ that may be of independent interest. 

\subsection{The image $\Psi(\mathcal{O}_0 \otimes \rho) \in D(Y)$}

Let $\mathcal{O}_0$ be the structure sheaf of the origin $0 \in
\mathbb{C}^n$. For any irreducible representation $\rho \in \irr G$
denote by $\mathcal{O}_0 \otimes \rho$ the $G$-sheaf where
$\mathcal{O}_{\mathbb{C}^n}$ acts on the first factor and $G$ on the
second.

The proposition below shows how to compute
$\Psi_{\mathcal{F}}(\mathcal{O}_0 \otimes \rho)$ for any
$\gnat$-family $\mathcal{F}$ on $Y$. In view of Proposition
\ref{prps-Fdual} we can compute
$\Psi_{\mathcal{F}^\vee[n]}(\mathcal{O}_0 \otimes \chi)$ by setting
$\mathcal{F}$ in the proposition below to be the dual family
$\tilde{\mathcal{F}}$. 

\begin{prps} \label{prps-preimages-of-point-reps-general-case} Let
$\mathcal{F}$ be any $\gnat$-family over $Y$.  The image of
$\mathcal{O}_0 \otimes \rho$ under the integral transform
$\Psi_{\mathcal{F}}:\; D^G(\mathbb{C}^n) \rightarrow D(Y)$ is the
$G$-invariant part of the complex \begin{align}
\label{eqn-preimages-of-point-reps-general-case-complex} \wedge^{n}
\givrep^\vee \otimes \pi_{Y*}\mathcal{F} \otimes \rho
\xrightarrow{\delta_n} \dots \xrightarrow{\delta_{k+1}} \wedge^{k}
\givrep^\vee \otimes \pi_{Y*}\mathcal{F} \otimes \rho
\xrightarrow{\delta_k} \dots \xrightarrow{\delta_1}
\pi_{Y*}\mathcal{F} \otimes \rho \end{align} where the rightmost term
lies in degree zero. The differentials $\delta_k$ are defined for $m_i
\in \givrep^\vee$, $s \in \pi_{Y *}\mathcal{F}$ and $v \in \rho$ by
\begin{align}\label{eqn-preimage-of-a-point-rep-differential-general}
(m_1 \wedge \dots \wedge m_k) \otimes s \otimes v \; \mapsto \;
\sum_{i = 1}^{k} (-1)^{i}  (\dots m_{i-1} \wedge m_{i+1} \dots)
\otimes (m_i \cdot s) \otimes v \end{align} \end{prps}
 
\begin{proof} To compute $\Psi_{\mathcal{F}}(\mathcal{O}_0 \otimes
\rho) = \pi_{Y*}(\mathcal{F} \overset{\lder}{\otimes}_{Y \times
\mathbb{C}^n} \pi^*_{\mathbb{C}^n}(\mathcal{O}_0 \otimes \rho))^G$ we
first resolve $\mathcal{O}_0 \otimes \rho$ via a Koszul resolution.
Observe that elements of $\givrep^\vee$ are the non-constant linear
functions on $\mathbb{C}^n$ and so any basis of $\givrep^\vee$ 
generates the ideal of $0 \in \mathbb{C}^n$. 
We therefore obtain the complex \begin{align} \wedge^n \givrep^\vee \otimes
\mathcal{O}_{\mathbb{C}^n} \otimes \rho \rightarrow \dots \rightarrow
\wedge^k \givrep^\vee \otimes \mathcal{O}_{\mathbb{C}^n} \otimes \rho
\rightarrow \dots \rightarrow \mathcal{O}_{\mathbb{C}^n} \otimes \rho
\end{align} with the differential maps \begin{align}
\label{eqn-differentials-in-images-of-point-reps} (m_1 \wedge \dots
\wedge m_k) \otimes s \otimes v \quad \mapsto \quad \sum_{i = 1}^{k}
(-1)^{i+1} (\dots m_{i-1} \wedge m_{i+1} \dots) \otimes  m_i \cdot s \otimes v 
\end{align} 
where $m_i \in \givrep^\vee$, $s \in \mathcal{O}_{\mathbb{C}^n}$ and 
$v \in \rho$. 

Pulling back to $Y \times \mathbb{C}^n$ and tensoring with
$\mathcal{F}$ we get \begin{align} \wedge^n \givrep^\vee \otimes
\mathcal{F} \otimes \rho \rightarrow \dots \rightarrow \wedge^k
\givrep^\vee \otimes \mathcal{F} \otimes \rho \rightarrow \dots
\rightarrow \mathcal{F} \otimes \rho \end{align} The result now
follows by applying $\pi_{Y *}$.  \end{proof}

\subsubsection{The abelian case}
\label{section-inverse-transform-abelian-case}

If the group $G$ is abelian then we can describe the results above
more explicitly. This is because we can decompose any $\gnat$-family
$\pi_{Y*} \mathcal{F}$ on $Y$ into eigensheaves as follows. 

Any $\gnat$-family $\mathcal{F}$ defines a Hilbert-Chow morphism
$\pi_\mathcal{F}:\; Y \rightarrow\mathbb{C}^n/G$ which sends any $y
\in Y$ to the $G$-orbit in $\mathbb{C}^n$ that the fiber
$\mathcal{F}_{|y}$ is supported on. We can then define
the notion of $G$-Weil divisors on $Y$ as in
\cite{Logvinenko-Natural-G-Constellation-Families}, Section 2. Any
$\gnat$-family on $Y$ can be decomposed as $\pi_{Y*} \mathcal{F} =
\bigoplus_{\chi \in G^\vee} \mathcal{L}(-D_\chi)$ for a unique set of
$G$-Weil divisors $\{D_{\chi}\}$ with $D_{\chi_0} = 0$ (see
\cite{Logvinenko-Natural-G-Constellation-Families}, Section 3.1). A good thing
to remember is that $G$ acts on $\mathcal{L}(D_\chi)$ by the character $\chi$ and
hence on $\mathcal{L}(-D_\chi)$ by $\chi^{-1}$. 

Fix a basis of $G$-eigenvectors for $\givrep^\vee$ and denote it by
$x_1, \dots, x_n$. This choice determines a $G$-eigenvector basis $I_k
= \{x_{i_1} \wedge \dots \wedge x_{i_k}: i_1 < \dots < i_k \}$ for
each of the spaces $\wedge^k \givrep^\vee$. Denote by $\kappa(x_i)$
the character with which $G$ acts on $x_i$.  More generally, given any
element $e = x_{i_1} \wedge \dots \wedge x_{i_k}$ of $I_k$ denote by
$\kappa(e)$ the character $\kappa(x_{i_1})\dots \kappa(x_{i_k})$, and
similarly for Laurent monomials in $x_i$. 

We can now give a refinement of Proposition
\ref{prps-preimages-of-point-reps-general-case}. 

\begin{prps} \label{prps-preimages-of-point-reps-abelian-case} Suppose
$G$ is abelian, let $\chi$ be a character of $G$ and let $\mathcal{F}
= \bigoplus_{\chi \in G^\vee} \mathcal{L}(-D_\chi)$ be a
$\gnat$-family. Then $\Psi_{\mathcal{F}}(\mathcal{O}_0 \otimes \chi)
\in D(Y)$ is equivalent to the complex \begin{align}
\label{eqn-preimage-of-point-rep-complex-abelian} \bigoplus_{e \in
I_{n}} \mathcal{L}(-D_{\kappa(e) \chi}) \rightarrow \dots \rightarrow
\bigoplus_{e \in I_{k}} \mathcal{L}(-D_{\kappa(e) \chi}) \rightarrow
\dots \rightarrow \bigoplus_{e \in I_{0}} \mathcal{L}(-D_{\kappa(e)
\chi}) \simeq \mathcal{L}(-D_{\chi}) \end{align} where, for each $j
\in \{1,\dots,k\}$, the differential  \begin{align}
\label{eqn-different-summands-in-point-rep-preimage-complex}
\mathcal{L}(-D_{\kappa(x_{i_1} \wedge \dots \wedge x_{i_k}) \chi})
\rightarrow \mathcal{L}(- D_{\kappa(\dots x_{i_{j-1}} \wedge
x_{i_{j+1}} \dots) \chi}) \end{align} is given by $s \mapsto (-1)^j
x_{i_j}\cdot s$.  \end{prps} \begin{proof} From decomposition
$\mathcal{F} = \bigoplus_{\chi' \in G^\vee} \mathcal{L}(-D_{\chi'})$
and $\wedge^k \givrep^\vee = \bigoplus_{e \in I_k} \mathbb{C} e$ we
get \begin{align*} \wedge^{k} \givrep^\vee \otimes \pi_{Y*}\mathcal{F}
\otimes \chi = \bigoplus_{\chi' \in G^\vee, e \in I_{k}} \mathbb{C} e
\otimes \mathcal{L}(-D_{\chi'}) \otimes \chi \end{align*} The group
$G$ acts on $\mathbb{C} e \otimes \mathcal{L}(D_{\chi'}) \otimes \chi$
by the character $\kappa(e) (\chi')^{-1} \chi$, so the $G$-invariant
part consists of terms where $\chi' = \kappa(e)\chi$. This gives
\eqref{eqn-preimage-of-point-rep-complex-abelian}. The claim about the
differentials follows immediately from
\eqref{eqn-preimage-of-a-point-rep-differential-general}.  \end{proof}

To compute $\Psi_{\tilde{\mathcal{F}}}(\mathcal{O}_0 \otimes \chi)$ 
we need the result below which allows us to express the $G$-Weil
divisors defining the dual
family $\tilde{\mathcal{F}}$ in terms of those defining $\mathcal{F}$:

\begin{prps} \label{prps-the-dual-family}
If $\mathcal{F} = \bigoplus_{\chi \in G^\vee} \mathcal{L}(-D_\chi)$ then 
$\tilde{\mathcal{F}} = \bigoplus_{\chi \in G^\vee} \mathcal{L} (-D'_{\chi})$ 
where $D'_\chi = - D_{\chi^{-1}}$.
\end{prps} 

\begin{proof} For any two $G$-Weil divisors $A$ and $B$
there is a standard $\mathcal{O}_Y$-module isomorphism
$\mathcal{L}(B-A) \xrightarrow{\sim}
\shhomm_{\mathcal{O}_Y}(\mathcal{L}(A), \mathcal{L}(B))$. The sheaves
$\mathcal{L}(A), \mathcal{L}(B)$ and $\mathcal{L}(B-A)$ all come with
a natural embedding into the constant sheaf $K(\mathbb{C}^n)$ on $Y$
and the isomorphism sends a section of $\mathcal{L}(B-A)$ to the map
of the multiplication by that section inside of $K(\mathbb{C}^n)$.

In particular, this gives an $\mathcal{O}_Y$-module isomorphism
$$\phi:\; \bigoplus_{\chi \in G^\vee} \mathcal{L}(D_{\chi})
\xrightarrow{\sim} \shhomm_{\mathcal{O}_{Y}}(\bigoplus_{\chi \in
G^\vee}\mathcal{L}(-D_\chi), \mathcal{O}_Y) = \tilde{\mathcal{F}}.$$
We claim that $\phi$ is $\twalg$-equivariant, thus making it an
isomorphism of $\mathcal{O}_Y \otimes (\twalg)$ modules. Let $s$ be a
section of $\bigoplus_{\chi \in G^\vee} \mathcal{L}(D_{\chi^{-1}})$
and $t$ a section of $\bigoplus_{\chi \in G^\vee}
\mathcal{L}(-D_{\chi})$. Then for any $g \in G$ and $m \in \regring$
\begin{align*} &\phi(r \cdot s)(t) = rst = r \cdot (\phi(s)(t)) =
(r\cdot \phi(s))(t)  \\ &\phi(g \cdot s)(t) = (g\cdot s) t = g \cdot (
s (g^{-1}\cdot t)) \overset{(1)}{=} s (g^{-1}\cdot t) =
\phi(s)(g^{-1}\cdot t) = (g \cdot \phi(s))(t) \end{align*} where the
equality $(1)$ is due to the fact that $s (g^{-1} \cdot t)$ is an
element of $\mathcal{O}_Y \subseteq K(\mathbb{C}^n)$ and thus
$G$-invariant.
Finally, since $G$ acts by $\chi$ on $\mathcal{L}(D_\chi)$ it must 
be the summand $\mathcal{L}(-D'_{\chi^{-1}})$ of $\tilde{\mathcal{F}}$. 
This yields $D'_{\chi} = -D_{\chi^{-1}}$.
\end{proof}

\begin{cor} \label{cor-preimages-of-point-reps-abelian-case} Suppose
$G$ is abelian, $\chi$ is a character of $G$ and let $\mathcal{F} =
\bigoplus_{\chi \in G^\vee} \mathcal{L}(-D_\chi)$ be a $\gnat$-family.
Then $\Psi_{\mathcal{\tilde{F}}}(\mathcal{O}_0 \otimes \chi) \in D(Y)$
is equivalent to the complex \begin{align} \bigoplus_{e \in I_{n}}
\mathcal{L}(D_{\kappa(e)^{-1} \chi^{-1}}) \rightarrow \dots
\rightarrow \bigoplus_{e \in I_{k}} \mathcal{L}(D_{\kappa(e)^{-1}
\chi^{-1}}) \rightarrow \dots \rightarrow \bigoplus_{e \in I_{0}}
\mathcal{L}(D_{\kappa(e)^{-1} \chi^{-1}}) \simeq
\mathcal{L}(D_{\chi^{-1}}) \end{align} where, for each $j \in
\{1,\dots,k\}$, the differential \begin{align}
\mathcal{L}(D_{\kappa(x_{i_1} \wedge \dots \wedge x_{i_k})^{-1}
\chi^{-1}}) \rightarrow \mathcal{L}(D_{\kappa(\dots x_{i_{j-1}} \wedge
x_{i_{j+1}} \dots)^{-1} \chi^{-1}}) \end{align} is given by $s \mapsto
(-1)^j x_{i_j} \cdot s$.  \end{cor}

\subsection{Examples}

We end this section with a couple of examples. 

\begin{exmpl} Let $G$ be the group $\frac{1}{7}(1,6)$. That is, $G
\cong \Z/7\Z$ embedded in $\gsl_2(\mathbb{C})$ via $1 \mapsto \left(
\begin{smallmatrix} \xi & \\ & \xi^{6} \end{smallmatrix}\right)$.
Denote by $\chi_i$ the character of $G$ given by $i \mapsto \xi^i$. In
these terms, $\kappa(x_1) = \chi_6$ and $\kappa(x_2) = \chi_1$. Let
$Y$ be the minimal resolution of $\mathbb{C}^2/G$ and let $\mathcal{F}
= \bigoplus \mathcal{L}(-D_\chi)$ be any $\gnat$-family on $Y$. Then
$\Psi_{\tilde{\mathcal{F}}}(\mathcal{O}_0 \otimes \chi_3) \in D(Y)$ is
given by the (total complex of the) square \begin{align*} \xymatrix{ &
\mathcal{L}(D_{\chi_5}) \ar"2,3"^{-x_1 \cdot} &   \\
\mathcal{L}(D_{\chi_4}) \ar"1,2"^{x_2 \cdot} \ar"3,2"_{-x_1 \cdot}  &
& \mathcal{L}(D_{\chi_4}) \\ & \mathcal{L}(D_{\chi_3}) \ar"2,3"_{x_2
\cdot} & } \end{align*} By convention the rightmost term is sitting in
degree $0$.  \end{exmpl}

\begin{exmpl} Let $G$ be the subgroup $\frac{1}{13}(1,5,7)$, $Y$ a
crepant resolution of $\mathbb{C}^3/G$ and $\mathcal{F} = \bigoplus
\mathcal{L}(-D_\chi)$ a $\gnat$-family on $Y$. Then the
$\Psi_{\tilde{\mathcal{F}}}(\mathcal{O}_0 \otimes \chi_5) \in D(Y)$ is
given by the (total complex of the) cube \begin{align}
\label{exmpl-transform-of-reptwist-3d-case} \xymatrix{ &
\mathcal{L}(D_{\chi_7}) \ar"2,5"^<<{x_3 \cdot} \ar"3,5"_<<<{-x_2
\cdot} & & & \mathcal{L}(D_{\chi_9}) \ar"2,6"^{-x_1 \cdot} & \\
\mathcal{L}(D_{\chi_8}) \ar"1,2"^{-x_1 \cdot} \ar"2,2"^{x_2 \cdot}
\ar"3,2"_{- x_3 \cdot} & \mathcal{L}(D_{\chi_3}) \ar"1,5"^>>>>{x_3
\cdot} \ar"3,5"_>>>>>{-x_1 \cdot} & & & \mathcal{L}(D_{\chi_0})
\ar"2,6"^{-x_2 \cdot} & \mathcal{L}(D_{\chi_8}) \\ &
\mathcal{L}(D_{\chi_1}) \ar"1,5"^<<<{x_2 \cdot} \ar"2,5"_<<{-x_1
\cdot} & & & \mathcal{L}(D_{\chi_{2}}) \ar"2,6"_{x_3 \cdot} & }
\end{align} where, again, the rightmost term is sitting in degree $0$.
\end{exmpl}

\section{Cohomology of Skew-Commutative Cubes}
\label{section-cohomology-of-skew-commutative-cubes}

In this section we study the cohomology of abstract complexes of the
form \eqref{exmpl-transform-of-reptwist-3d-case} with the aim to use
this later to calculate the cohomology sheaves of $\Psi(\mathcal{O}_0
\otimes \chi)$.

Let $X$ be a smooth separated scheme over $\C$.  Let $S :=
\{1,\dots,n\}$ and denote by $V$ the set of subsets of $S$. We
identify $V$ with the vertices of an ($n$-dimensional) unit cube. 

By a \em cube of line bundles \rm we mean the data of a line bundle
$\L_v$ on $X$ for every vertex $v \in V$ and a morphism $\alpha_v^i:
\L_{v \cup i} \rightarrow \L_v$ for each edge $v \cup \{i\}
\rightarrow v$ of the cube. One can think of it as a representation of
a cube-shaped quiver with vertices $V$ and arrows $v \cup \{i\}
\rightarrow v$ into a graded sheaf $\bigoplus_{v \in V} \mathcal{L}_v$
over $X$. 

We say that a cube is \em skew-commutative \rm (resp. \em
commutative\rm), if each of its two-dimensional faces forms an
anti-commutative (resp. \em commutative\rm) square, i.e.  $$\alpha^i_v
\circ \alpha^j_{v \cup i} + \; \alpha^j_v \circ \alpha^i_{v \cup j} =
0 \quad\text{ or }\quad \alpha^i_v \circ \alpha^j_{v \cup i} - \;
\alpha^j_v \circ \alpha^i_{v \cup j} = 0, $$ respectively, for all $v
\in V$ with $i,j \not\in v$. Given a skew-commutative cube, we can
form its total chain complex $T^\bullet$ where $T^{-i} =
\oplus_{|I|=i} \L_I$. The differential $d$ is given by summing the
maps in the cube. A commutative cube can be turned into a
skew-commutative cube by sprinkling some minus signs as follows
\begin{align} \label{eqn-commutative-to-skew-commutative}
\alpha^{i_j}_{\{i_1 \leq \dots \leq i_k\}} \rightsquigarrow (-1)^j
\alpha^{i_j}_{\{i_1 \leq \dots \leq i_k\}}.  \end{align}

We now turn to the case $n=3$ and further assume that the maps
$\alpha_v^i$ are non-zero. The vanishing locus of $\alpha_v^i$ is then
an effective Weil divisor on $X$ which we denote by $D_v^i$. Note that
$\mathcal{O}_X(D_v^i) \cong \L_v \otimes \L_{v \cup i}^{\vee}$.  The
corresponding skew-commutative cube looks like: \begin{align*}
\xymatrix{ & \mathcal{L}_{23} \ar"2,5"^<<{\alpha^3_{2}}
\ar"3,5"_<<{\alpha^2_{3}} & & & \mathcal{L}_1 \ar"2,6"^{\alpha^1} & \\
\mathcal{L}_{123} \ar"1,2"^{\alpha^1_{23}} \ar"2,2"^{\alpha^2_{13}}
\ar"3,2"_{\alpha^3_{12}} & \mathcal{L}_{13}
\ar"1,5"^>>>>{\alpha^3_{1}} \ar"3,5"_>>>>>{\alpha^1_{3}} & & &
\mathcal{L}_2 \ar"2,6"^{\alpha^2} & \mathcal{L} \\ & \mathcal{L}_{12}
\ar"1,5"^<<{\alpha^2_{1}} \ar"2,5"_<<{\alpha^1_{2}} & & &
\mathcal{L}_3 \ar"2,6"_{\alpha^3} & } \end{align*}

\begin{lemma}\label{lem:cube_cohomology} Let $T^\bullet$ be the total
complex of the skew-commutative cube $\{\mathcal{L}, \alpha\}$
depicted above. Then: \begin{enumerate} \item $H^0(T^\bullet) \cong \L
\otimes \O_Z$ where $Z$ is the scheme theoretic intersection $D^1 \cap
D^2 \cap D^3$ \item $H^{-1}(T^\bullet)$ has a three step filtration
with successive quotients \begin{itemize} \item $\O_Z \otimes
\L_{12}(\gcd(D_1^2,D_2^1))$ where $Z$ is the scheme theoretic
intersection of $\gcd(D_1^2,D_2^1)$ and the effective part of $D^3 +
\lcm(\tD_3^1,\tD_3^2) - \tD_1^2 - D^1$ \item $\O_Z \otimes
\L_{13}(\gcd(D_1^3,D_3^1))$ where $Z$ is the scheme theoretic
intersection of $\gcd(D_1^3,D_3^1)$ and the effective part of $D^2 +
\lcm(D_2^1,\tD_2^3) - \tD_3^1 - D^3$ \item $\O_Z \otimes
\L_{23}(\gcd(D_2^3,D_3^2))$ where $Z$ is the scheme theoretic
intersection of $\gcd(D_2^3,D_3^2)$ and the effective part of $D^1 +
\lcm(D_1^2,D_1^3) - \tD_2^3 - D^2$ \end{itemize} where $\tD^i_j =
D^i_j - \gcd(D^i_j,D^j_i)$ \item $H^{-2}(T^\bullet) \cong \L_{123}(D)
\otimes \O_D  $ where $D = \gcd(D_{23}^1, D_{13}^2, D_{12}^3)$ \item
$H^{-3}(T^\bullet) \cong 0$ \end{enumerate} \end{lemma} \begin{proof}
We shall need to work locally. Let $p$ be an arbitrary point of $X$.
We fix a local generator for each $\mathcal{L}_e$, identifying it
with $\mathcal{O}_{X,p}$. Each map $\alpha_e^i$ becomes then an
endomorphism of $\mathcal{O}_{X,p}$, i.e. the multiplication by some
$f_e^i \in \mathcal{O}_{X,p}$. Then $f_e^i$ is a local generator of
$D_e^i$, and we use it to identify $\mathcal{O}_X(-D_e^i)$ and
$\mathcal{O}_X(D_e^i)$ with $\mathcal{O}_{X,p}$.  

(1) We can view $H^0(T^\bullet)$ as the cokernel of $\bigoplus_i
\O_X(-D^i) \otimes \L \hookrightarrow \O_X \otimes \L$, which is by
definition $\O_Z \otimes \L$ where $Z = \cap_i D^i$. 

(2) Assume without loss of generality that $D^1,D^2,D^3$ have no prime
divisor in common. We claim that the kernel of $d: T^{-1} \rightarrow
T^0$ is generated by $\beta_1 = \frac{1}{h^{12}}(f^2,f^1,0), \beta_2 =
\frac{1}{h^{13}} (f^3,0,f^1)$ and $\beta_3 = \frac{1}{h^{23}}
(0,f^3,f^2)$ where $h^{ij} = \gcd(f^i,f^j)$. Suppose
$af^1+bf^2+cf^3=0$. We must have $\gcd( \frac{f^2}{h^{12}},
\frac{f^3}{h^{13}})$ divide $a$, for it divides $f_1 a$ and no prime
divisor of it can divide $f_1$ by definitions of $h^{12}$ and
$h^{13}$. Thus we can obtain $a$ as a linear combination of
$\frac{f^2}{h^{12}}$ and $\frac{f^3}{h^{13}}$, i.e. we can combine
$\beta_1$ and $\beta_2$ to get an element of the form $(a,b',c')$.
Then $(0,b'-b,c'-c)$ is a multiple of $\beta_3$ - for it must also lie
in the kernel of $d: T^{-1} \rightarrow T^0$, which means $f^2 (b'-b)
+ f^3 (c'-c) = 0$. This shows the claim.

For any $i \ne j$ denote by $g^{ij}$ the greatest common divisor of
$f^i_j$ and $f^j_i$ and let $\tf^i_j = f^i_j/g^{ij}$. By
skew-commutativity of the cube $f^j_i f^i + f^i_j f^j = 0$. Therefore
$\tf^i_j = f^i/h^{ij}$. We now rewrite $\beta_1 = (\tf^2_1, \tf^1_2,
0)$, $\beta_2 = (\tf^3_1, 0, \tf^1_3)$ and $\beta_3 =
(0,\tf^3_2,\tf^2_3)$. On the other hand, the image of $d$ in $T^1$ is
generated by $(f^2_1,f^1_2,0)$, $(f^3_1,0,f^1_3)$ and $(0,f^3_2,
f^2_3)$.  In particular, this means that if $D_i^j$ and $D_j^i$ have
no common divisors for any $i \ne j$ then $H^{-1}(T^\bullet) = 0$. 

Consider the filtration $\mbox{Im}(d) = F^0 \subset F^1 \subset F^2
\subset F^3 = \mbox{ker}(d)$ inside $T^{-1}$ where \begin{itemize}
\item $F^0$ is generated by $g^{12} \cdot (\tf^2_1, \tf^1_2,0)$,
$g^{13} \cdot (\tf^3_1, 0, \tf^1_3)$ and $g^{23} \cdot (0, \tf^3_2,
\tf^2_3)$ \item $F^1$ is generated by $g^{12} \cdot (\tf^2_1,
\tf^1_2,0)$, $g^{13} \cdot (\tf^3_1, 0, \tf^1_3)$ and $(0, \tf^3_2,
\tf^2_3)$ \item $F^2$ is generated by $g^{12} \cdot (\tf^2_1,
\tf^1_2,0)$, $(\tf^3_1, 0, \tf^1_3)$ and $(0, \tf^3_2, \tf^2_3)$ \item
$F^3$ is generated by $(\tf^2_1, \tf^1_2,0)$, $(\tf^3_1, 0, \tf^1_3)$
and $(0, \tf^3_2, \tf^2_3)$ \end{itemize} We show that globally the
quotient $F^2/F^1$ is isomorphic to $\mathcal{O}_Z \otimes
\L_{13}(\gcd(D_1^3,D_3^1))$, where $Z$ is the intersection of
$\gcd(D_1^3,D_3^1)$ and the effective part of $D^2 +
\lcm(D_2^1,\tD_2^3) - \tD_3^1 - D^3$ (the other two quotients are
computed similarly). 

We can combine $g^{12} \cdot (\tf_1^2,\tf_2^1,0) = (f_1^2, f_2^1,0)$
and $(0,\tf_2^3,\tf_3^2)$ to give a multiple of $(\tf_1^3,0,\tf_3^1)$.
The smallest such multiple is $$\frac{\tf_2^3}{\gcd(f_2^1,\tf_2^3)}
(f_1^2,f_2^1,0) - \frac{f_2^1}{\gcd(f_2^1,\tf_2^3)}
(0,\tf_2^3,\tf_3^2) = \frac{\tf_3^2 \cdot f_2^1}{\tf_3^1 \cdot
\gcd(f_2^1, \tf_2^3)} (\tf_1^3,0,-\tf_3^1).$$  Let $h$ be the regular
part of $\frac{\tf_3^2 \cdot f_2^1}{\tf_3^1 \cdot \gcd(f_2^1,
\tf_2^3)}$, it is well defined since $\mathcal{O}_{X,p}$ is UFD.
Then, locally, $F^2/F^1$ is $\mathcal{O}_{X,p}/(g^{13},h) \cdot
(\tf_1^3,0,\tf_3^1)$. Consider now the global injection $$
\L_{13}(\gcd(D_1^3,D_3^1)) \xrightarrow{\alpha^3_1 \oplus 0 \oplus
\alpha^1_3} \L_1 \oplus \L_2 \oplus \L_3.$$ Locally, its image in
$\L_1 \oplus \L_2 \oplus \L_3$ is just $\mathcal{O}_{X,p} \cdot
(\tf_1^3,0,\tf_3^1)$. It therefore only remains to show that the ideal
$(g^{13},h)$ defines $Z$ in $\mathcal{O}_{X,p}$. This follows since
$(g^{13})$ defines $\gcd(D^1_3,D^3_1)$ and $(h)$ -- the effective part
of \begin{eqnarray*} \tD_3^2 + D_2^1 - \tD_3^1 - \gcd(D_2^1, \tD_2^3)
&=& \tD_2^3 - D^3 + D^2 + D_2^1 - \gcd(D_2^1, \tD_2^3) - \tD_3^1 \\
&=& D^2 + \lcm(D_2^1, \tD_2^3) - D^3 - \tD_3^1, \end{eqnarray*} where
the first equality follows since $\tD_3^2 + D^3 = \tD_2^3 + D^2$ by
the skew-commutativity. 

(3) Suppose, at first, that $D_{23}^1, D_{13}^2$ and $D_{12}^3$ have
no divisor in common. We claim that, locally, if $(a,b,c) \in T^{-2}$
lies in the kernel of $d$ then it lies in $d(T^{-3})$. Let $g$ be a
prime factor of $f_{23}^1$ and suppose (without loss of generality)
that $g$ does not divide $f_{12}^3$. Let $p,q$ be the largest
non-negative integers such that $g^p|f_{23}^1$ and $g^q|f_2^3$.
Commutativity of faces implies $f_{23}^1 \cdot f_2^3 = f_{12}^3 \cdot
f_2^1$ so $g^{p+q}$ divides $f_2^1$. But since $(a,b,c)$ lies in the
kernel $af_2^3-cf_2^1=0$. This means $g^{p+q}$ divides $af_2^3$. Since
the local ring is a UFD we get $g^p|a$. Since this is true for any $g$
we get $f_{23}^1|a$. A similar argument shows $f_{13}^2|b$ and
$f_{12}^3|c$ so $(a,b,c)$ lies in the image of $d$.  

More generally, if $D$ is the largest common divisor, the we can
factor $d: T^{-3} \rightarrow T^{-2}$ as $\L_{123} \rightarrow
\L_{123}(D)$ followed by a map which has no common divisors. By the
above, the image of $\L_{123}(D)$ under this map equals the kernel of
$d:T^{-2} \rightarrow T^{-1}$. Thus $H^{-2}(T^\bullet)$ is precisely
the cokernel of $\L_{123} \hookrightarrow \L_{123}(D)$ which is $\O_D
\otimes \L_{123}(D)$. 
 
(4) Since any non-zero map between line bundles is injective
$\L_{123}$ injects into (say) $\L_{12}$ so $H^{-3}(T^\bullet) \cong
0$.  \end{proof} 

Lemma \ref{lem:cube_cohomology} makes it clear that given a
$\gnat$-family $\mathcal{F} = \bigoplus \mathcal{L}(-D_\chi)$ in order
to calculate the cohomology sheaves of $\Psi_\mathcal{F}(\mathcal{O}_0
\otimes \chi)$ we need to calculate the vanishing loci $D^i_v$ of all
maps $\alpha^i_v$ in   all the corresponding skew-commutative cubes.
By Proposition \ref{prps-preimages-of-point-reps-abelian-case} this
amounts to calculating the vanishing loci of maps \begin{align*}
\mathcal{L}(- D_\chi) &\rightarrow \mathcal{L}(- D_{\chi
\kappa(x_i)^{-1}}) \\ s &\mapsto x_i \cdot s \end{align*} for all
$\chi$ and $x_i$. We learn to do this in our next section. 
 
\section{The Associated Representation of the McKay Quiver}
\label{section-the-mckay-quiver}

In this section, we take $G$ to be a finite abelian subgroup of
$\gsl_3(\mathbb{C})$ such that $\mathbb{C}^3/G$ has a single isolated
singularity at the origin. Let $Y \rightarrow \C^3/G$ be any crepant
resolution and $\mathcal{F} = \bigoplus \mathcal{L}(-D_\chi)$ be
any $\gnat$-family on $Y$. We fix an exceptional divisor $E \in
\Except(Y)$ and generically along $E$ study the behavior of the McKay
quiver representation $\mckquiv_{\mathcal{F}}$. 

\subsection{The McKay quiver of $G \subset \gsl_3(\mathbb{C})$ and its
planar embedding}

\label{section-the-mckay-quiver-of-G-and-its-planar-embedding}

By a \it quiver \rm we mean a vertex set $Q_0$, an arrow set $Q_1$ and
a pair of maps $h \colon Q_1 \rightarrow Q_0$ and $t \colon Q_1
\rightarrow Q_0$ giving the head $hq \in Q_0$ and the tail $tq \in
Q_0$ of each arrow $q \in Q_1$. 

\begin{defn} The \it McKay quiver $\mckquiv$ \rm of (a not necessarily 
abelian) $G \subset \gsl_3(\C)$
is the quiver whose vertex set $Q_0$ are the irreducible
representations $\rho$ of $G$ and whose arrow set $Q_1$ has $\dim
\homm_G (\rho_i \otimes \givrep^\vee , \rho_j)$ arrows going from the
vertex $\rho_i$ to the vertex $\rho_j$.  
\end{defn}

As we assumed $G$ to be abelian 
$\givrep^\vee \cong \bigoplus \mathbb{C} x_i$ as $G$-representations.
Let $\chi_i$ and $\chi_j$ be any pair of characters of $G$.
We then see that $\homm_G (\chi_i \otimes \givrep^\vee,
\chi_j)$ contains a copy of $\C$ for each $k \in \{1,2,3\}$ satisfying
$\chi_j = \kappa(x_k) \chi_i$. Thus each vertex $\chi \in G^\vee$ of $\mckquiv$
has $3$ arrows emerging from it and going to vertices $\kappa(x_k)
\chi$ for $k = 1,2,3$.  We denote the arrow from $\chi$ to
$\kappa(x_k) \chi$ by $(\chi, k)$ and say further that it is an
$x_k$-arrow. For an example see Section
\ref{section-example-the-group}.
 
The fact that $G \subseteq \gsl_3(\mathbb{C})$ allows us to embed its
McKay quiver in a two dimensional real torus. We briefly recall this
embedding as constructed by Craw and Ishii in \cite{Craw-Ishii-02}.

Consider the maximal torus $(\mathbb{C}^*)^3 \subset
\gsl_3(\mathbb{C})$ containing $G$. We have an exact sequence of
abelian groups: \begin{align}\label{seq-g-torus} \xymatrix{ 0 \ar[r] &
G \ar[r] & (\mathbb{C}^*)^3 \ar[r] & T \ar[r] & 0 } \end{align} where
$T$ is the quotient torus which acts on the quotient space
$\mathbb{C}^3/G$. By applying $\homm(\bullet, \mathbb{C}^*)$ to
$\eqref{seq-g-torus}$ we obtain an exact sequence \begin{align}
\label{seq-monom-char} \xymatrix{ 0 \ar[r] & M \ar[r] & \mathbb{Z}^3
\ar[r]^{\rho} & G^\vee \ar[r] & 0 } \end{align} where $\mathbb{Z}^3$
is thought of as the lattice of exponents of Laurent monomials. Given
$m = (k_1, k_2, k_3) \in \mathbb{Z}^3$ we write $x^{m}$ for $x_1^{k_1}
x_2^{k_2} x_3^{k_3}$. Note that $\rho$ is the weight map, that is
$x^{m}(g\cdot\bv) = \rho(m)(g)\; x^{m}(\bv)$ for any $\bv \in
\mathbb{C}^3$. By definition $g\cdot x^m(\bv) = x^m(g^{-1}\cdot\bv)$,
so $\rho = \kappa^{-1}$, where $\kappa$ is the map introduced in
Section $\ref{section-inverse-transform-abelian-case}$, which maps
every Laurent monomial to the character $G$ acts on it with. $M$ is
then the sublattice in $\mathbb{Z}^3$ of (exponents of) $G$-invariant
Laurent monomials. As $G \subseteq \gsl_3(\mathbb{C)}$ we have
$(1,1,1) \in M$, i.e. $x_1x_2x_3$ is an invariant monomial.  Take $H =
\mathbb{Z}^3 / \mathbb{Z}(1,1,1)$ and $M' = M/\mathbb{Z}(1,1,1)$. Then
$\eqref{seq-monom-char}$ induces \begin{align} 0 \rightarrow M'
\rightarrow H \rightarrow G^\vee \rightarrow 0 \end{align} For every
Laurent monomial $x^{m}$ for some $m \in \mathbb{Z}^3$ we shall write
$[x^m]$ for the class of $m$ in $H$, e.g.  $[x_1 x_2^2]$ for the class
of $(1,2,0)$.   

\begin{defn} The universal cover $U$ of $\mckquiv$ is the quiver whose
vertex set are the elements of $H$ and whose arrow set is $$\{(h, h +
[x_i]) \;|\; h \in H, i \in 1,2,3\}$$ \end{defn} We have a natural
`embedding' of $U$ into $H \otimes \mathbb{R} \cong \mathbb{R}^2$,
where the arrow $(h,h + [x_i])$ is identified with the line segment
$\{h + \lambda [x_i] \;|\; \lambda \in (0,1) \}$.  For illustrative
purposes, we fix a specific isomorphism:  \begin{align*} \phi_H:\; H
\otimes \mathbb{R} \rightarrow \mathbb{R}^2:\quad \begin{cases} [x_1]
\mapsto (\frac{\sqrt{3}}{2}, -\frac{1}{2}) \\ [x_3] \mapsto (0, 1)
\end{cases} \end{align*}

\begin{figure}[h] \begin{center}
\includegraphics[scale=0.15]{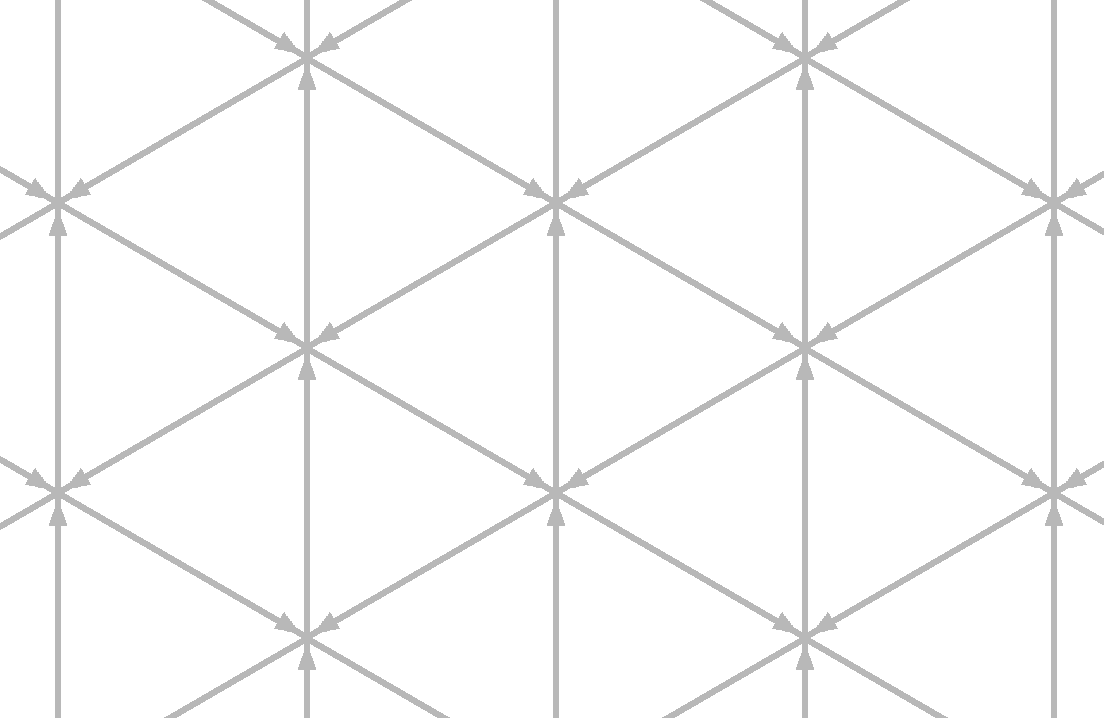} \end{center}
\caption{\label{figure-1} Tessellation of $\mathbb{R}^2$ by $U$}
\end{figure}

As shown in Figure \ref{figure-1} the image of $U$ under $\phi_H$ is a
tessellation of $\mathbb{R}^2$ into regular triangles. The quotient
$T_G = H \otimes \mathbb{R} /M'$ is a two dimensional real torus.
$\mckquiv$ is the quotient of $U$ by the action of $M'$ and thus is a
tessellation of $T_G$ into regular triangles. The boundary of a
triangle is the cycle $$\chi \xrightarrow{x_{\sigma(1)}}
\kappa(x_{\sigma(1)}) \chi \xrightarrow{x_{\sigma(2)}}
\kappa(x_{\sigma(1)}x_{\sigma(2)}) \chi \xrightarrow{x_{\sigma(3)}}
\kappa(x_{\sigma(1)}x_{\sigma(2)}x_{\sigma(3)}) \chi = \chi$$ where
$\sigma \in S_3$ is a permutation.  There are altogether $2|G|$
triangles in $\mckquiv$. When depicting $T_G$ in diagrams we draw its
fundamental domain in $\mathbb{R}^2$. The depiction of the McKay
quiver of $G = 1/13(1,5,7)$, as embedded into $T_G$, is illustrated in
Figure \ref{figure-2} (Section \ref{section-example-the-group}).

\subsection{The associated representation $\mckquiv_\mathcal{F}$ of
the McKay quiver}

\label{section-associated-representation}

A \it representation of a quiver \rm over $\spec \mathbb{C}$ is a
graded vector space $\bigoplus_{i \in Q_0} V_i$ together with a
collection of linear maps $\{\alpha_q\colon V_{tq} \rightarrow
V_{hq}\}_{q \in Q_1}$.  Over an arbitrary scheme $S$ this translates
into a graded locally free sheaf $\bigoplus_{i \in Q_0}
\mathcal{N}_{i}$ and a collection of morphisms $\{\alpha_q\colon
\mathcal{N}_{tq} \rightarrow \mathcal{N}_{hq}\}_{q \in Q_1}$.

\begin{defn} Given a $\gnat$-family $\mathcal{F} = \bigoplus_{\chi \in
G^\vee} \mathcal{L}(-D_\chi)$ on $Y$ \em the associated representation
$\mckquiv_\mathcal{F}$ of the McKay quiver over $Y$ is the sheaf
$\mathcal{F}$ where the summand $\mathcal{L}(-D_{\chi^{-1}})$ is
graded by $\chi$ (since $G$ acts by $\chi$ on
$\mathcal{L}(-D_{\chi^{-1}})$). The maps $\alpha_{\chi, x_k}$ are
defined by \begin{align} \label{eqn-gcon-rep-definition} \alpha_{\chi,
x_k}:\; \mathcal{L}(-D_{\chi^{-1}}) \rightarrow
\mathcal{L}(-D_{\kappa(x_k)^{-1} \chi^{-1}}), \; s \mapsto x_k \cdot
s.  \end{align} where $x_k$ act via $\twalg$-module structure of
$\mathcal{F}$.  \end{defn}

Denote by $B_{\chi,x_k}$ the effective Weil divisor on $Y$ where
$\alpha_{\chi, x_k}$ vanishes. As shown in Section 4.6 of
\cite{Logvinenko-DerivedMcKayCorrespondenceViaPureSheafTransforms},
$B_{\chi,x_k}$ is given by the formula \begin{align}
\label{eqn-zero-divisors-formula} B_{\chi,x_k} = D_{\chi^{-1}} + (x_i)
- D_{\kappa(x_k)^{-1} \chi^{-1}} \end{align}

\begin{prps}
\label{prps-zero-divisors-are-supported-on-except-set-plus-hyperplane}
Each $B_{\chi,x_k}$ is supported on the union of the exceptional locus
$\Except(Y)$ of $Y \rightarrow \C^3/G$ and the strict transform of the
surface $x_k^{|G|}=0$ in $\C^3/G$.  \end{prps} \begin{proof} Observe
that it follows from the formula \ref{eqn-zero-divisors-formula} that
\begin{align} \label{eqn-principal-divisor-of-x_i-via-zero-divisors}
\sum_{\chi \in G^\vee} B_{\chi,x_k} = \sum_{\chi \in G^\vee}
D_{\chi^{-1}} + |G|(x_k) - \sum_{\chi \in G^\vee} D_{\chi^{-1}} =
|G|(x_k) = (x_k^{|G|}).  \end{align} As each $B_{\chi,x_k}$ is
effective, its support is contained in that of the principal Weil
divisor $(x_k^{|G|})$ on $Y$. This is precisely the strict transform
of the surface $(x_k^{|G|})$ in $\C^3/G$ union the exceptional locus.
\end{proof}

The following simple result, which generalizes Corollary 10.2 of
\cite{Craw-Ishii-02}, shows that $G$ being a subgroup of
$\gsl_3(\mathbb{C})$ together with $Y \rightarrow \C^3/G$ being
crepant impose a strong restriction on the multiplicities of the
exceptional divisors in  $B_{\chi,i}$.

\begin{prps} \label{prps-just-one-out-of-three} Let $E$ be an
irreducible exceptional divisor of $Y$ and $\mckquiv_\mathcal{F}$ a
$\gnat$-family. In any triangle of $\mckquiv$, one of the three maps
corresponding to $\mckquiv_\mathcal{F}$ vanishes along $E$ with
multiplicity one while the other two do not vanish along $E$.
Consequently, $E$ occurs in any $B_{\chi,i}$ with multiplicity zero or
one.  \end{prps} \begin{proof} Suppose the vertices of the triangle
are $\chi$, $\kappa(x_1) \chi$ and $\kappa(x_1 x_2) \chi$. Then by
formula \ref{eqn-principal-divisor-of-x_i-via-zero-divisors}
\begin{align} B_{\chi, x_1} + B_{\kappa(x_1) \chi, x_2} +
B_{\kappa(x_1 x_2) \chi, x_3} = (x_1x_2x_3).  \end{align} This is a
$G$-invariant monomial since $G \subseteq \gsl_3(\C)$. 

Since any $B_{\chi,x_i}$ is effective it suffices to show that
$v_E(x_1 x_2 x_3) = 1$ where $v_E$ is the valuation at $E$. This
is in fact true for any crepant (and thus monomial) valuation
of $K(Y)$ - see \cite{ItReid96}, in particular Step 5 of the proof
of Theorem 1.4.
\end{proof}

\begin{figure}[h] \begin{center}
\includegraphics[scale=0.18]{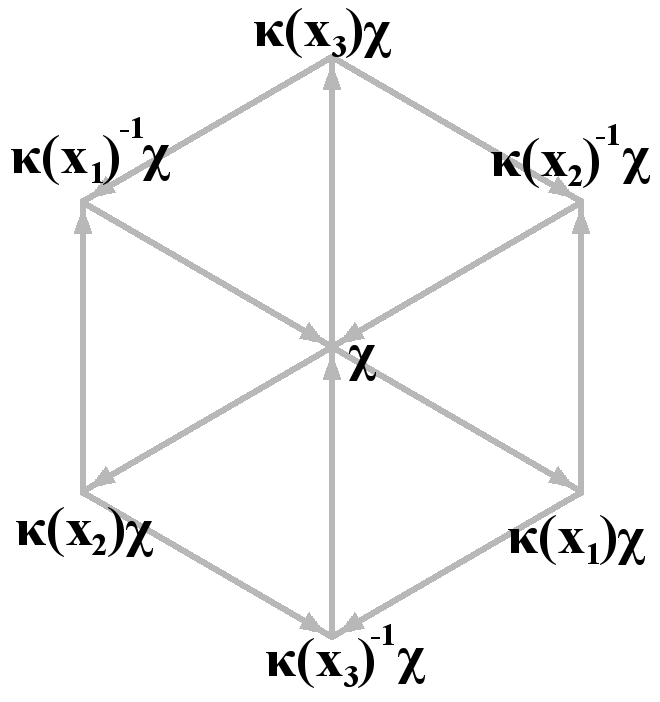} \end{center}
\caption{\label{figure-04} The subquiver $\hex(\chi)$.} \end{figure}

The relevance of the associated representation to the context of this
paper is as follows. Let $\chi$ be a character of $G$. Recall the
definition of the $3$-dimensional cube quiver in Section
\ref{section-cohomology-of-skew-commutative-cubes}. Define 
$\hex(\chi)$ to be the subquiver of the McKay quiver which contains
the six triangles which share the vertex $\chi$. This is depicted in
Figure \ref{figure-04}.  
The restriction of the representation $\mckquiv_{\mathcal{F}}$ to
$\hex(\chi)$ gives a representation $\hex(\chi)_\mathcal{F}$. In the
language of Section \ref{section-cohomology-of-skew-commutative-cubes}
this is a commutative cube of line bundles (which has been projected
onto the plane). We can turn it into a skew-commutative cube as per
\eqref{eqn-commutative-to-skew-commutative}. The following is then an
immediate consequence of Proposition
\ref{prps-preimages-of-point-reps-abelian-case}.

\begin{prps} \label{prps-image-of-chi-twisted-O0-is-hex-chi-inverse}
Let $\chi$ be a character of $G$, $Y \rightarrow \C^3/G$ a crepant
resolution and $\mathcal{F}$ a $\gnat$-family on $Y$. Then
$\Psi_{\mathcal{F}}(\O_0 \otimes \chi) \in D(Y)$ is the total complex
of the skew-commutative cube of line bundles induced by
$\hex(\chi^{-1})_\mathcal{F}$.  \end{prps}

\subsection{Behavior of $\mckquiv_\mathcal{F}$ along $E$}
 
Given a $\gnat$-family $\mathcal{F}$ we would like to understand better 
the divisors $B_{\chi, x_k}$ where 
the maps $\alpha_{\chi, x_k}$ of $\mckquiv_\mathcal{F}$ vanish. In 
this section we fix one irreducible exceptional divisor $E$ and 
study the behaviour of $\mckquiv_{\mathcal{F}}$ and its
subrepresentations $\hex(\chi)_{\mathcal{F}}$ generically along $E$.

\begin{prps} \label{prps-painting-of-hex-chi} Fix a character $\chi
\in G^\vee$. On the quiver $\hex(\chi)$ we color an arrow black if the
zero divisor of the corresponding map $\alpha$ in
$\hex(\chi)_\mathcal{F}$ contains $E$, otherwise we color it grey.
Then all the possible colorings of $\hex(\chi)$ are depicted in
Figures \ref{figure-05}-\ref{figure-08}.  
\begin{figure}[h]
\begin{center} \includegraphics[scale=0.10]{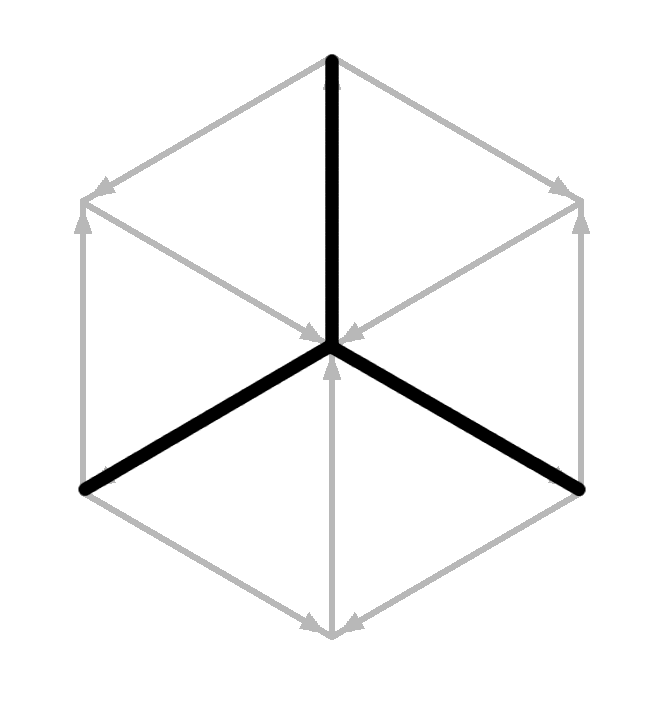}
\includegraphics[scale=0.10]{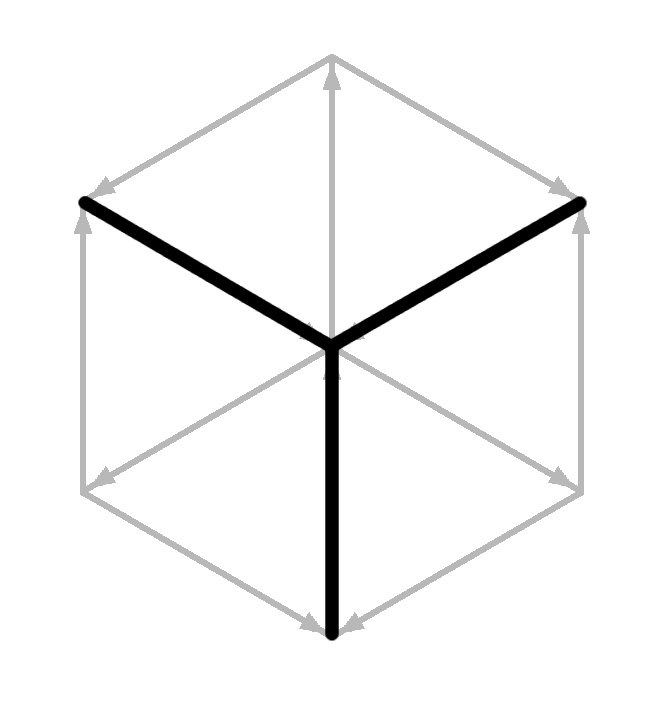}
\caption{\label{figure-05} The $(3,0)$-sink and the $(0,3)$-sink.}
\end{center} \begin{center}
\includegraphics[scale=0.10]{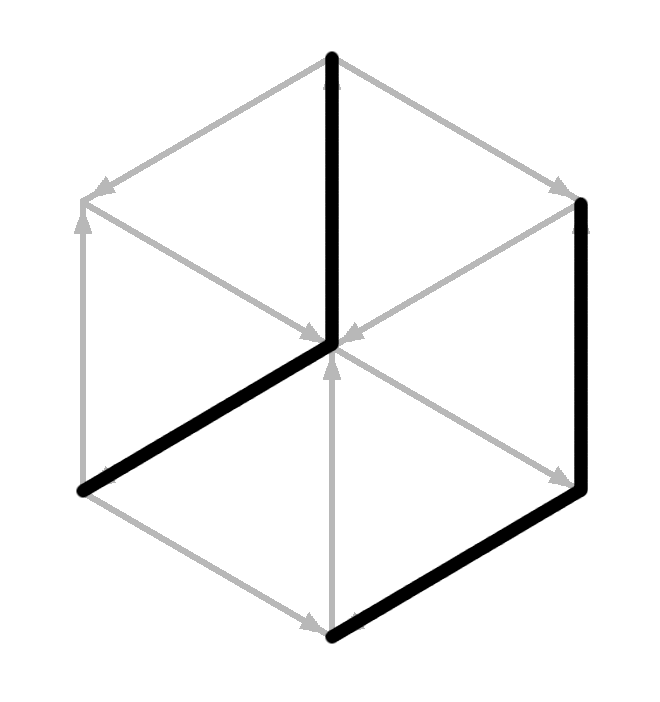}
\includegraphics[scale=0.10]{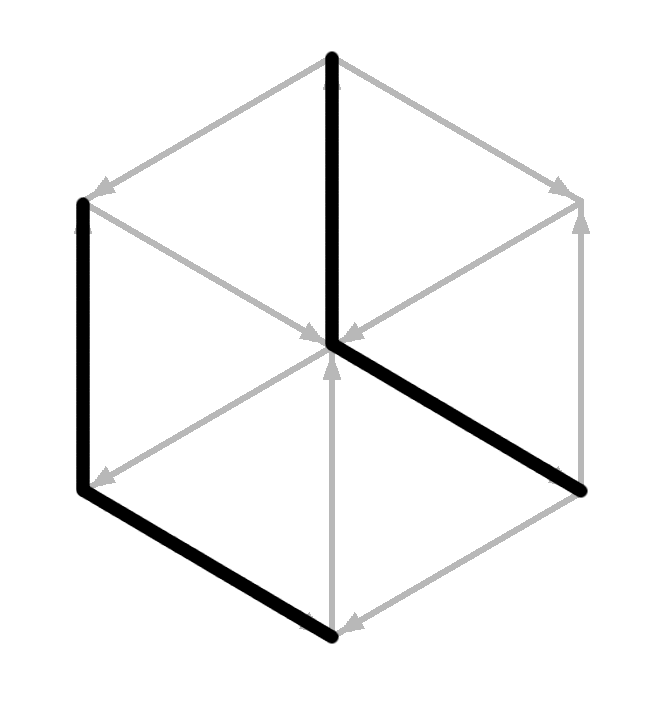}
\includegraphics[scale=0.10]{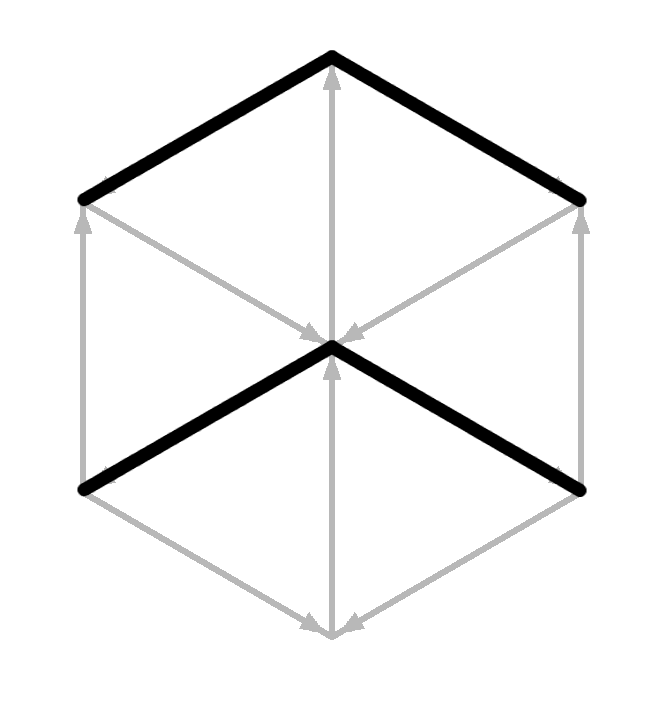}
\includegraphics[scale=0.10]{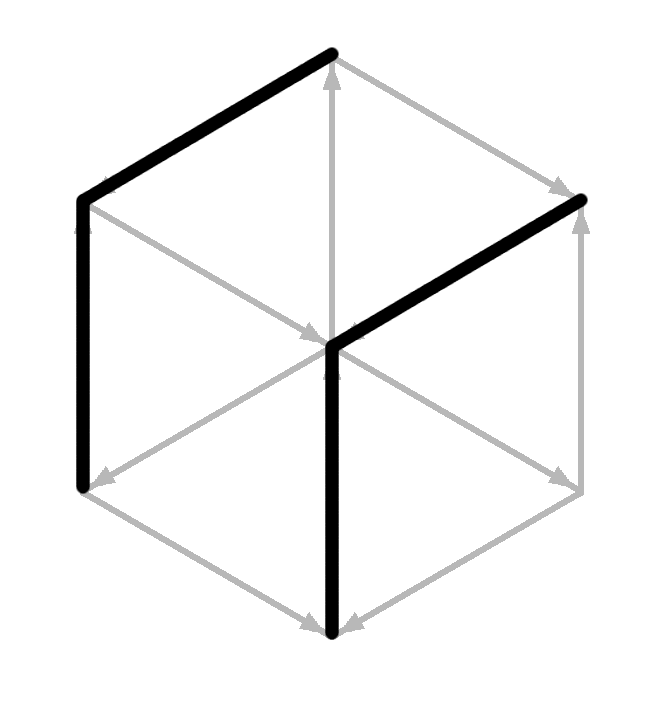}
\includegraphics[scale=0.10]{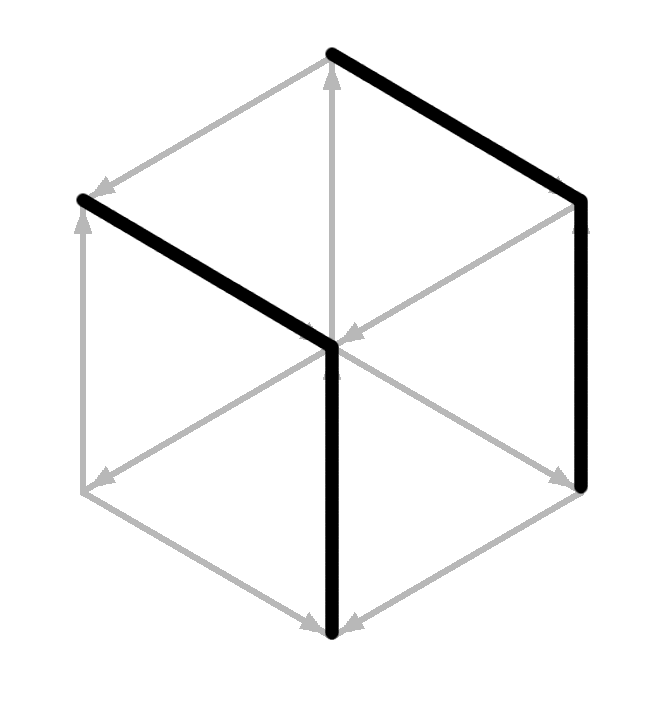}
\includegraphics[scale=0.10]{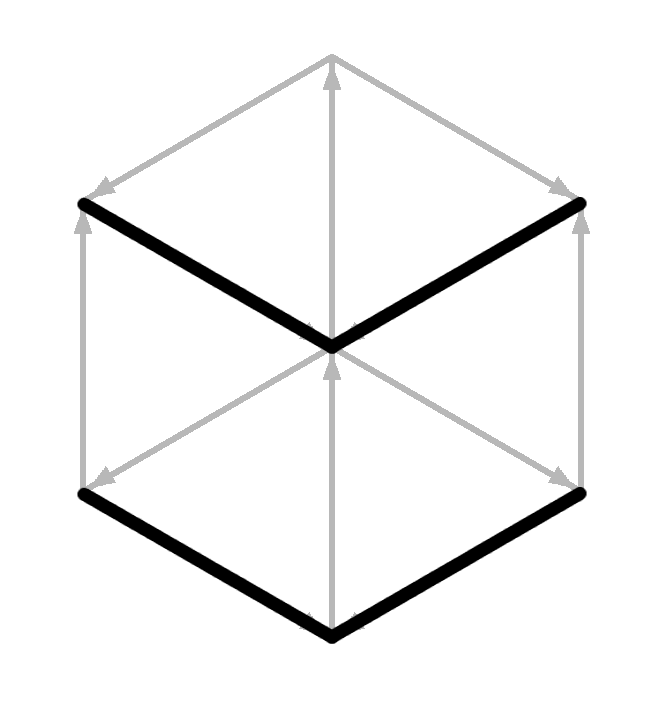}
\caption{\label{figure-06} The $x_1$-, $x_2$- and $x_3$-
$(1,0)$-charges, the $x_1$-, $x_2$- and $x_3$- $(0,1)$-charges}
\end{center} \begin{center}
\includegraphics[scale=0.10]{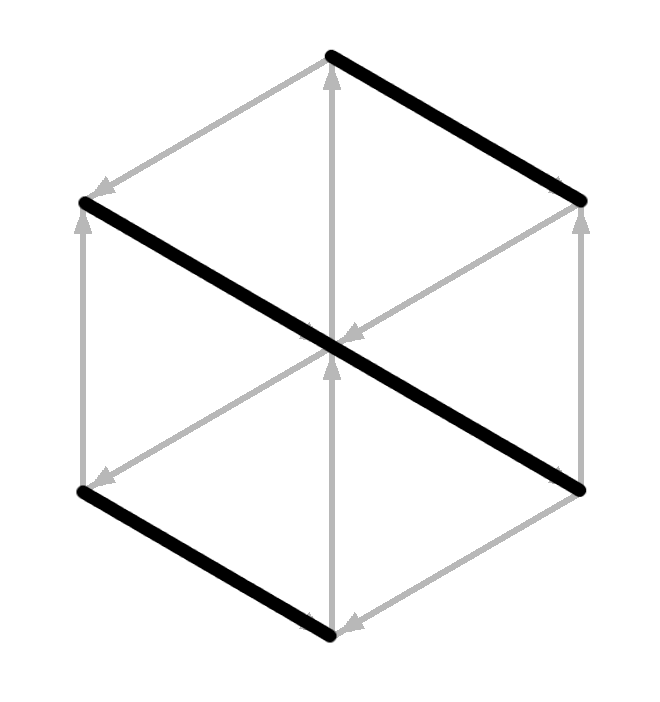}
\includegraphics[scale=0.10]{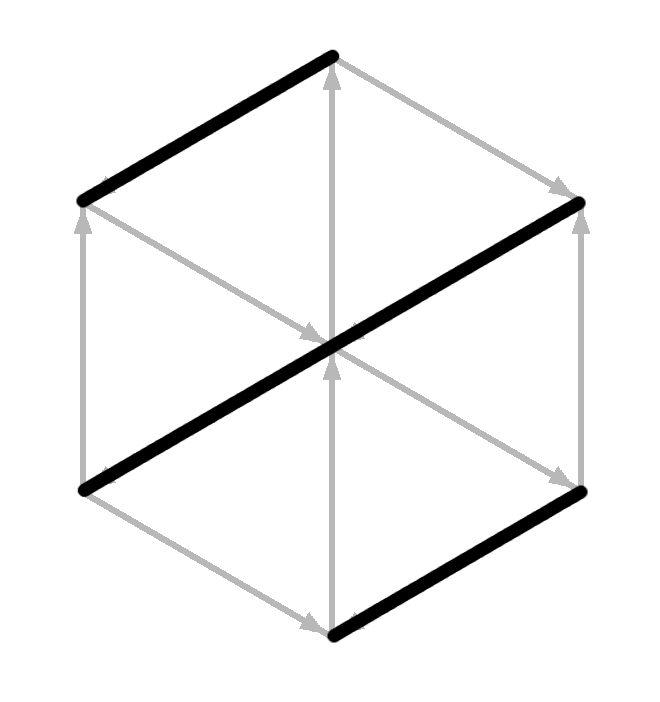}
\includegraphics[scale=0.10]{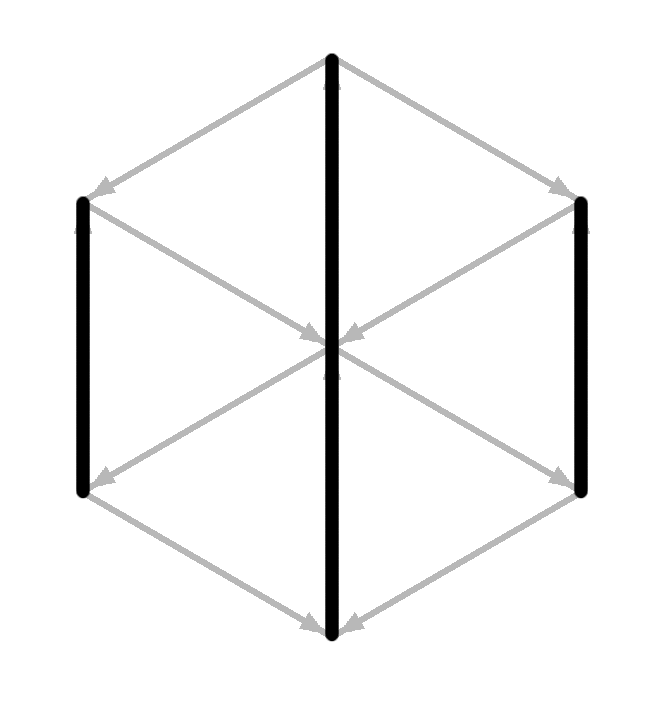}
\caption{\label{figure-07} The $x_1$-tile, the $x_2$-tile and the
$x_3$-tile.} \end{center} \begin{center}
\includegraphics[scale=0.10]{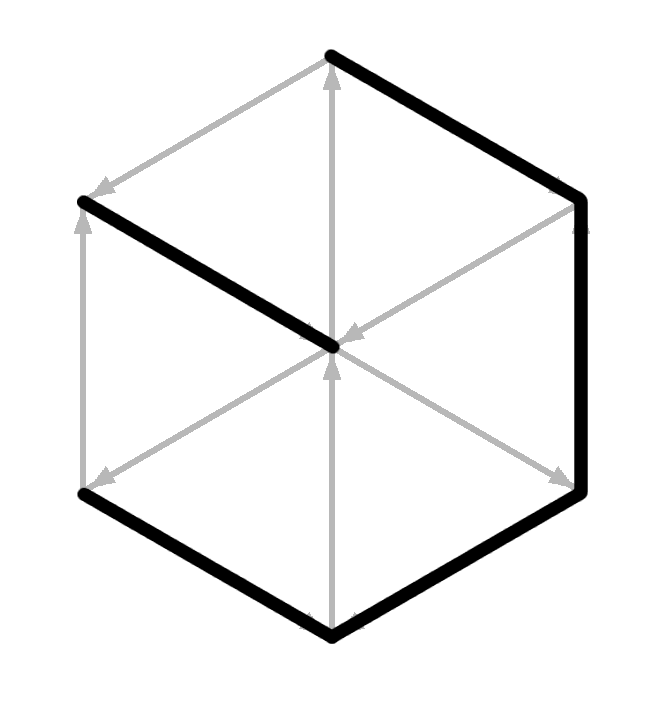}
\includegraphics[scale=0.10]{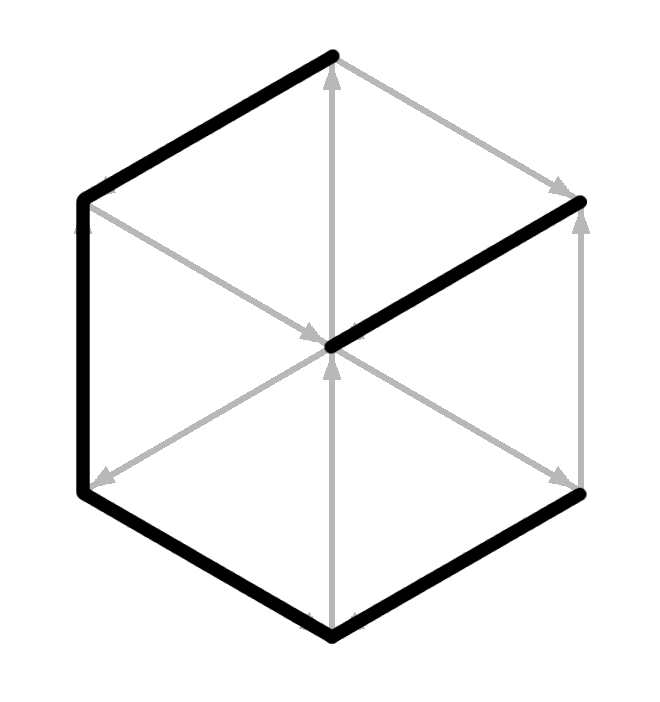}
\includegraphics[scale=0.10]{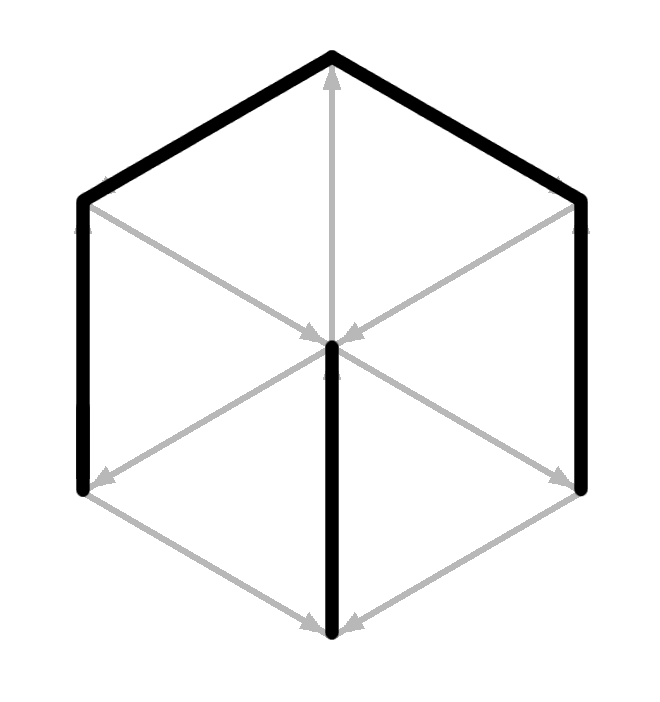}
\includegraphics[scale=0.10]{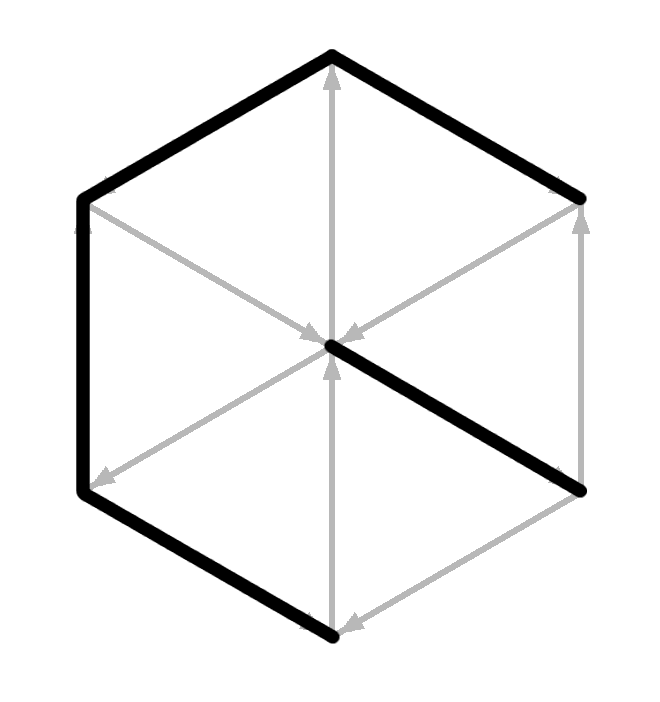}
\includegraphics[scale=0.10]{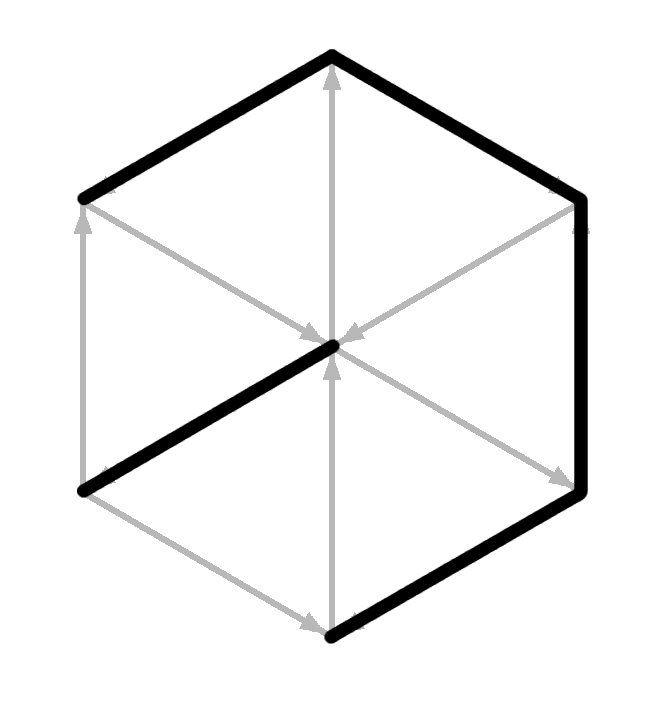}
\includegraphics[scale=0.10]{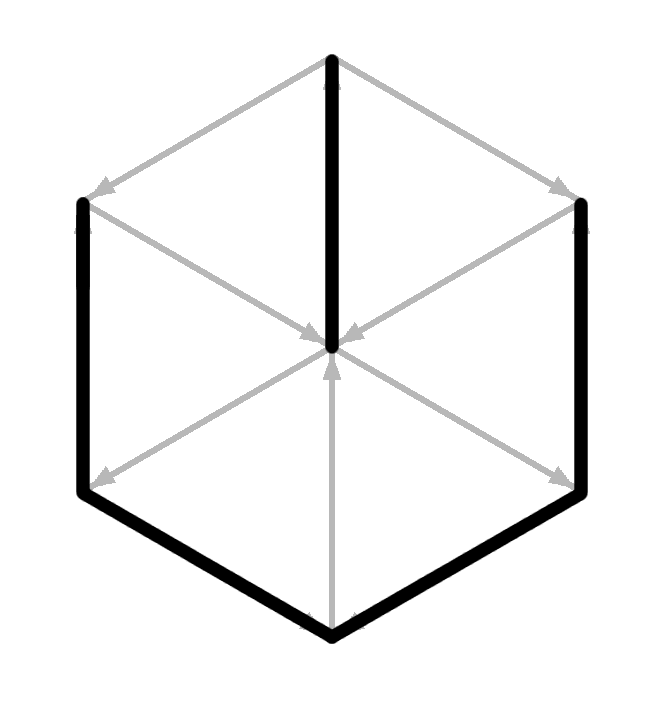}
\caption{\label{figure-08} The $x_1$-, $x_2$- and $x_3$-
$(1,2)$-sources, the $x_1$-, $x_2$- and $x_3$- $(2,1)$-sources.}
\end{center} \begin{center}
\includegraphics[scale=0.10]{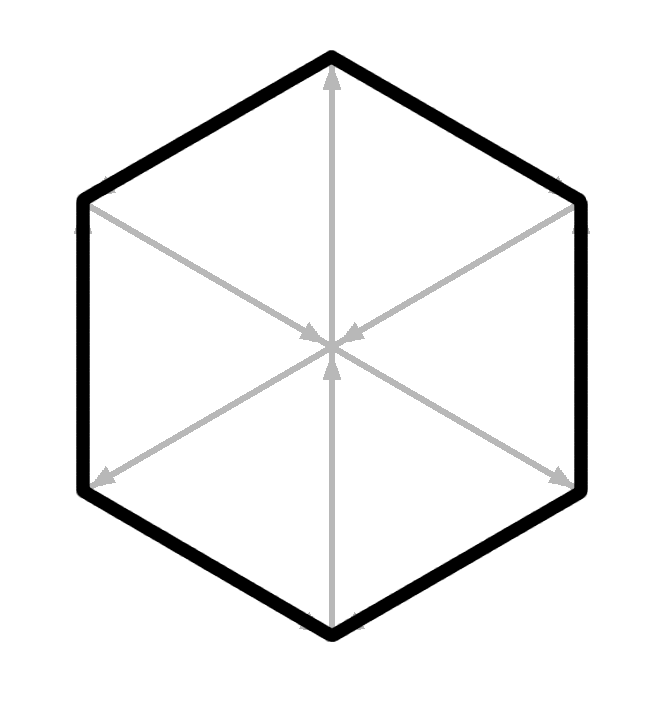}
\caption{\label{figure-09} The $(3,3)$-source} \end{center}
\end{figure} 
\end{prps} 
\begin{proof} 
It follows from Proposition $\ref{prps-just-one-out-of-three}$ that in
each of the six triangles making up $\hex(\chi)$ exactly one out of
three arrows is colored black. It is a simple combinatorial exercise
to check that all the possible configurations satisfying this
condition are depicted in the Figures \ref{figure-05}-\ref{figure-08}.  
\end{proof}

The terminology used in Figures \ref{figure-05}-\ref{figure-08} will
be justified shortly as a consequence of the following combinatorial
observations. 

\begin{cor} \label{cor-charges} Let $\chi$ be a character of $G$.
\begin{itemize} \item If $\chi$ is an $x_i$-$(1,0)$-charge, then
$\kappa(x_i)\chi$ is either another $x_i$-$(1,0)$-charge or a
$(3,0)$-sink. While $\kappa(x_i)^{-1}\chi$ is either another
$x_i$-$(1,0)$-charge, an $(3,3)$-source, an $x_i$-$(1,2)$-source or an
$x_j$-$(2,1)$-source, where $j \neq i$.  \item If $\chi$ is an
$x_i$-$(0,1)$-charge, then $\kappa(x_i)^{-1}\chi$ is either another
$x_i$-$(0,1)$-charge or a $(0,3)$-sink, while $\kappa(x_i)\chi$ is
either another $x_i$-$(0,1)$-charge, an $(3,3)$-source, an
$x_i$-$(2,1)$-source or an $x_j$-$(1,2)$-source, where $j \ne i$.
\end{itemize} \end{cor}

As $\kappa(x_i) \chi$ is the next vertex after $\chi$ in the direction
of $x_i$-arrows and $\kappa(x_i)^{-1}\chi$ is the previous one,
Corollary \ref{cor-charges} may be interpreted as saying that charges
propagate from a source to a sink in a straight line. The
$x_i$-$(1,0)$-charges propagate in the direction of $x_i$-arrows,
while $x_i$-$(0,1)$-charges propagate against the direction of
$x_i$-arrows. We say that $x_i$ is the \em orientation \rm of
$x_i$-$(1,0)$-charge and that $(1,0)$ is its \em type \rm (and the
same for the other types of vertices).  

\begin{cor} \label{cor-sources} Let $\chi$ be a character of $G$.
\begin{itemize} \item If $\chi$ is a $(3,3)$-source, then for $i \in
\{ 1,2,3 \}$ the vertex $\kappa(x_i) \chi$ is either an
$x_i$-$(1,0)$-charge or a $(3,0)$-sink, while the vertex
$\kappa(x_i)^{-1} \chi$ is either an $x_i$-$(0,1)$-charge or a
$(0,1)$-sink.  \item If $\chi$ is an $x_i$-$(1,2)$-source, then the
vertex $\kappa(x_i) x_i$ is either an $x_i$-$(1,0)$-charge or a
$(3,0)$-sink, while each of the two vertices $\kappa(x_j)\chi$, for $j
\neq i$, is either $x_j$-$(0,3)$-charge or a $(0,3)$-sink.  \item If
$\chi$ is an $x_i$-$(2,1)$-source, then the vertex $\kappa(x_i)^{-1}
x_i$ is either an $x_i$-$(0,1)$-charge or a $(0,3)$-sink, while each
of the two vertices $\kappa(x_j)\chi$, for $i \neq j$, is either
$x_j$-$(3,0)$-charge or a $(3,0)$-sink.  \end{itemize} \end{cor}

\begin{cor} \label{cor-sinks} Let $\chi$ be a character of $G$.
\begin{itemize} \item If $\chi$ is a $(3,0)$-sink, then each of the
vertices $\kappa(x_i)^{-1}\chi$ for $i=1,2,3$ is either an
$x_i$-$(1,0)$-charge, $(3,3)$-source, $x_i$-$(1,2)$-source or an
$x_j$-$(2,1)$-source, where $j \neq i$.	\item If $\chi$ is an
$(0,3)$-sink, then each of the vertices $\kappa(x_i) \chi$ for
$i=1,2,3$ is either an $x_i$-$(0,1)$-charge, $(3,3)$-source,
$x_i$-$(2,1)$-source or an $x_j$-$(1,2)$-source, where $j \neq i$.
\end{itemize} \end{cor}

In other words, sources emit charges which propagate in a straight
line to a sink. An $(a,b)$-source emits $a$ charges of type $(1,0)$
(which propagate in the direction of the arrows of the quiver) and $b$
charges of type $(0,1)$ (which propagate against the direction of the
arrows). Similarly, an $(a,b)$-sink receives $a$ charges of type
$(1,0)$ and $b$ charges of type $(0,1)$.

Corollaries \ref{cor-charges}-\ref{cor-sinks} imply that given a
\gnat-family $\sF$ and an exceptional divisor $E$ we can draw a graph
on top of the McKay quiver $\mckquiv$ as follows. The vertices are the
sources and sinks (with respect to the coloring of $\mckquiv$ keeping
track of the vanishing of $E$) while the edges are the straight lines
along which charges propagate from a source to a sink. We call it the
\em sink-source graph \rm and denote it $\sinksource_{\mathcal{F},E}$. 

The sink-source graph subdivides the torus $T_H$ into regions. The
vertices interior to these regions are all tile vertices. Since two
tiles can be adjacent only if they have the same orientation, we
conclude that inside any one given region all the tiles have the same
orientation $x_i$ for some $i \in 1,2,3$.

\begin{exmpl} Let $G$ be the group $\frac{1}{13}(1,5,7)$, $Y$ be
$G$-$\hilb(\mathbb{C}^3)$ and let $\tilde{\mathcal{M}}$ be the dual of
the universal family of $G$-clusters (cf. Section
$\ref{section-worked-example}$).

In Figure \ref{figure-10} the thick lines depict the sink-source graph
$\sinksource_{\tilde{\mathcal{M}},E_7}$ of the exceptional divisor
$E_7$ (see Section \ref{section-example-the-resolution}) drawn on the
torus $T_H$ (left). The lift of the graph to the universal cover $H
\otimes \mathbb{R}^2$ is drawn on the right.

\begin{figure}[h] \begin{center}
\includegraphics[scale=0.37]{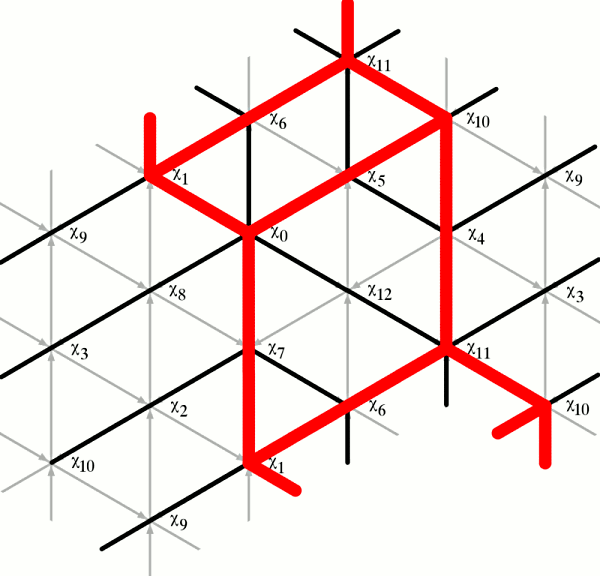} \hspace{0.5cm}
\includegraphics[scale=0.37]{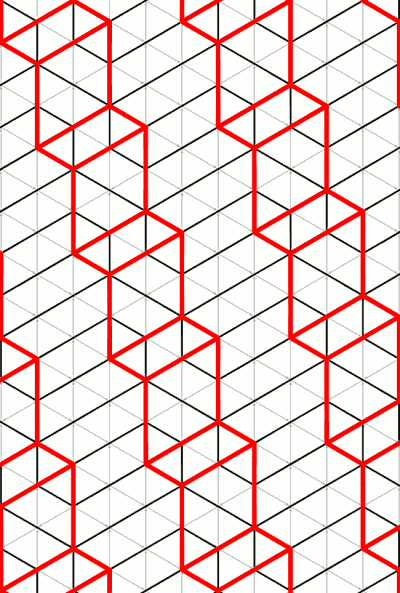}
\caption{\label{figure-10} The sink-source graph of $E_7$ on $T_H$
(left) and on $H \otimes \mathbb{R}^2$(right)} \end{center}
\end{figure}

Looking at the left diagram, there are two sources (the
$x_2$-$(1,2)$-source $\chi_{10}$ and the $x_2$-$(2,1)$-source
$\chi_1$) and two sinks (the $(3,0)$-sink $\chi_0$ and the
$(0,3)$-sink $\chi_{11}$).  \end{exmpl}

\subsection{Further Properties of $\mckquiv_{\mathcal{F}}$}

We now restrict our attention more to the case where $Y = G$-
$\hilb(\mathbb{C}^3)$ and $\mathcal{F} = \tilde{\mathcal{M}}$ is the
dual of the universal family of $G$-clusters.

Given a single $G$-cluster $\mathcal{V}$ its associated representation
$\mckquiv_\mathcal{V}$ is the vector space underlying the
corresponding $\twalg$-module $V = \Gamma(\mathcal{V})$.  The grading
$\oplus_\chi V_\chi$ is given by the natural decomposition of $V$ into
irreducible $G$-representations. The map corresponding to an arrow
$(\chi,x_k)$ is given by the multiplication map $V_\chi
\xrightarrow{\cdot x_k} V_{\chi \kappa(x_k)}$. 

We begin with a simple observation about $G$-clusters.  
\begin{prps}
\label{prps-gcluster-no-three-arrows} If $\mathcal{V}$ is a single
$G$-cluster then in the associated representation of $\mathcal{V}$ one
of the arrows leading into the vertex $\chi \ne \chi_0$ must be
non-zero. If $\mathcal{V}$ is the dual of a single $G$-cluster then
one of the arrows leading out of $\chi$ must be non-zero.  \end{prps}
\begin{proof} Since a $G$-cluster is of the form $\mathcal{O}_Z$ for a
subscheme $Z \subset \mathbb{C}^3$, the $\twalg$-module $V$
corresponding to $\mathcal{V}$ is of the form $\regring/I$ for some
$G$-invariant ideal $I$.  Therefore it is generated by $1$ as an
$\regring$-module. In particular, for any $\chi \neq \chi_0$ there
exists a monomial $1 \ne m \in R$ such that $m \cdot 1$ is a non-zero
element of $V_\chi$. Writing  $m = x_i m'$ for some $x_i$ we find that
the arrow labeled $x_i$ leading into $\chi$ must be a non-zero map. 
 
The dual claim follows since taking duals corresponds to reversing all
the arrows (and relabeling $\chi$ by $\chi^{-1}$).  \end{proof}

\begin{cor}\label{cor-gcluster-no-three-arrows} Let $Y =
G$-$\hilb(\mathbb{C}^3)$ and $\tilde{\mathcal{M}}$ the dual of the
universal family of $G$-clusters. If $\chi \in G^\vee$ is non-trivial
then the intersection of the three divisors
$B_{\chi,1}, B_{\chi,2}$ and $B_{\chi,3}$ associated to
$\mckquiv_{\tilde{\mathcal{M}}}$ is empty.  \end{cor} \begin{proof}
Over any point belonging to the intersection of $B_{\chi,1},
B_{\chi,2}$ and $B_{\chi,3}$ we get the dual of a $G$-cluster for
which all the maps leading out of $\chi$ are zero (remember that
$B_{\chi,i}$ records the locus where the map $\alpha_{\chi,i}$ is
zero). This cannot happen by Proposition
\ref{prps-gcluster-no-three-arrows} so the intersection must be empty.
\end{proof}

\begin{prps} \label{prps-ghilb-one-(0,3)-sink} Let $Y$ be
$G$-$\hilb(\C^3)$, $E \in \Except(Y)$ an irreducible exceptional
divisor and $\mathcal{M}$ the universal family of $G$-clusters. Then
the sink-source graph $\sinksource_{\mathcal{M},E}$ contains exactly
one $(0,3)$-sink (given by the vertex $\chi_0$). Similarly, the
sink-source graph $\sinksource_{\tilde{\mathcal{M}},E}$ contains
exactly one $(3,0)$-sink (given by the vertex $\chi_0$).  \end{prps}
\begin{proof} It suffices to prove the first claim as the dual claim
follows similarly. If $\chi$ is a $(0,3)$-sink in
$\sinksource_{\mathcal{M},E}$ then in the associated representation of
a fiber of $\mathcal{M}$ over any point in $E$ all three arrows
leading into $\chi$ are zero-maps. This, by Proposition
\ref{prps-gcluster-no-three-arrows}, can only happen if $\chi =
\chi_0$.  

It remains to show that $\chi_0$ is a $(0,3)$-sink in $\sinksource_{
\mathcal{M},E}$. If not then there exists a point $y \in E$ such that
in the associated representation of $\mathcal{M}|_y$ one of the three
arrows leading into $\chi_0$ is not a zero map. Since the same is also
true for every $\chi \neq \chi_0$ (by Proposition
\ref{prps-gcluster-no-three-arrows}), we can find a path of non-zero
maps in the McKay quiver which starts and ends at the same vertex.
Since each of these maps is multiplication by non-zero number their
composition is also non-zero. So this path gives a monomial $1 \ne m
\in \regring$ which is invariant under $G$ (since the path starts and
ends at the same vertex) such that $m \cdot 1 \neq 0$ in
$\mathcal{M}|_{y}$.

This is impossible since the resolution $Y \xrightarrow{\pi} \C^3/G$
is the Hilbert-Chow map of $\mathcal{M}$ and so for any monomial $1
\ne m \in \regring^G$ and $s \in \mathcal{M}|_{y}$ we have $m \cdot s
= m(\pi(y)) s = m(0) s = 0$. Here $\pi(y) = 0$ since $y$ lies in the
exceptional locus of $Y$.  \end{proof}

The next result is a consequence of \ref{prps-ghilb-one-(0,3)-sink}.
It will be heavily used to derive properties of the associated
representation $\mckquiv_{\tilde{\mathcal{M}}}$. 

\begin{prps} \label{prps-one-3,0-one-0,3} For any $Y$, $\mathcal{F}$
and $E$, if the sink-source graph $\sinksource_{\mathcal{F},E}$
contains exactly one $(3,0)$-sink then there are exactly three
possibilities for its remaining vertices: \begin{enumerate} \item
\label{case-one-(3,3)-emitter} One $(0,3)$-sink and one
$(3,3)$-source.  \item \label{case-one-(1,2)-one-(2,1)-emitter} One
$(0,3)$-sink, one $(1,2)$- and one $(2,1)$-source.  \item
\label{case-three-(1,2)-emitters}  Two $(0,3)$-sinks and three
$(1,2)$-sources.  \end{enumerate} \end{prps} \begin{proof} If there
exists only one $(3,0)$-sink, all the sources together must emit
exactly three charges of type $(1,0)$. This leaves us with three
possibilities for the sources: either one $(3,3)$-source, or one
$(1,2)$- and one $(2,1)$-source, or three $(1,2)$-sources. The
statement about sinks follows by counting in each case the number of
$(0,1)$-charges emitted.  \end{proof}

\begin{figure}[h] \centering \subfigure[One $(1,2)$-source and one
$(2,1)$-source] { \label{figure-11a}
\includegraphics[scale=0.18]{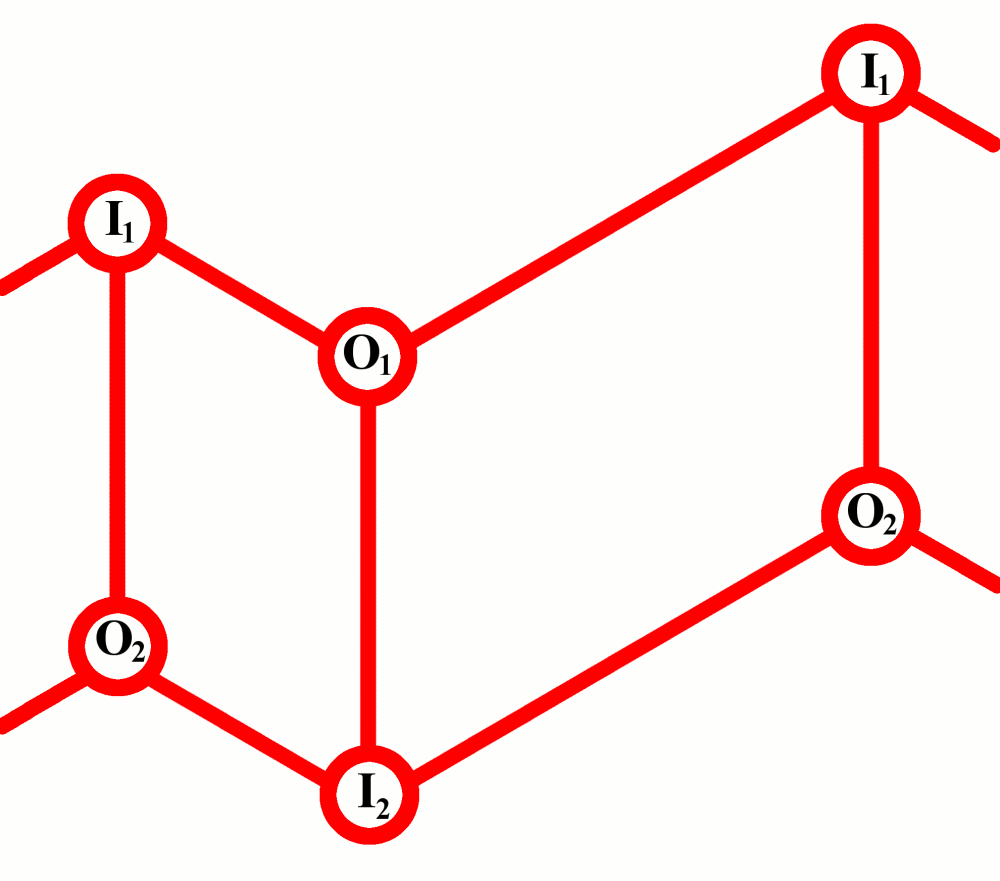} } \hspace{1cm}
\subfigure[A single $(3,3)$-source] { \label{figure-11b}
\includegraphics[scale=0.14]{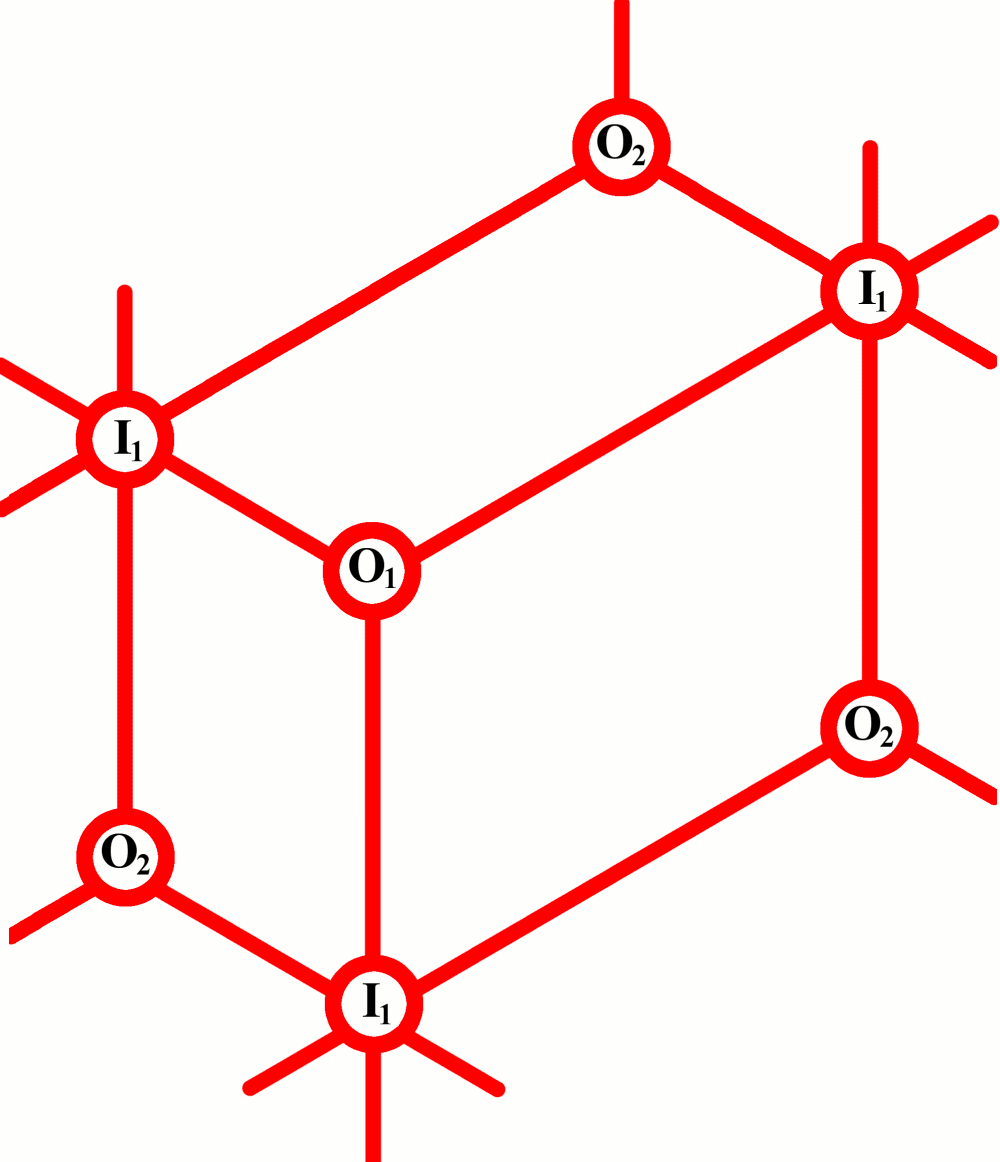} } \\ \subfigure[Three
$(1,2)$-sources] { \label{figure-11c}
\includegraphics[scale=0.22]{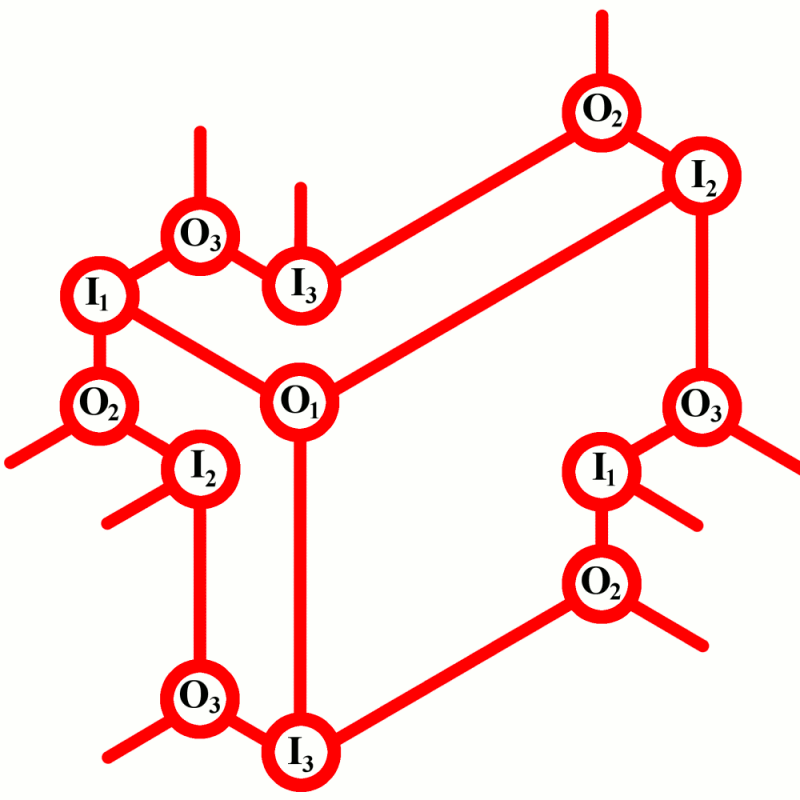} } \caption{Sink-source
graphs in the case of a single $(3,0)$-sink} \label{figure11}
\end{figure}

\begin{cor} \label{cor-shield-shaped-sas-graph} If the sink-source
graph $\sinksource_{\mathcal{F}, E}$ contains exactly one $(3,0)$-sink
then on the torus $T_H$ it looks like
either Figure $\ref{figure-11a}$, Figure $\ref{figure-11b}$ or Figure
\ref{figure-11c} (up to the lengths of the sides and in the case of
Figure \ref{figure-11a} up to a rotation by $\pm 2\pi/3$).  \end{cor}
\begin{proof} On the universal cover $H \otimes \mathbb{R}^2$ take any
vertex projecting to the $(3,0)$-sink $O_1$ in $T_H$ and denote it by
$A$.  Follow the three $(1,0)$-charge lines which enter $A$ back to
their respective sources. Denote the source emitting the
$x_i$-$(1,0)$-charge line by $B_i$. We obtain, up to the actual
lengths of $A B_i$, the picture in Figure \ref{figure-12}.
\begin{figure}[h] \begin{center}
\includegraphics[scale=0.28]{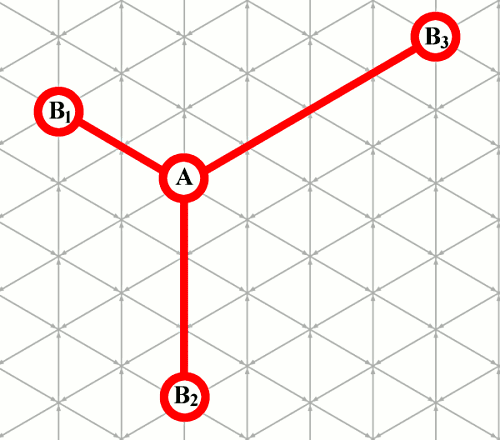}
\caption{\label{figure-12} The lift of $(3,0)$-sink and its three
charge lines to $H \otimes \mathbb{R}^2$} \end{center} \end{figure}

By  Proposition \ref{prps-one-3,0-one-0,3}, $B_1$, $B_2$ and $B_3$
project in $T_H$ either to the same $(3,3)$-source or to a
$(1,2)$-source and a $(2,1)$-source or to three distinct
$(1,2)$-sources. We treat here the case of a $(1,2)$-source and a
$(2,1)$-source, which leads to Figure \ref{figure-11a} and its
rotations. The other two cases are analogous and lead to
configurations in Figures \ref{figure-11b} and \ref{figure-11c}.

By our assumption exactly two of $B_i$ must project to same vertex in
$T_H$. Assume they are $B_1$ and $B_3$ - two other possibilities lead
to rotations of the resulting configuration.  Then the image of $B_1$
and $B_3$ in $T_H$ emits an $x_1$-$(1,0)$-charge and an
$x_2$-$(1,0)$-charge. This means it must be an $x_3$-$(2,1)$-source
(check with the list of source vertices in Figure \ref{figure-08}).
Since we have one $(3,0)$-sink and one $(0,3)$-sink, the two sources
must emit together one $x_i$-$(1,0)$- and one $x_i$-$(0,1)$-charge for
$i \in 1,2,3$. The source $B_2$ projects to  has therefore to be
an $x_3$-$(1,2)$-source. 

Consider now the $x_1$-$(0,1)$-charge line emitted by $B_2$ and the
$x_3$-$(0,1)$-charge line emitted by $B_1$. We claim that the vertex
$B_{12}$ at the intersection of these lines projects to the
$(0,3)$-sink in $T_H$. If not then one of the two charge lines would
have to terminate in the $(0,3)$-sink before intersecting the other.
Then, from Figure $\ref{figure-12}$, one of the other charge lines
ending in this sink would have to come from within the parallelogram
$AB_1 B_2 B_{12}$.  As this charge line can't cross $AB_1$ or $AB_2$,
its source must also lie within this parallelogram. This is impossible
- every vertex in edges $A B_1$ and $A B_2$ projects to a distinct
  vertex in $T_H$, and therefore every vertex in the parallelogram $A
B_1 B_2 B_{12}$ must also project to a distinct vertex. But by
assumption there are only two source vertices in $T_H$ and $B_1$ and
$B_2$ already project to those. Therefore no other vertex within $A
B_1 B_2 B_{12}$ can project to a source vertex.  

Thus $B_{12}$ projects to the $(0,3)$-sink and drawing in all the
remaining charge lines yields the configuration on Figure
\ref{figure-11a}.  \end{proof}

We now present a few more properties of the associated representation
of $\tilde{\mathcal{M}}$ -- the dual of the universal family of
$G$-clusters.

\begin{figure}[h]
\centering
\subfigure[ ]
{
    \label{figure-13a}
    \includegraphics[scale=0.12]{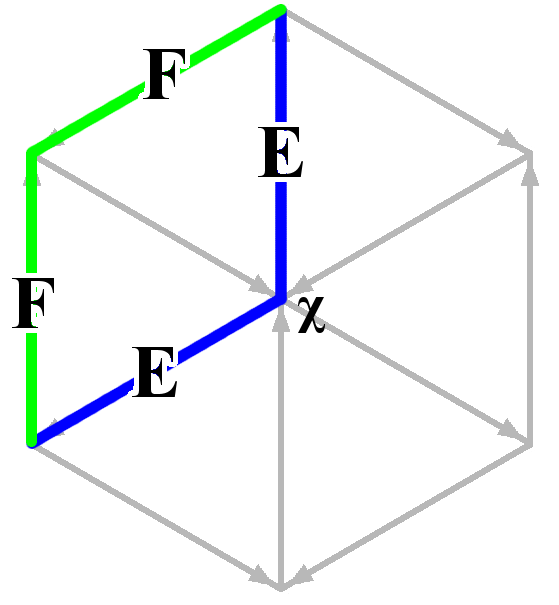}
}
\hspace{1cm}
\subfigure[ ]
{
    \label{figure-13b}
    \includegraphics[scale=0.12]{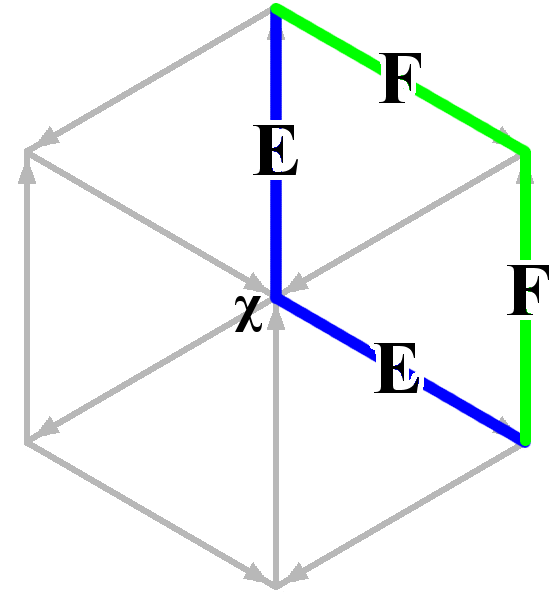}
}
\hspace{1cm}
\subfigure[ ]
{
    \label{figure-13c}
    \includegraphics[scale=0.12]{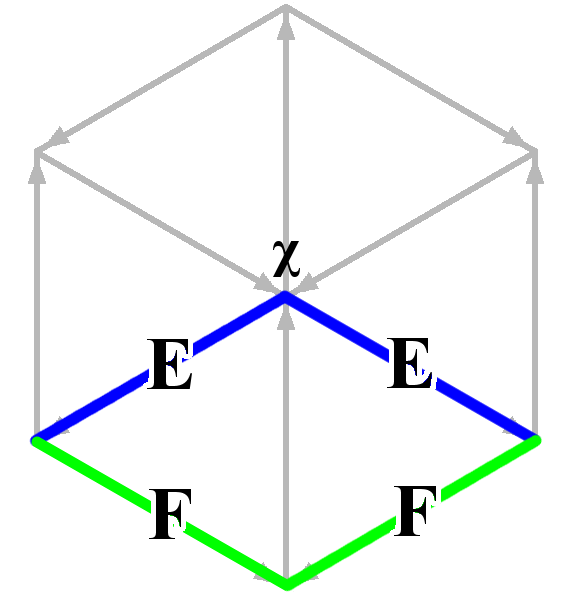}
} \caption{The configurations from Proposition
\ref{prps-opposite-oriented-charges-do-not-pass-through-each-other}}
\label{figure-13} 
\end{figure}

\begin{prps}
\label{prps-opposite-oriented-charges-do-not-pass-through-each-other}
Let $Y$ be $G$-$\hilb(\mathbb{C}^3)$ and $\tilde{\mathcal{M}}$ be the dual 
of the universal family of $G$-clusters. Let $E$ and $F$ be two
irreducible exceptional divisors on $Y$. Suppose $\chi \in G^\vee$ 
such that in $\mckquiv_{\tilde{\mathcal{M}}}$ we have
$E \subseteq B_{\chi, i}, B_{\chi,j}$, while
$F \subseteq B_{\chi\kappa(x_i),j}, B_{\chi\kappa(x_j),i}$ for some
$i \neq j$, cf. Figure \ref{figure-13}. Then $E$ and $F$ do not intersect. 
\end{prps}
\begin{proof}
We give the proof under the assumption that the sink-source graph 
$\sinksource_{\tilde{\mathcal{M}},F}$ is the one depicted 
in Figure \ref{figure-11a}. The cases of
$\sinksource_{\tilde{\mathcal{M}},F}$ being
as in Figures \ref{figure-11b} and \ref{figure-11c} are similar and
simpler, each requiring only a part of the argument presented below.

We use the notation in Figure $\ref{figure-11a}$ to label the vertices
of $\sinksource_{\tilde{\mathcal{M}},F}$: $O_1$ is the $(3,0)$-sink,
$O_2$ the $(0,3)$-sink, $I_1$ the $x_3$-$(2,1)$-source and $I_2$ the
$x_3$-$(1,2)$-source. By Proposition \ref{prps-ghilb-one-(0,3)-sink}
$O_1 = \chi_0$ is also the unique $(3,0)$-sink in
$\sinksource_{\tilde{\mathcal{M}},E}$. 

We first deal with the configuration in Figure $\ref{figure-13c}$.  The proof will require working with both sink-source graphs $\sinksource_{\tilde{\mathcal{M}},E}$ and $\sinksource_{\tilde{\mathcal{M}},F}$ which we superimpose on the torus $T_H$. 

First note that the vertex $\chi$ must be an $x_3$-$(1,0)$-charge for $E$.
All such charges lie on the $x_3$-$(1,0)$-charge line in 
$\sinksource_{\tilde{\mathcal{M}},E}$ which terminates at $O_1$.  
We claim that the source in $\sinksource_{\tilde{\mathcal{M}},F}$ for this charge line must lie below the vertex $\kappa(x_3)^{-1} I_2$ - see Figure \ref{figure-11a}. If this were not the case then every vertex which is an $x_3$-$(1,0)$-charge for $E$ would be either an $x_3$-$(1,0)$ charge or the source $I_2$ for $F$. This would render Figure $\ref{figure-13c}$ impossible. Hence $\kappa(x_3)^{-1} I_2$ is an $x_3$-$(1,0)$-charge for $E$ and we have the configuration depicted in Figure \ref{figure-14}.
\begin{figure}[h]
\begin{center}
\includegraphics[scale=0.15]{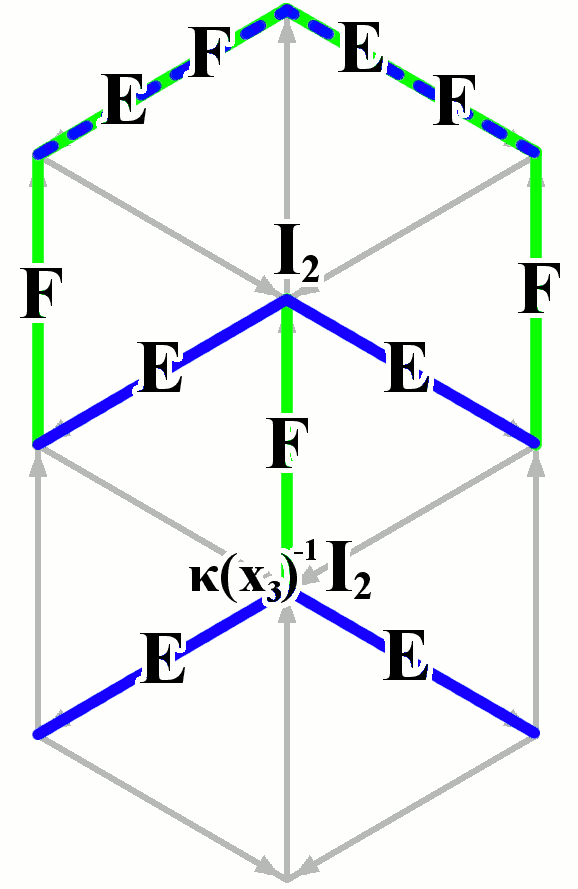}
\caption{\label{figure-14} Partial divisor configuration for $E$
and $F$ near $\kappa(x_3)^{-1} I_2$}
\end{center}
\end{figure}
Since all three arrows leading out of $\kappa(x_3)^{-1} I_2$
are marked by either $E$ or $F$ we must have $E \cap F = \emptyset$
by Corollary \ref{cor-gcluster-no-three-arrows}.

Next we assume the configuration in Figure $\ref{figure-13b}$. 
This also deals with the configuration in Figure $\ref{figure-13a}$
as the argument there is identical.

As we assumed that $\sinksource_{\tilde{\mathcal{M}}, F}$ looks 
like in Figure $\ref{figure-11a}$, for $F$ vertex $\chi$ can 
only be an $x_2$-$(0,1)$-charge or the $x_3$-$(1,2)$-source $O_2$.
So $\chi$ must lie on the $x_2$-$(0,1)$-charge line $I_2 O_2$ 
in $\sinksource_{\tilde{\mathcal{M}},F}$. On 
the other hand $\chi$ must lie on the $x_2$-$(1,0)$-charge 
in $\sinksource_{\tilde{\mathcal{M}},E}$ which terminates at $O_1$
(remember that $O_1$ is also the $(3,0)$-sink for $E$).
Looking at Figure \ref{figure-11a} one sees that on the torus $T_H$
the path along this charge line from $\chi$ to $O_1$ must
pass through both $I_2$ and $I_1$. 
Since charge lines for $E$ can't cross each other, the source 
in $\sinksource_{\tilde{\mathcal{M}},E}$
for the $x_1$-$(1,0)$-charge (resp. the $x_3$-$(1,0)$-charge) line
must lie \em strictly \rm between $O_1$ and $I_1$ (resp. between $O_1$
and $I_2$). Denoting this source by $P_1$ (resp. $P_2$) we obtain the 
configuration in Figure \ref{figure-15}. 
\begin{figure}[h]
\begin{center}
\includegraphics[scale=0.18]{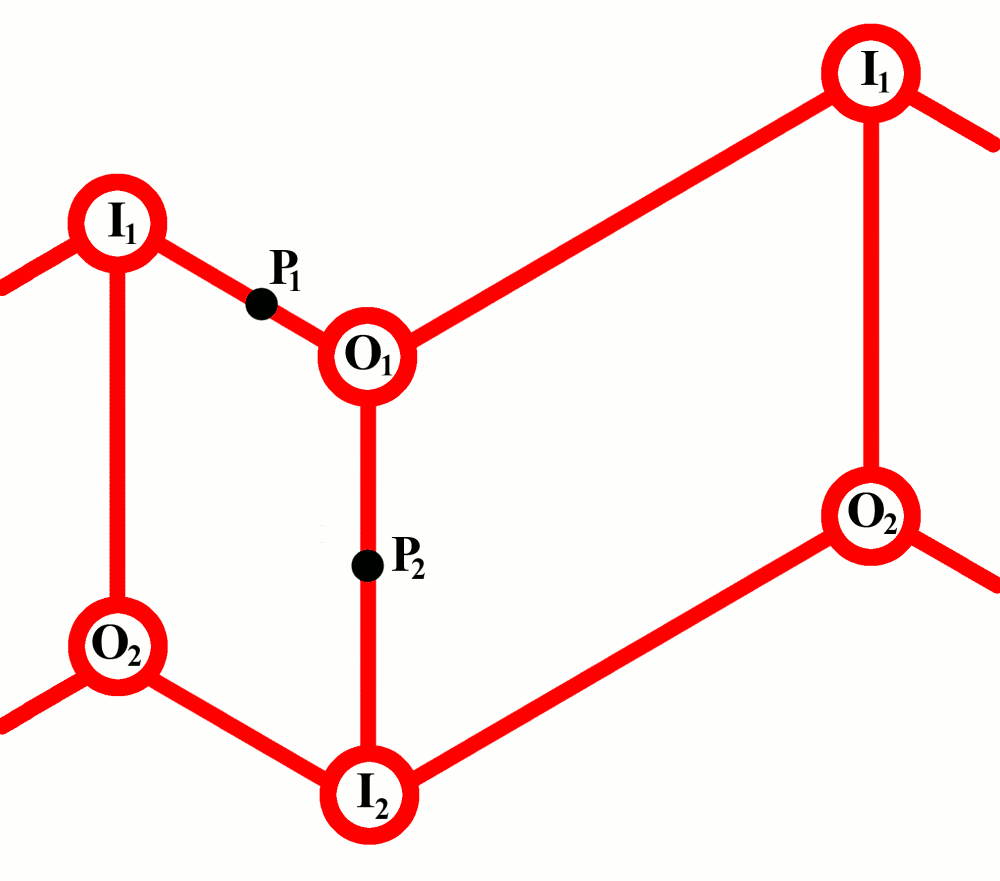}
\caption{\label{figure-15} The sink-source graph
$\sinksource_{\tilde{\mathcal{M}}, F}$ with 
the sources $P_1$ and $P_2$ for $E$ }
\end{center}
\end{figure}

Suppose $\sinksource_{\tilde{\mathcal{M}},E}$ is shaped as in Figures
\ref{figure-11a} or \ref{figure-11b}. Then it has at most two sources.
Therefore the $x_2$-$(1,0)$-charge line in $\sinksource_{\tilde{\mathcal{M}},E}$
must originate at either $P_1$ or $P_2$. Suppose it originates
at $P_1$, the other case is analogous. Then $P_1$ is an 
$x_3$-$(2,1)$-source for $E$ since this is the unique source 
which emits both $x_2$-$(1,0)$- and $x_1$-$(1,0)$-charges.
We claim that $\kappa(x_2) P_1$ lies within 
the left parallelogram $O_1 I_2 O_2 I_1$ in Figure \ref{figure-15}. 
Then $\kappa(x_2) P_1$ has to be an $x_2$ tile for $F$ which gives
the configuration in Figure \ref{figure-16}, whence 
$E \cap F = \emptyset$ by Corollary \ref{cor-gcluster-no-three-arrows}.
Observe - as $I_1$ is an $x_2$-$(1,0)$-charge for $E$, $P_1$ can not 
be the vertice $\kappa(x_1) I_1$, by direct inspection.
Also, the side $O_1 P_2 I_2$ is at least two edges long as 
it has to accommodate three vertices. These two observations
show that, as claimed, $\kappa(x_2) P_1$ must lie in the interior of 
$O_1 I_2 O_2 I_1$.

\begin{figure}[h]
\begin{center}
\includegraphics[scale=0.15]{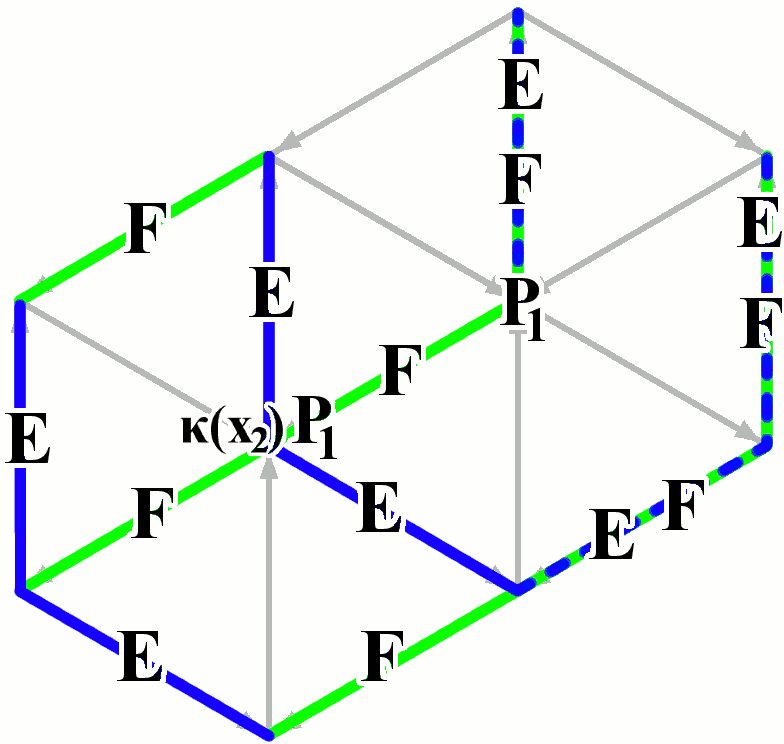}
\caption{\label{figure-16} Divisor configuration for $E$ and $F$ near
$\kappa(x_2) P_1$}
\end{center}
\end{figure}

Finally, suppose $\sinksource_{\tilde{\mathcal{M}},E}$ is shaped as in
Figure \ref{figure-11c}. Then $x_2$-$(1,0)$-charge line in
$\sinksource_{\tilde{\mathcal{M}},E}$ originates at an
$x_2$-$(1,2)$-source distinct from $P_1$ and $P_2$. Denote it by
$P_3$. But then looking at Figure \ref{figure-11c} (and keeping 
in mind that there $P_1$ is $O_1$, $P_2$ is $O_3$ and $P_3$ is $O_2$)
we see that the vertex $\kappa(x_2) P_3$ must lie within the
parallelogram formed by $O_1 P_1$ and $O_1 P_2$. Therefore, from
Figure $\ref{figure-15}$ it certainly lies within the left
parallelogram $O_1 I_2 O_2 I_1$. We conclude that $\kappa(x_2) P_3$
is an $x_2$-tile in $\sinksource_{\tilde{\mathcal{M}},F}$. This gives
a configuration similar to that in Figure \ref{figure-16}, showing 
again that $E \cap F = \emptyset$ by Corollary 
\ref{cor-gcluster-no-three-arrows}. \end{proof}

\begin{figure}[h]
\centering
\subfigure[ ]
{
    \label{figure-25a}
    \includegraphics[scale=0.12]{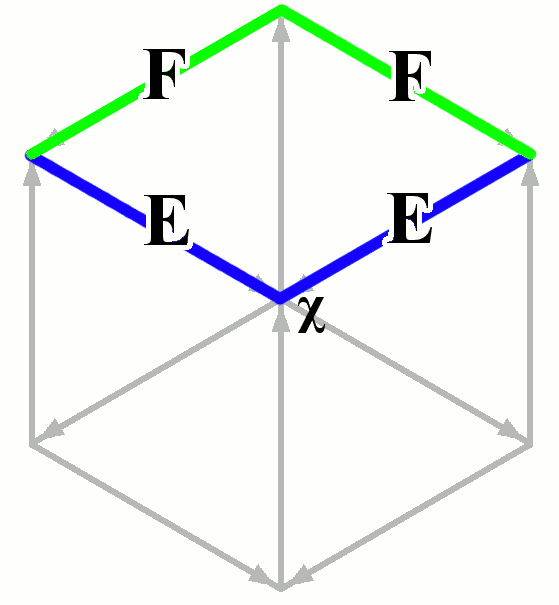}
}
\hspace{1cm}
\subfigure[ ]
{
    \label{figure-25b}
    \includegraphics[scale=0.12]{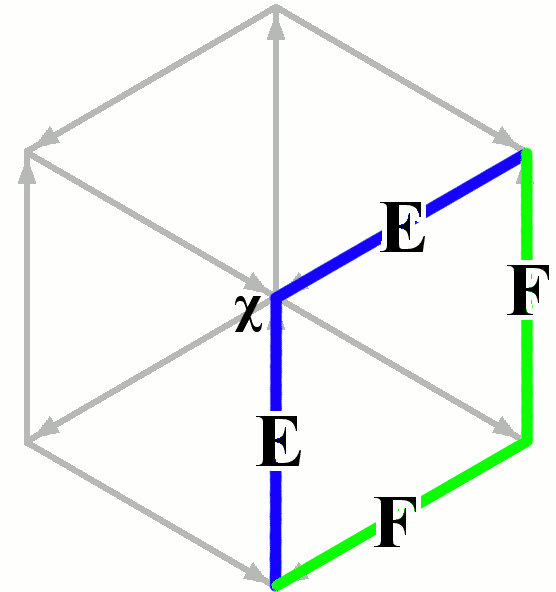}
}
\hspace{1cm}
\subfigure[ ]
{
    \label{figure-25c}
    \includegraphics[scale=0.12]{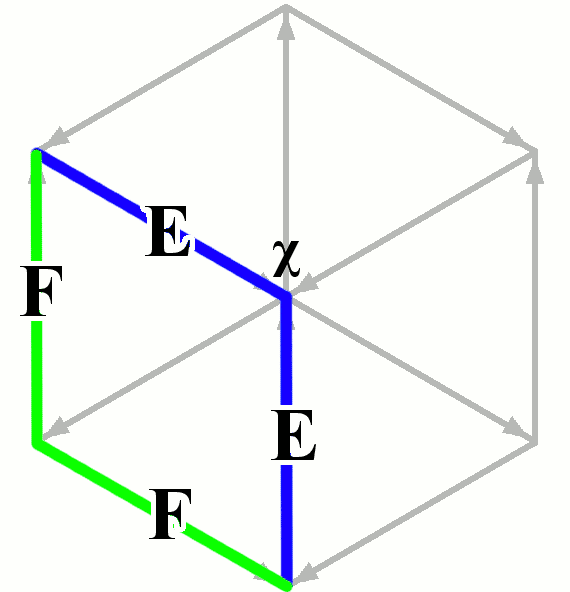}
} \caption{The configurations from Corollary 
\ref{cor-opposite-oriented-charges-do-not-pass-through-each-other}}
\label{figure-25} 
\end{figure}

\begin{cor}\label{cor-opposite-oriented-charges-do-not-pass-through-each-other}
Let $Y$ be $G$-$\hilb(\mathbb{C}^3)$ and $\tilde{\mathcal{M}}$ be 
the dual of the universal family of $G$-clusters. Let $E$ and $F$ be 
two irreducible exceptional divisors on $Y$. Suppose $\chi \in G^\vee$ 
such that in the associated representation of $\tilde{\mathcal{M}}$ we have
$E \subseteq B_{\chi \kappa(x_i)^{-1}, i} , B_{\chi \kappa(x_j)^{-1},
j}$ and $F \subseteq B_{\chi \kappa(x_k), i}, B_{\chi \kappa(x_k),
j}$ for some $i \neq j \neq k$, cf. Figure \ref{figure-25}. Then $E$
and $F$ do not intersect.
\end{cor}
\begin{proof}
If in the configuration in Figure \ref{figure-25a} we re-center on the vertex $\kappa(x_3)\chi$, we obtain in $\hex(\kappa(x_3)\chi)$ the configuration in Figure \ref{figure-13a}. By Proposition \ref{prps-opposite-oriented-charges-do-not-pass-through-each-other} we obtain $E \cap F = \emptyset$. The other two cases follow similarly.  
\end{proof}

\begin{figure}[h]
\centering
\subfigure[ ]
{
    \label{figure-17a}
    \includegraphics[scale=0.12]{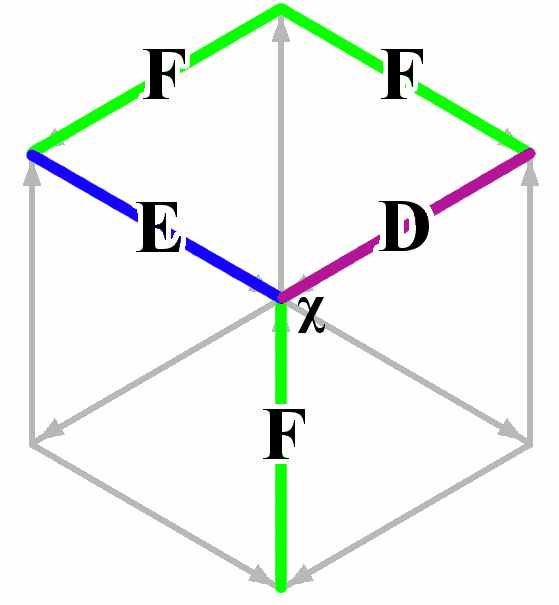}
}
\hspace{1cm}
\subfigure[ ]
{
    \label{figure-17b}
    \includegraphics[scale=0.12]{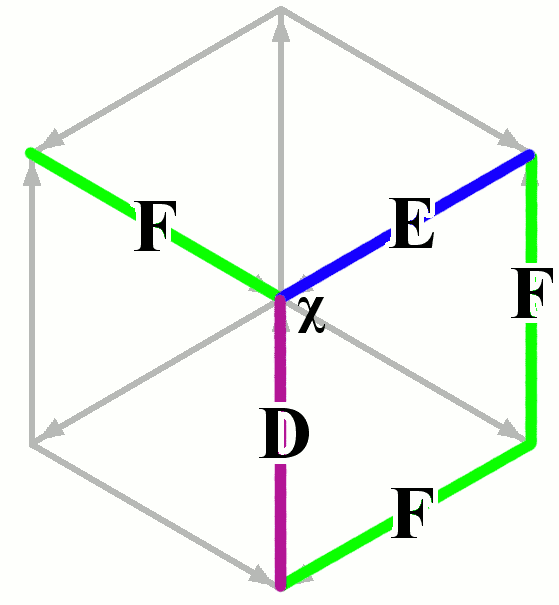}
}
\hspace{1cm}
\subfigure[ ]
{
    \label{figure-17c}
    \includegraphics[scale=0.12]{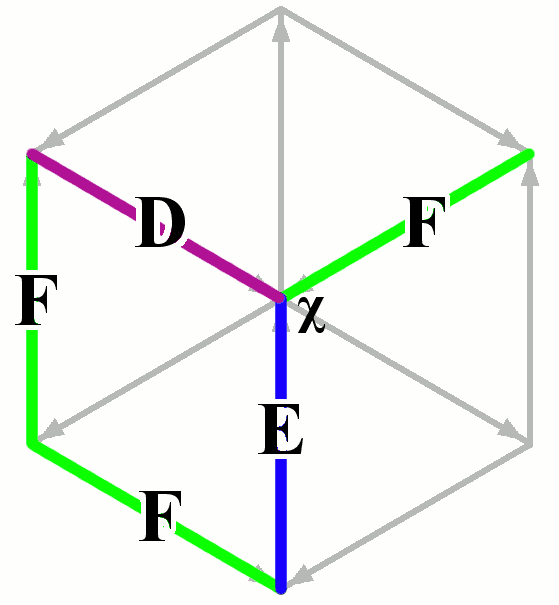}
} \caption{ The configurations from Proposition \ref{prps-no-comets}}
\label{figure-17} 
\end{figure}

\begin{prps}
\label{prps-no-comets}
Let $Y$ be $G$-$\hilb(\mathbb{C}^3)$ and $\tilde{\mathcal{M}}$ be the dual 
of the universal family of $G$-clusters. Let $D$, $E$ and $F$ be 
irreducible exceptional divisors on $Y$. Suppose $\chi \in G^\vee$ such 
that in $\mckquiv_{\tilde{\mathcal{M}}}$ we have
$D \in B_{\chi\kappa(x_i)^{-1}, i}$, $E \in B_{\chi\kappa(x_j)^{-1},j}$
and $F \in B_{\chi\kappa(x_k),i}, B_{\chi\kappa(x_k),j},
B_{\chi\kappa(x_k)^{-1},k}$ for some $i \neq j \neq k$, cf. 
Figure \ref{figure-17}. Then $D \cap E \cap F = \emptyset$. 
\end{prps}
\begin{proof}
We give the proof for the configuration in Figure $\ref{figure-17a}$, 
the other two cases are proven similarly. 

Comparing Figure $\ref{figure-17a}$ to Figures
$\ref{figure-05}$-$\ref{figure-09}$ we find that 
in $\sinksource_{\tilde{\mathcal{M}}, F}$ the vertex $\chi$ can only be the $x_3$-$(1,2)$-source. 
Therefore $\sinksource_{\tilde{\mathcal{M}}, F}$ is as depicted in
Figure \ref{figure-11a} (no rotation necessary) or in Figure
$\ref{figure-11c}$. Assume the former. Using the same notation as 
in Figure \ref{figure-11a}, note that $\chi = I_2$. 

Denote by $P_1$ the source of $x_3$-$(1,0)$-charge line for $E$. 
Since this charge line also terminates at $O_1$ the source $P_1$ must lie
below $O_1$. If it lies below $I_2$ then the configuration in Figure
\ref{figure-14} would occur, implying $E \cap F = \emptyset$ (by
Corollary \ref{cor-gcluster-no-three-arrows}). $P_1$ cannot equal
$I_2$ since then $\hex(I_2)$ would not look as in Figure
\ref{figure-17a}. So we can assume $P_1$ lies strictly between $O_1$
and $I_2$.

Suppose now that $P_1 \neq \kappa(x_3) I_2$. Then 
$\kappa(x_3)^{-1} P_1$ is an $x_3$-$(1,0)$-charge in
$\sinksource_{\tilde{\mathcal{M}}, F}$. Also, $P_1$ must be either 
an $x_3$-$(1,2)$-source, an $x_1$-$(2,1)$-source or an $x_2$-$(2,1)$-source
in $\sinksource_{\tilde{\mathcal{M}},E}$. If it is an 
$x_3$-$(1,2)$-source we get the configuration depicted in
Figure \ref{figure-18}, which by Corollary 
\ref{cor-gcluster-no-three-arrows} means $E \cap F = \emptyset$.
\begin{figure}[h]
\begin{center}
\includegraphics[scale=0.15]{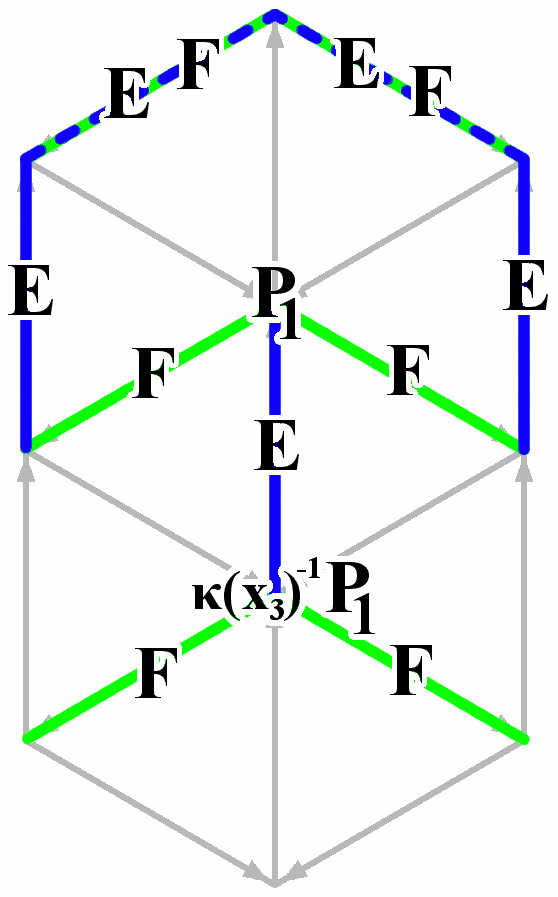}
\caption{\label{figure-18} Partial divisor configuration for $E$ and $F$ near $P_1$}
\end{center}
\end{figure}
If $P_1$ is an $x_2$-$(2,1)$-source, then we claim that 
$\kappa(x_1) P_1$ is an interior vertex of the right parallelogram 
in Figure \ref{figure-11a} and so we have the configuration in 
Figure $\ref{figure-19}$, implying $E \cap F = \emptyset$. The
argument for $P_1$ being an $x_1$-$(2,1)$ source is the same 
but with the vertex $\kappa(x_2) P_1$.
\begin{figure}[h]
\begin{center}
\includegraphics[scale=0.15]{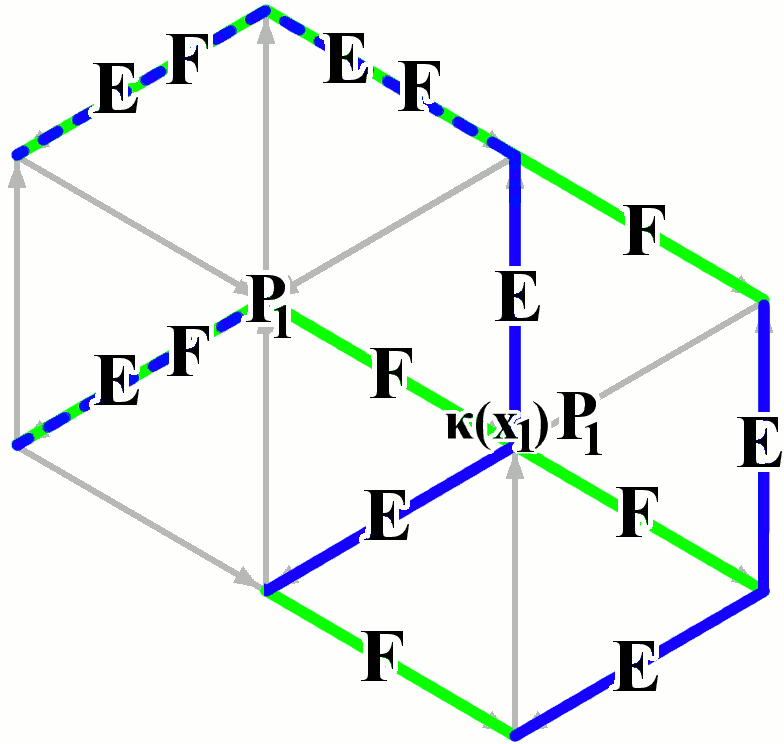}
\caption{\label{figure-19} Divisor configuration for $E$ and $F$ near $\kappa(x_1) P_1$}
\end{center}
\end{figure}

It suffices to show that the side $O_1 I_1$ is more than one edge
long, as then $\kappa(x_1) P_1$ can't lie on the side $I_1 O_2$
and has to be an interior vertex. Observe that 
the $x_1$-$(1,0)$-charge line for $E$
begins at $P_1$ and terminates at $O_1$, and so it 
must pass through $I_1$ (see Figure \ref{figure-11a}). So 
$I_1$ is an $x_2$-$(1,0)$-charge for $E$. Therefore the 
side $O_1 I_1$ is more than one edge long -- a $(3,0)$-source
$O_1$ and an $x_2$-$(1,0)$-charge $I_1$ could not be adjacent 
in such a way. 

Finally, suppose $P_1 = \kappa(x_3) I_2$. Recall that in 
$\sinksource_{\tilde{\mathcal{M}}, E}$ the vertex $P_1$ is
either an $x_3$-$(1,2)$-source, an $x_1$-$(2,1)$-source or 
an $x_2$-$(2,1)$-source. The first and the third
options are incompatible with the data of Figure \ref{figure-17a} 
since $\chi = I_2$. So $P_1$ is an $x_1$-$(2,1)$-source for $E$.
Now the $x_2$-$(1,0)$-charge line for $E$ begins at $P_1$
and terminates at $O_1$ so, looking at Figure \ref{figure-11a}, 
it must pass through $I_1$. We conclude that $I_1$ is 
an $x_2$-$(1,0)$-charge for $E$. 

Repeating the entire argument from the start with $D$ instead of $E$ 
yields that either $D \cap F = \emptyset$ or $I_1$ is an 
$x_1$-$(1,0)$-charge for $D$. So $I_1$ is an $x_1$-$(1,0)$-charge 
for $D$, an $x_2$-$(1,0)$-charge for $E$ and
the $x_3$-$(1,2)$-source for $F$ (its original definition). 
This gives the configuration in Figure $\ref{figure-20}$, 
implying $D \cap E = \emptyset$ by Corollary \ref{cor-gcluster-no-three-arrows}. 
\begin{figure}[h]
\begin{center}
\includegraphics[scale=0.15]{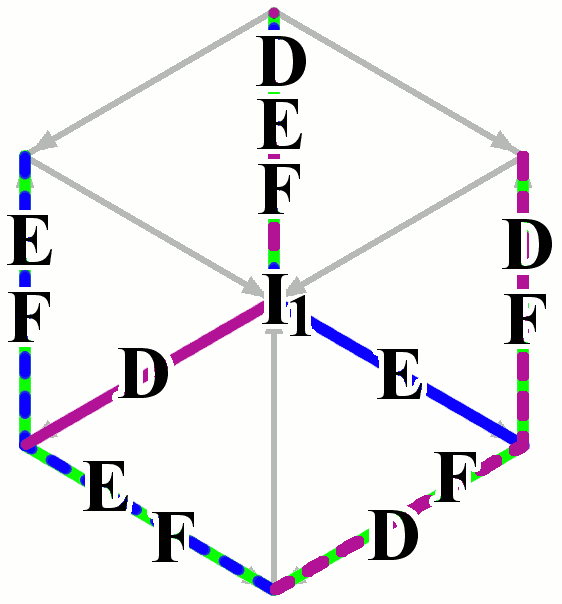}
\caption{\label{figure-20} Divisor configuration for $D$, $E$, and $F$ near $I_1$}
\end{center}
\end{figure}

The argument for $\sinksource_{\tilde{\mathcal{M}}, F}$ being as in Figure \ref{figure-11c} is very similar, so we only give a brief outline. We again denote by $P_1$ the source of $x_3$-$(1,0)$-charge line for $E$. But now the case $P_1 = \kappa(x_3) I_3$ leads immediately to a contradicting configuration at $\kappa(x_2)^{-1} I_2$. The difficulty is when $P_1$ is an $x_1$-$(2,1)$-source lying somewhere higher between $O_1$ and $I_3$ while $\kappa(x_2) P_1$ is $I_2$. It turns out in such case the $x_1$-$(1,2)$-source for $E$ has to be $\kappa(x_1)I_1 $ (otherwise $\chi = I_3$ would be an $x_3$-tile for $E$ contradicting Figure $\ref{figure-17a}$).  Denoting by $P_2$ the source of $x_3$-$(1,0)$-charge line for $D$ and repeating the same argument as for $E$, we now have deal with the case where $P_2$ is an $x_2$-$(2,1)$-source lying somewhere on $O_1 I_3$ and $\kappa(x_1) P_2 = I_1$. But then $P_2 = \kappa(x_1)^{-1} I_1$ is an $x_2$-$(2,1)$-source for $D$ and $\kappa(x_1)I_1$ is a $x_1$-$(1,2)$-source for $E$, which creates a contradicting configuration at the vertex $I_1$ lying between them. 
\end{proof}

\begin{figure}[h]
\centering
\subfigure[ ]
{
    \label{figure-21a}
    \includegraphics[scale=0.12]{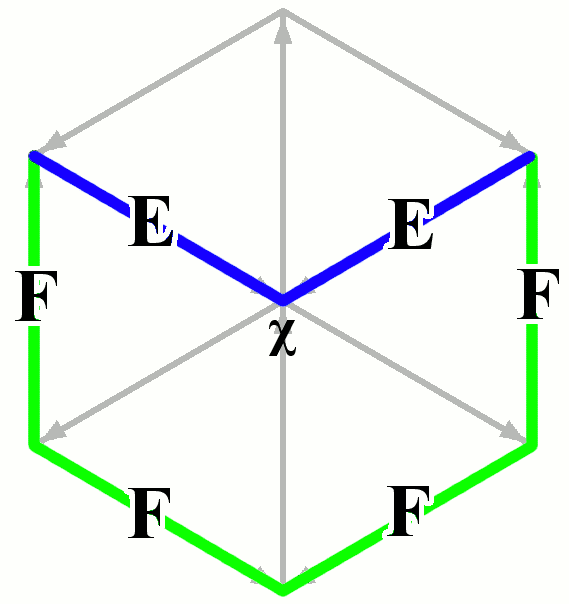}
}
\hspace{1cm}
\subfigure[ ]
{
    \label{figure-21b}
    \includegraphics[scale=0.12]{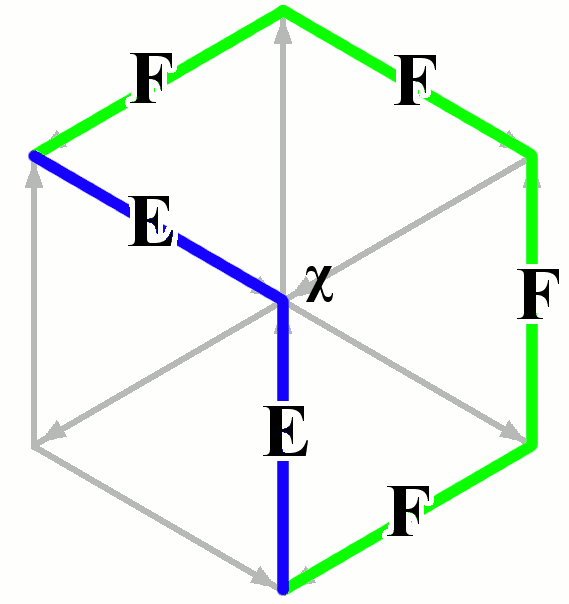}
}
\hspace{1cm}
\subfigure[ ]
{
    \label{figure-21c}
    \includegraphics[scale=0.12]{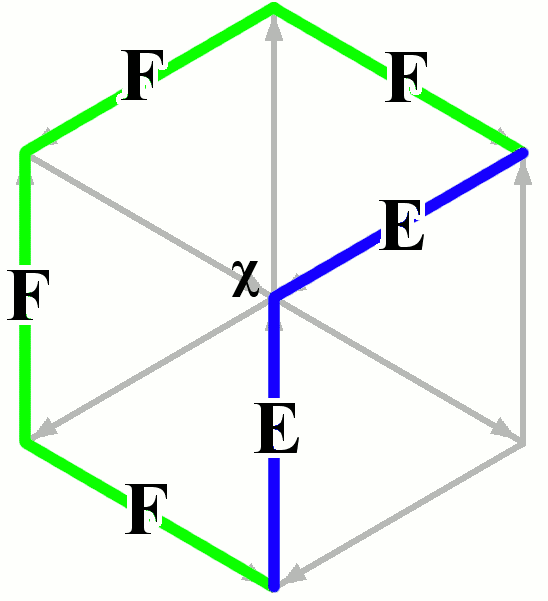}
} \caption{The configurations from Proposition \ref{prps-no-shields}}
\label{figure-21} 
\end{figure}

\begin{prps}
\label{prps-no-shields}
Let $Y$ be $G$-$\hilb(\mathbb{C}^3)$ and $\tilde{\mathcal{M}}$ be the dual 
of the universal family of $G$-clusters. Let $E$ and $F$ be 
irreducible exceptional divisors on $Y$. Suppose there exists 
$\chi \in G^\vee$ such that for $\tilde{\mathcal{M}}$ we have
$E \in B_{\chi\kappa(x_i)^{-1}, i}, B_{\chi\kappa(x_j)^{-1},j}$
and $F \in B_{\chi\kappa(x_i),j}, B_{\chi\kappa(x_i),k}, B_{\chi\kappa(x_j),i},
B_{\chi\kappa(x_j),k}$ for some
$i \neq j \neq k$, cf. Figure \ref{figure-21}. Then $E \cap F = \emptyset$. 
\end{prps}
\begin{proof}
We give the proof for the configuration in Figure $\ref{figure-21a}$. 
The other two are identical.

Comparing Figure $\ref{figure-21a}$ with Figures $\ref{figure-05}$-$\ref{figure-09}$ we conclude that in $\sinksource_{\tilde{\mathcal{M}},F}$ the vertex $\chi$ must be a $(3,3)$-source or an $x_3$-$(2,1)$-source. 

Suppose the latter. Then $\sinksource_{\tilde{\mathcal{M}},F}$ is 
as depicted in Figure \ref{figure-11a} (no rotation necessary).
Observe that $\chi = I_1$.  
Let $P_1$ and $P_2$ denote the sources of the $x_1$- and $x_2$-$(1,0)$-charge 
lines for $E$, respectively. From Figure $\ref{figure-21a}$ 
we see that $I_1$ can neither be an $x_1$-$(1,0)$-charge nor
$x_2$-$(1,0)$-charge for $E$. Therefore 
$P_1$ (resp. $P_2$) lies on the 
$x_1$-$(1,0)$-charge line (resp. the $x_2$-$(1,0)$-charge line)
as depicted in Figure \ref{figure-22}. 
\begin{figure}[h]
\begin{center}
\includegraphics[scale=0.21]{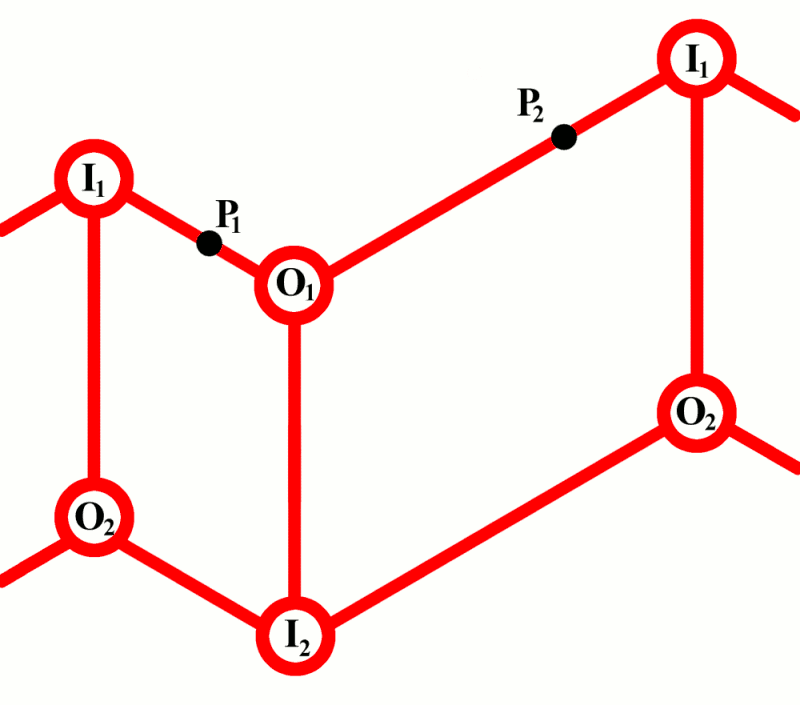}
\caption{\label{figure-22} The sink-source graph
$\sinksource_{\tilde{\mathcal{M}}, F}$ with sources $P_1$ and $P_2$ for $E$ }
\end{center}
\end{figure}

Since $P_1 \neq P_2$ the graph $\sinksource_{\tilde{\mathcal{M}},E}$
is as in \ref{figure-11a} (up to rotation) or as in \ref{figure-11c}.
In former case we find that either $P_1$ is an $x_1$-$(1,2)$-source
or $P_2$ is an $x_2$-$(1,2)$-source and in latter case we find that
both of these are true. Suppose $P_1$ is an $x_1$-$(1,2)$-source, 
the other case is similar. Since $I_1 = \chi$ the vertex $I_1
\kappa(x_1)$ can't be an $x_1$-$(1,2)$-source for $E$ by inspection of
Figure $\ref{figure-17a}$.  So $P_1 \neq I_1 \kappa(x_1)$ and hence
$\kappa(x_1)^{-1} P_1$ is an $x_1$-$(1,0)$-charge for $F$ and get the
configuration depicted in Figure $\ref{figure-23}$. As usual, it
implies $E \cap F = \emptyset$ by Corollary
\ref{cor-gcluster-no-three-arrows}. 

\begin{figure}[h] \begin{center}
\includegraphics[scale=0.15]{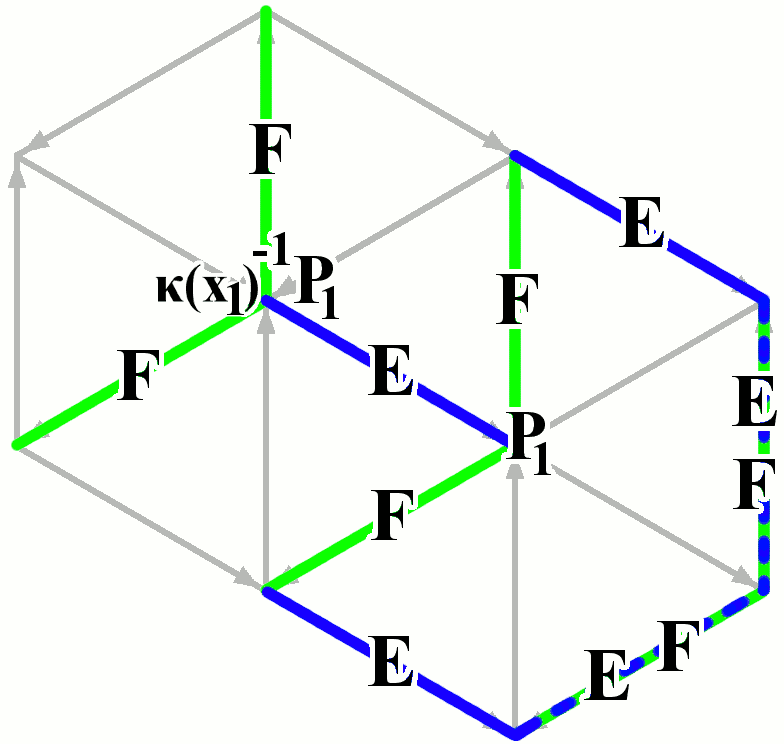}
\caption{\label{figure-23} Divisor configuration for $E$ and $F$ 
near $P_1$}
\end{center}
\end{figure}
\begin{figure}[h]
\begin{center}
\includegraphics[scale=0.18]{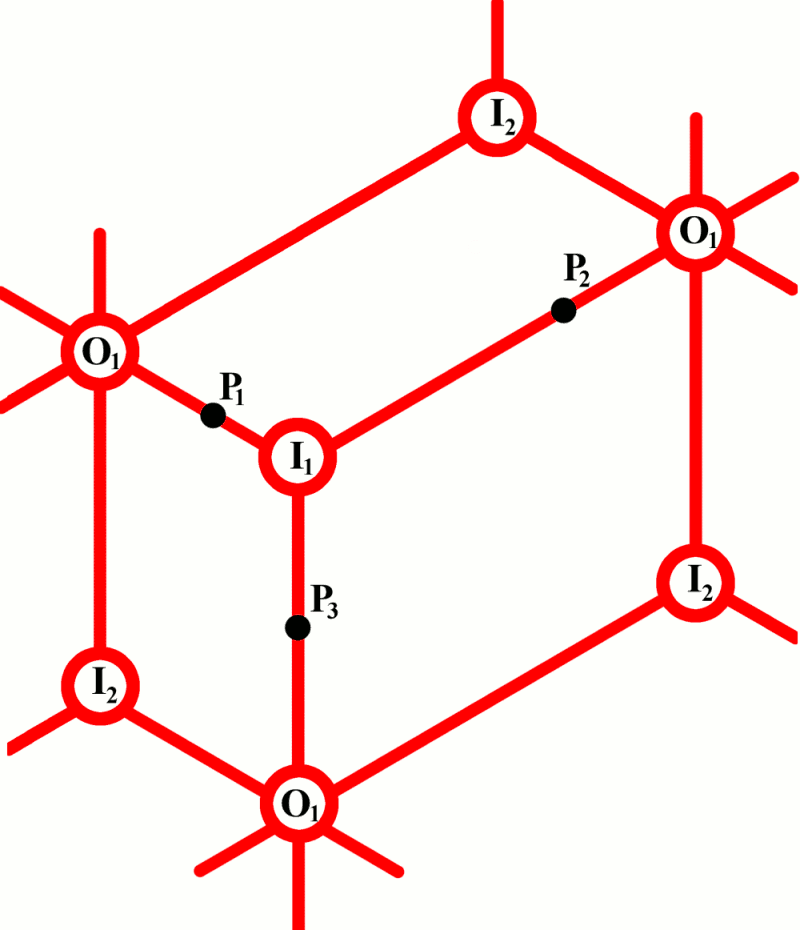}
\caption{\label{figure-24} The sink-source graph
$\sinksource_{\tilde{\mathcal{M}}, F}$ with 
sources $P_1$,$P_2$ and $P_3$ for $E$}
\end{center}
\end{figure}

It remains to consider the case when $\chi$ is a $(3,3)$-source for
$F$. This time the sink-source graph $\sinksource_{\tilde{M},
\mathcal{F}}$ is as in Figure \ref{figure-11b}.  Reasoning as above,
the vertices $P_1$, $P_2$ and $P_3$ are configured as in Figure
\ref{figure-24}. However, observe that in both Figures
\ref{figure-11b} and \ref{figure-11c} the three connected regions
together form a fundamental domain of the McKay quiver. So if
$\sinksource_{\tilde{\mathcal{M}},E}$ has three sources $P_1,P_2,P_3$
as in Figure \ref{figure-24} the corresponding fundamental domain will
be strictly contained in the fundamental domain defined by 
$\sinksource_{\tilde{\mathcal{M}},F}$ (contradiction). 
\end{proof}

\section{Proof of the main results}
\label{section-main-results}

\begin{proof}[Proof of Theorem \ref{theorem-images-are-pure-sheaves}]

Recall that $Y$ is $G$-$\hilb(\mathbb{C}^3)$, $\tilde{\mathcal{M}}$ is
the dual of the universal family of $G$-clusters on $Y$ and $\Psi:
D^G(\C^3) \rightarrow D(Y)$ is the integral transform defined by 
$\tilde{\mathcal{M}}$. We fix $\chi \in G^\vee$ and look at the
associated representation $\mckquiv_{\tilde{\mathcal{M}}}$ and
its subrepresentation
$\hex(\chi^{-1})_{\tilde{\mathcal{M}}}$. Then the total complex $T^\bullet$
of the skew-commutative cube corresponding to
$\hex(\chi^{-1})_{\tilde{\mathcal{M}}}$ gives $\Psi(\O_0 \otimes \chi)
\in D(Y)$ (Proposition
\ref{prps-image-of-chi-twisted-O0-is-hex-chi-inverse}). 

If $\tilde{\mathcal{M}} = \bigoplus_{\chi \in G^\vee} \mathcal{L}(-M'_\chi)$ then in the notation of Section \ref{section-cohomology-of-skew-commutative-cubes}, the skew-commutative cube $T^\bullet$ corresponding to $\hex(\chi^{-1})_{\tilde{\mathcal{M}}}$ is described as follows. For any subset $v$ of $\{1,2,3\}$ we have $\mathcal{L}_{v} = \mathcal{L}(-M'_{\kappa(v)\chi})$ where $\kappa(v) = \prod_{i \in v} \kappa(x_i)$. The divisors $D^i_v$ record where the maps $\alpha^i_v$ in the cube vanish (see Figure \ref{figure-26}). 
\begin{figure}[h]
\begin{center}
\includegraphics[scale=0.17]{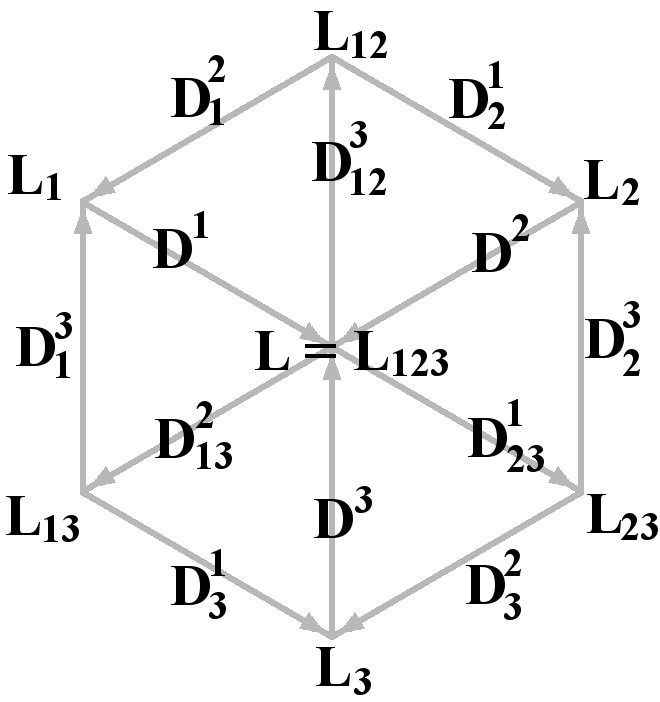}
\caption{\label{figure-26} The skew-commutative cube corresponding to $\hex(\chi^{-1})$. } 
\end{center}
\end{figure}

By Lemma \ref{lem:cube_cohomology} the only way $H^{-2}(T^\bullet) \ne
0$ is if $gcd(D_{23}^1, D_{13}^2, D_{12}^3) \neq 0$. This is only the
case if $\chi = \chi_0$ (Corollary
\ref{cor-gcluster-no-three-arrows}). On the other hand, if $\chi =
\chi_0$, then $\chi$ is a $(0,3)$-source in
$\sinksource_{\tilde{\mathcal{M}},E}$ for every irreducible
exceptional divisor $E$ on $Y$. So every $D^i_{jk}$ is equal to the
whole exceptional set $\Except(Y)$ while every $D^i_j$ and $D^i$ is
zero. By Lemma \ref{lem:cube_cohomology} this means $H^{-2}(T^\bullet)
= \mathcal{L}_{123}(\Except(Y))\otimes \mathcal{O}_{\Except(Y)}$ and
$H^{-1}(T^\bullet) = H^0(T^\bullet) = 0$. Since $M'_{\chi_0} = 0$ we
have $\mathcal{L}_{123} = \mathcal{L}(-M'_{\chi_0}) = \mathcal{O}_Y$
and $\Psi(\O_0 \otimes \chi) = \mathcal{O}(\Except(Y)) \otimes
\mathcal{O}_{\Except(Y)}$. 

For the rest of the argument we assume $\chi \ne \chi_0$ in which case
$T^\bullet$ is supported in degrees $0$ and $-1$. We will need two
facts. One is that  
$$\Ext^0(T^\bullet, T^\bullet) \cong \Ext^0_G(\O_0 \otimes \chi, \O_0 \otimes \chi) \cong \C$$
since $\Psi$ is an equivalence. The other fact is the standard exact triangle
\begin{equation}\label{std-triangle}
\sW[1] \rightarrow T^\bullet \rightarrow \O_{D^1 \cap D^2 \cap D^3}(\sL)
\end{equation}
where $\sW = H^{-1}(T^\bullet)$. 

Suppose $H^0(T^\bullet) \ne 0$. By Lemma \ref{lem:cube_cohomology} this means $D^1 \cap D^2 \cap D^3 \ne 0$ where each $D^i \in \Except(Y)$. Let $L$ be an irreducible component of $D^1 \cap D^2 \cap D^3 \ne 0$.  There are three possibilities: 
\begin{enumerate}
\item \label{item-isolated-points} $L$ is an isolated point 
\item \label{item-codim-1}  $L$ is an intersection $E_1 \cap E_2$ of two irreducible exceptional divisors 
\item \label{item-divisors} $L$ is an exceptional divisor
\end{enumerate}
Below we show that case (\ref{item-isolated-points}) cannot occur while in the other cases $\ext^1(\mathcal{O}_L \otimes \mathcal{L}', \sW [1]) = 0$ for any line bundle $\mathcal{L}'$. This implies $\ext^1(\mathcal{O}_{D^1 \cap D^2 \cap D^3} \otimes \mathcal{L}, \sW[1]) = 0$ because we can repeatedly apply an exact sequence of the form 
$$\O_{L_2} \otimes \mathcal{I}_{L_1} \rightarrow \O_{L_1 \cup L_2} \rightarrow \O_{L_1}$$
where $L_1$ and $L_2$ are components of $D^1 \cap D^2 \cap D^3$. Notice that if we break up $D^1 \cap D^2 \cap D^3$ into components in the right order then we can ensure that $\O_{L_2} \otimes \mathcal{I}_{L_1}$ is a line bundle on $L_2$ so that by induction $\ext^1(\O_{L_2} \otimes \mathcal{I}_{L_1}, \sW[1]) = \ext^1(\O_{L_1}, \sW[1]) = 0$ which implies $\ext^1(\O_{L_1 \cup L_2}, \sW[1]) = 0$. 

The fact that $\ext^1(\O_{D^1 \cap D^2 \cap D^3}(\sL), \sW[1]) = 0$
implies $T^\bullet \cong \sW[1] \oplus \O_{D^1 \cap D^2 \cap D^3}
(\sL)$. This means $\sW = 0$ for otherwise $\mbox{dim}
\Ext^0(T^\bullet, T^\bullet) \ge 2$. We conclude that if
$H^0(T^\bullet) \neq 0$ then $H^1(T^\bullet) = 0$. As by
assumption $H^i(T^\bullet) \neq 0$ only for $i = 0, -1$, this
completes the proof of Theorem \ref{theorem-images-are-pure-sheaves}.  

Case (\ref{item-isolated-points}). The only way $D^1 \cap D^2 \cap D^3$ can contain an isolated point is if each $D^i$ contains an (irreducible) exceptional divisor $E_i$ such that $E_1 \cap E_2 \cap E_3$ is a point $p$. Then $H^0(T^\bullet)$ is the direct sum of $\O_p$ and some sheaf. If $p \not\in \supp(\hmone)$ then $T^\bullet$ breaks up as the direct sum of $\O_p$ and some complex which contradicts the fact that $\Ext^0(T^\bullet, T^\bullet) \cong \C$. Thus $p$ belongs to the support of $\hmone$. 

By Lemma \ref{lem:cube_cohomology} the support of $\hmone$ is contained in the union of $\gcd(D_1^2,D_2^1)$, $\gcd(D_1^3,D_3^1)$ and $\gcd(D_2^3,D_3^2)$. Without loss of generality suppose $p \in \gcd(D_1^2,D_2^1)$. Since $p$ belongs to precisely three exceptional divisors ($E_1,E_2$ and $E_3$), one of them must belong to $\gcd(D_1^2,D_2^1)$. It cannot be $E_1$ or $E_2$ since $D^1 + D_1^2 = D^2 + D_2^1$ contain $E_1 + E_2$. Therefore $E_3$ belongs to $\gcd(D_1^2,D_2^1)$. This is impossible by Proposition \ref{prps-no-comets}. 

Case (\ref{item-codim-1}). 
Since $L = E_1 \cap E_2$ is a component of $D^1 \cap D^2 \cap D^3$ we can assume (without loss of generality) that $D^1$ and $D^2$ contain $E_1$ while $D^3$ contains $E_2$. To show $\Ext^1(\O_L \otimes \mathcal{L}', \sW[1]) = 0$ it suffices to show that $\Ext^2(\O_L \otimes \mathcal{L}', \O_C) = 0$ for any irreducible component $C$ of $\supp(\sW)$.  Denote by $i: C \rightarrow Y = G$-$\hilb (\mathbb{C}^3)$ the inclusion.  Then by adjunction $$\Ext^2(\O_L \otimes \mathcal{L}', \O_C) \cong \Ext^2(i^* \O_L \otimes \mathcal{L}', \O_C)$$
so if $L \cap C = \emptyset$ then $i^* \O_L = 0$ and we are done. So we can assume $L \cap C$ is non-empty. 

Suppose $C \subset Y$ has codimension one (i.e. is an irreducible exceptional divisor). If $L \subset C$ then by Lemma \ref{lem:cube_cohomology} some $\gcd(D_i^j,D_j^i)$ must contain either $E_1$ or $E_2$. None of them can contain $E_1$ since $D^1$ and $D^2$ contain $E_1$ and so some $\gcd(D_i^j, D_j^i)$ must contain $E_2$. Since $D^3$ contains $E_2$ this means it must be $\gcd(D_1^2, D_2^1)$ that contains $E_2$. By Corollary \ref{cor-opposite-oriented-charges-do-not-pass-through-each-other} this implies $E_1 \cap E_2 = \emptyset$ (contradiction). 

Suppose $L \cap C = E_1 \cap E_2 \cap C$ is zero dimensional (i.e. a point). By \ref{lem:cube_cohomology} some $\gcd(D_i^j,D_j^i)$ must contain $C$. By Corollary \ref{cor-opposite-oriented-charges-do-not-pass-through-each-other} $\gcd(D_1^2,D_2^1)$ cannot contain $C$ since this would mean $C \cap E_1 = \emptyset$. So, without loss of generality, we can assume $\gcd(D_2^3,D_3^2)$ contains $C$. Since $C \in \mbox{supp} H^{-1}(T^\bullet)$, using Lemma \ref{lem:cube_cohomology}, either $D^1$ contains $C$ or $\lcm(D_1^2,D_1^3)$ contain $C$.

If $D^1$ contains $C$ then by Proposition \ref{prps-no-comets} $E_1 \cap E_2 \cap C = \emptyset$ (contradiction). So we assume $\lcm(D_1^2,D_1^3)$ contains $C$. Since neither $D^1$ nor $D^2$ contain $C$ this means either both $D_1^2$ and $D_2^1$ cannot $C$ or neither contain $C$. They cannot both contain $C$ because if they did then $C \cap E_1$ would be empty by Corollary $\ref{cor-opposite-oriented-charges-do-not-pass-through-each-other}$. Thus $D_1^2$ does not contain $C$ and since $\lcm(D_1^2,D_1^3)$ contains $C$ this means $D_1^3$ must contains $C$. Finally, $D^3$ does not contain $C$ so this means $D_3^1$ must contain $C$. We conclude that $D_2^3, D_3^2, D_1^3$ and $D_3^1$ all contain $C$. By Proposition \ref{prps-no-shields} this means $C \cap E_1 = \emptyset$ (contradiction). 

Finally, suppose $C \subset Y$ has codimension two so that $C = C_1 \cap C_2$ where $C_1, C_2 \in \Except(Y)$. The fact $L \cap C$ is non-empty means that (without loss of generality) $C_1$ equals $E_1$ or $E_2$. Notice that $\gcd(D_i^j,D_j^i)$ cannot contain $C_1$ (to deal with $\gcd(D_1^2,D_2^1)$ we use \ref{cor-opposite-oriented-charges-do-not-pass-through-each-other}). It follows by Lemma \ref{lem:cube_cohomology} that $C_2$ must be contained in some $\gcd(D_i^j,D_j^i)$. By \ref{cor-opposite-oriented-charges-do-not-pass-through-each-other} $C_2$ cannot be contained in $\gcd(D_1^2,D_2^1)$ so, without loss of generality, let's suppose it is contined in $\gcd(D_2^3,D_3^2)$. By Lemma \ref{lem:cube_cohomology} $C_1$ must be contained in 
$$D^1 + \lcm(D_1^2,D_1^3) - \tilde{D_2^3} - D^2$$
since $C_1 \cap C_2 \in \mbox{supp} H^{-1}(T^\bullet)$. Since $E_1$ belongs to $D^2$ this means $E_1$ cannot belong to $D^1 + \lcm(D_1^2,D_1^3) - \tilde{D_2^3} - D^2$. So $C_1 = E_2$. But then $E_2$ belongs to $\tilde{D_2^3} + D^2$ so again $C_1 = E_2$ cannot belong to $D^1 + \lcm(D_1^2,D_1^3) - \tilde{D_2^3} - D^2$ (contradiction). 

Case (\ref{item-divisors}). 
The only way $\Ext^1(\O_L(\sL), \sW[1])$ is non-zero is if $L$ intersects (without loss of generality) $\gcd(D_1^2,D_2^1)$. Then $D^1$ and $D^2$ would contain $L$ while $\gcd(D_1^2,D_2^1)$ would contain an exceptional divisor $F$ such that $L \cap F \ne \emptyset$. This is impossible by Corollary \ref{cor-opposite-oriented-charges-do-not-pass-through-each-other}. 
\end{proof}

\section{A worked example}
\label{section-worked-example}

In this section we give a worked example of how to compute images
$\Psi(\mathcal{O}_0 \otimes \chi)$ for a given abelian subgroup 
of $\gsl_3(\mathbb{C})$. 

\subsection{The group $G$ and its McKay quiver} 
\label{section-example-the-group}

We set $G$ to be the group $\frac{1}{13}(1,5,7)$. That is, the image 
in $\gsl_3(\mathbb{C})$ of group $\mu_{13}$ of $13$th roots of unity 
under the embedding
$\xi \mapsto 
\left( \begin{smallmatrix}
\xi^1 & & \\
& \xi^5 & \\
& & \xi^7
\end{smallmatrix} \right)
$. We denote by $\chi_{i}$ the character of $G$ induced by $\xi \mapsto \xi^i$.
Then $\rho(x_1) = \chi_1$, $\rho(x_2) = \chi_5$, $\rho(x_3) = \chi_7$, 
therefore $\kappa(x_1) = \chi_{12}$, $\kappa(x_2) = \chi_{8}$ and
$\kappa(x_3) = \chi_6$.

\begin{figure}[h]
\begin{center}
\includegraphics[scale=0.35]{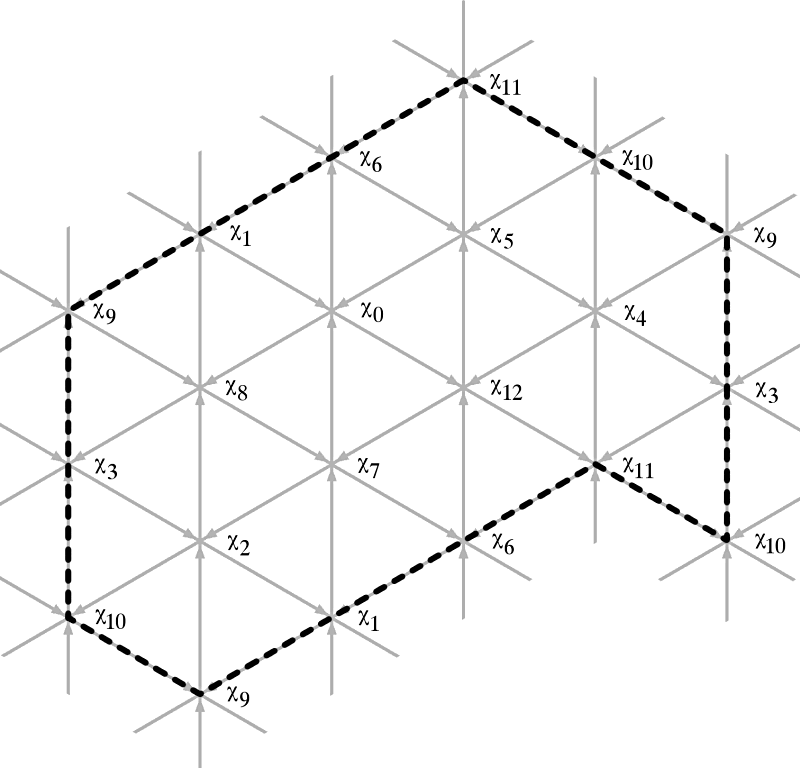}
\end{center}
\caption{\label{figure-2} A fundamental domain of the McKay quiver of $G = \frac{1}{13}(1,5,7)$.}
\end{figure}

As explained in Section
\ref{section-the-mckay-quiver-of-G-and-its-planar-embedding}, for 
each $\chi_i \in G^\vee$ the McKay quiver $\mckquiv$ has $3$ arrows 
emerging from the vertex $\chi_i$ and going to vertices $\kappa(x_j)
\chi_i$ for $j = 1,2,3$. E.g. from the vertex $\chi_3$ there are 
arrows going to vertices $\kappa(x_1) \chi_{3} = \chi_{2}$,
$\kappa(x_2) \chi_{3} = \chi_{11}$ and $\kappa(x_3) \chi_{3} =
\chi_{9}$. Choosing a particular fundamental domain in the planar embedding 
$\phi_H$ of the universal cover $U$ of $\mckquiv$ into
$\mathbb{R}^2$, we can depict $\mckquiv$ as shown in 
Figure \ref{figure-2}. Gluing together the edges of the diagram 
appropriately one obtains the real $2$-dimensional torus $T_G$ as 
tessellated by the McKay quiver. 

\subsection{The resolution $G$-Hilb$(\mathbb{C}^3)$ and the universal family $\mathcal{M}$ of $G$-clusters}
\label{section-example-the-resolution}
Let $Y = G$-$\hilb(\mathbb{C}^3)$.  For any abelian subgroup $G
\subset \gsl_3(\C)$ $Y$ is a toric variety whose toric fan can be
computed as described in
\cite{Craw-AnexplicitconstructionoftheMcKaycorrespondenceforAHilbC3},
Section 2. In our case we obtain the
$3$-dimensional fan in $\mathbb{R}^3$ as depicted in Figure
\ref{figure-03}. As usual, we drew the junior simplex $\Delta$ -  
the two-dimensional section of the
fan cut out by the plane $l_1 + l_2 + l_3 = 1$.
The generators $e_i$ of the one-dimensional cones of the fan all lie
in $\Delta$ and have the following coordinates:
\begin{align} \label{eqn-e_i}
\begin{tabular}{l l l}
$e_1 = (1,0,0)$  &  $e_2 = (0,1,0)$  &  $e_3 = (0,0,1)$ \\
$e_4 = \frac{1}{13}(1,5,7)$  &  $e_5 = \frac{1}{13}(2,10,1)$ 
& $e_6 = \frac{1}{13}(3,2,8)$ \\
$e_7 = \frac{1}{13}(4,7,2)$  &  $e_8 = \frac{1}{13}(6,4,3)$ 
& $e_9 = \frac{1}{13}(8,1,4)$ 
\end{tabular} 
\end{align} 

To each $k$-dimensional cone $\left<e_{i_1}, \dots, e_{i_k}\right>$ 
in the toric fan of $Y$ there corresponds a $(3-k)$-dimensional torus
orbit in $Y$. We denote this orbit by $S_{i_1,\dots,i_k}$ and its closure 
by $E_{i_1, \dots, i_k}$. The set of all irreducible 
exceptional divisors on $Y$ is then $\{E_4,E_5,E_6,E_7,E_8,E_9\}$. 

On $Y = G$-$\hilb(\mathbb{C}^3)$ the universal family 
$\mathcal{M} = \bigoplus \mathcal{L}(-M_\chi)$ of $G$-clusters 
$Y$ is the \it maximal shift $\gnat$-family \rm (cf. 
\cite{Logvinenko-Natural-G-Constellation-Families}, Prop. 3.17).
Each $G$-Weil divisor $M_\chi$ is of the form 
$\sum_{E_i \in \Except(Y)} m_{\chi,i} E_i$ 
(\cite{Logvinenko-Natural-G-Constellation-Families}, Prop. 3.15). 
As shown in
\cite{Logvinenko-Families-of-G-constellations-over-resolutions-of-quotient-singularities},
Example 4.21 the coefficients $m_{\chi,i} \in \mathbb{Q}$ can be computed via 
the formula
\begin{align} \label{eqn-maximal-shift-set}
        m_{\chi,i} = 
        \inf \{ e_i(m) \; | \; m \in \mathbb{Z}^3_{\geq 0} \text{ and } \rho(x^m) = \chi \},
\end{align} 

To compute $m_{\chi, i}$ it suffices to test $e_i(m_1,m_2,m_3)$ for
minimality only on $(m_1, m_2, m_3) \in \mathbb{Z}^3$ for which $0
\leq m_1, m_2, m_3 \leq 13$. This is because $x_1^{13}, x_2^{13}$ and
$x_3^{13}$ are $G$-invariant, while $e_i(m_1 + 13, m_2, m_3) \geq
e_i(m_1, m_2, m_3)$ and similarly for the other two. It then becomes a
simple computer computation.  

For example, the minimal value of $e_7(-)$ on $(m_1,m_2,m_3)$ with
$\rho(x_1^{m_1}x_2^{m_2}x_3^{m_3}) = \chi_{10}$
and $0 \leq m_1, m_2, m_3 \leq 13$ can be verified to be $1
\frac{1}{13}$. E.g. $\rho(x_1^3 x_3) = \chi_1^3 \chi_7 = \chi_{10}$ 
and
\begin{align*}
e_7(x_1^3 x_3) = e_7(3,0,1) = \frac{4}{13} \cdot 3 + \frac{7}{13} \cdot 0 + \frac{2}{13} \cdot 1 = 1 \frac{1}{13}.
\end{align*}
Hence $m_{\chi_{10}, 7} = 1 \frac{1}{13}$.  

Repeating this for other $\chi$ and $e_i$, we compute $m_{\chi,i}$ to be:
\begin{align} \label{eqn-maximal-shift-family}
\text{
\begin{tiny}
\begin{tabular}{|r|r|r|r|r|r|r||r|r|r|r|r|r|r|}
\hline
$\chi \setminus i $ & $ 4 $ & $ 5 $ & $ 6 $ & $ 7 $ & $ 8 $ & $ 9 $ & 
$\chi \setminus i $ & $ 4 $ & $ 5 $ & $ 6 $ & $ 7 $ & $ 8 $ & $ 9 $ \\ 
\hline
& & & & & & & & & & & & & \\
$\chi_{0}$ & $ 0 $ & $ 0 $ & $ 0 $ & $ 0 $ & $ 0 $ & $ 0 $ &
$\chi_{7}$ & $  \frac{7}{13} $ & $  \frac{1}{13} $ & $ \frac{8}{13} $ & 
$ \frac{2}{13} $ & $  \frac{3}{13} $ & $ \frac{4}{13}$ \\
& & & & & & & & & & & & & \\
$\chi_{1}$ & $  \frac{1}{13} $ & $  \frac{2}{13} $ & $ 
 \frac{3}{13} $ & $  \frac{4}{13} $ & $  \frac{6}{13} $ & $ 
\frac{8}{13}$ &
$\chi_{8}$ & $  \frac{8}{13} $ & $  \frac{3}{13} $ & $ 
 \frac{11}{13} $ & $  \frac{6}{13} $ & $  \frac{9}{13} $ & $ 
 \frac{12}{13}$ \\
& & & & & & & & & & & & & \\
$\chi_{2}$ & $  \frac{2}{13} $ & $  \frac{4}{13} $ & $ \frac{6}{13} $ & 
$ \frac{8}{13} $ & $  \frac{12}{13} $ & $ \frac{3}{13}$ &
$\chi_{9}$ & $  \frac{9}{13} $ & $  \frac{5}{13} $ & $ 1 \frac{1}{13} $ & 
$ \frac{10}{13} $ & $  1 \frac{2}{13} $ & $ \frac{7}{13}$ \\
& & & & & & & & & & & & & \\
$\chi_{3}$ & $  \frac{3}{13} $ & $  \frac{6}{13} $ & $ \frac{9}{13} $ & 
$ \frac{12}{13} $ & $  1 \frac{5}{13} $ & $ \frac{11}{13}$ &
$\chi_{10}$ & $  \frac{10}{13} $ & $  \frac{7}{13} $ & $ \frac{4}{13} $ & 
$ 1 \frac{1}{13} $ & $  \frac{8}{13} $ & $ \frac{2}{13}$ \\
& & & & & & & & & & & & & \\
$\chi_{4}$ & $  \frac{4}{13} $ & $  \frac{8}{13} $ & $ \frac{12}{13} $ & 
$ 1 \frac{3}{13} $ & $ \frac{11}{13} $ & $ \frac{6}{13}$ &
$\chi_{11}$ & $ \frac{11}{13} $ & $  \frac{9}{13} $ & $ \frac{7}{13} $ & 
$ 1 \frac{5}{13} $ & $  1 \frac{1}{13} $ & $ \frac{10}{13}$ \\
& & & & & & & & & & & & & \\
$\chi_{5}$ & $  \frac{5}{13} $ & $  \frac{10}{13} $ & $ \frac{2}{13} $ & 
$ \frac{7}{13} $ & $  \frac{4}{13} $ & $ \frac{1}{13}$ & 
$\chi_{12}$ & $  \frac{12}{13} $ & $  \frac{11}{13} $ & $ \frac{10}{13} $ & 
$ \frac{9}{13} $ & $  \frac{7}{13} $ & $ \frac{5}{13}$ \\
& & & & & & & & & & & & & \\
$\chi_{6}$ & $  \frac{6}{13} $ & $  \frac{12}{13} $ & $ \frac{5}{13} $ & 
$ \frac{11}{13} $ & $  \frac{10}{13} $ & $ \frac{9}{13}$ &
& & & & & & \\
& & & & & & & & & & & & & \\
\hline 
\end{tabular}
\end{tiny}
}
\end{align}

\subsection{Divisors of the associated representation of the dual family $\tilde{\mathcal{M}}$}

By Proposition \ref{prps-the-dual-family} the dual $\tilde{\mathcal{M}}$ of
$\mathcal{M} = \bigoplus(-M_\chi)$ is the family 
$\oplus \mathcal{L}(-M'_\chi)$ where $M'_\chi = - M_{\chi^{-1}}$. So to
compute the divisors of zeroes $B_{\chi,i}$ of the maps
$\alpha_{\chi,i}$ in the 
associated representation of the McKay quiver $\mckquiv_{\tilde{\mathcal{M}}}$ 
we must set $D_\chi = - M_{\chi^{-1}}$ in formula
$\ref{eqn-zero-divisors-formula}$. Thus we obtain the formula 
\begin{align} \label{eqn-zero-divisors-of-the-dual-family}
B_{\chi,i} = M_{\chi \rho(x_i)^{-1}} + (x_i) - M_{\chi}
\end{align}

The principal $G$-Weil divisors $(x_i)$ of the basic monomials can be 
computed via the formula $(x_i) = \sum e_i(x_i) E_i$ (cf. 
\cite{Logvinenko-Families-of-G-constellations-over-resolutions-of-quotient-singularities}, Prop. 3.2). For an example, $e_4 = \frac{1}{13}(1,5,7)$
so $e_4(x_1) = \frac{1}{13}$, $e_4(x_2) = \frac{5}{13}$ and $e_4(x_3)
= \frac{7}{13}$ and obtain: 
\begin{footnotesize}
\begin{align*} 
(x_1) = E_1 + \frac{1}{13} E_4 + \frac{2}{13} E_5 + \frac{3}{13} E_6 + 
\frac{4}{13} E_7 + \frac{6}{13} E_8 + \frac{8}{13} E_9 \\
(x_2) = E_2 + \frac{5}{13} E_4 + \frac{10}{13} E_5 + \frac{2}{13} E_6 + 
\frac{7}{13} E_7 + \frac{4}{13} E_8 + \frac{1}{13} E_9 \\
(x_3) = E_3 + \frac{7}{13} E_4 + \frac{1}{13} E_5 + \frac{8}{13} E_6 + 
\frac{2}{13} E_7 + \frac{3}{13} E_8 + \frac{4}{13} E_9
\end{align*}
\end{footnotesize}

Let us now compute, for example, zero divisor $B_{\chi_5,1}$. Using
formula \eqref{eqn-zero-divisors-of-the-dual-family} we have
\begin{align*}
& B_{\chi_5,1} = M_{\chi_{5}\chi_{1}^{-1}} + (x_1) +  M_{\chi_{5}} =  M_{\chi_4} + (x_1) - M_{\chi_{5}} \\
& 
\setlength{\extrarowheight}{0.2cm}
\begin{array}{r r r r r r r r r r r r r r r r}
= & & & & \frac{4}{13} E_4 & + & \frac{8}{13} E_5 & + & \frac{12}{13} E_6 & 
+ & 1 \frac{3}{13}E_7 & + & \frac{11}{13} E_8 & + & \frac{6}{13}E_9 & \\
& + & E_1 & + & 
\frac{1}{13} E_4 & + & \frac{2}{13} E_5 & + & \frac{3}{13} E_6 & 
+ & \frac{4}{13} E_7 & + & \frac{6}{13} E_8 & + & \frac{8}{13}E_9 & \\
& & & - & 
\frac{5}{13} E_4 & - & \frac{10}{13} E_5 & - & \frac{2}{13} E_6 &
- & \frac{7}{13} E_7 & - & \frac{4}{13} E_8 & - &
\frac{1}{13} E_9 & \\
= & & E_1 & & & & & + & E_6 & + & E_7 & + & E_8 & + & E_9 &
\end{array}
\end{align*}

Computing similarly the remaining zero divisors $B_{\chi,k}$ of $\mckquiv_{\tilde{\mathcal{M}}}$ we obtain:
\begin{align} \label{eqn-divisors-of-zeroes}
\text{
\begin{footnotesize}
\begin{tabular}{l l}
$ B_{\chi_{0}, 1} = E_1 + E_4 + E_5 + E_6 + E_7 + E_8 + E_9 $ & 
$ B_{\chi_{7}, 1} = E_1 + E_5 + E_7 + E_8 + E_9 $ \\
$ B_{\chi_{0}, 2} = E_2 + E_4 + E_5 + E_6 + E_7 + E_8 + E_9 $ & 
$ B_{\chi_{7}, 2} = E_2 + E_5 + E_7 + E_8 $ \\
$ B_{\chi_{0}, 3} = E_3 + E_4 + E_5 + E_6 + E_7 + E_8 + E_9 $ & 
$ B_{\chi_{7}, 3} = E_3 $ \\
$ B_{\chi_{1}, 1} = E_1 $ & 
$ B_{\chi_{8}, 1} = E_1 $ \\
$ B_{\chi_{1}, 2} = E_2 + E_4 + E_5 + E_6 + E_7 + E_8 $ &
$ B_{\chi_{8}, 2} = E_2 + E_5 + E_7 + E_8 $ \\
$ B_{\chi_{1}, 3} = E_3 + E_4 + E_6 $ &
$ B_{\chi_{8}, 3} = E_3 $ \\
$ B_{\chi_{2}, 1} = E_1 + E_9 $ &
$ B_{\chi_{9}, 1} = E_1 + E_9 $ \\
$ B_{\chi_{2}, 2} = E_2 + E_4 + E_5 + E_7 $ &
$ B_{\chi_{9}, 2} = E_2 + E_5 + E_7 $ \\
$ B_{\chi_{2}, 3} = E_3 + E_4 + E_6 + E_9 $ &
$ B_{\chi_{9}, 3} = E_3 $ \\
$ B_{\chi_{3}, 1} = E_1 $ &
$ B_{\chi_{10}, 1} = E_1 +  E_6 + E_8 + E_9 $ \\
$ B_{\chi_{3}, 2} = E_2 + E_4 + E_5 + E_7 $ &
$ B_{\chi_{10}, 2} = E_2 + E_5 $ \\
$ B_{\chi_{3}, 3} = E_3 + E_4 + E_6 $ &
$ B_{\chi_{10}, 3} = E_3 +  E_6 + E_8 + E_9 $ \\
$ B_{\chi_{4}, 1} = E_1 + E_8 + E_9 $ &
$ B_{\chi_{11}, 1} = E_1 $ \\
$ B_{\chi_{4}, 2} = E_2 + E_4 + E_5 $ &
$ B_{\chi_{11}, 2} = E_2 + E_5 $ \\
$ B_{\chi_{4}, 3} = E_3 + E_4 $ &
$ B_{\chi_{11}, 3} = E_3 + E_6 $ \\
$ B_{\chi_{5}, 1} = E_1 + E_6 + E_7 + E_8 + E_9 $ &
$ B_{\chi_{12}, 1} = E_1 + E_7 + E_8 + E_9 $ \\
$ B_{\chi_{5}, 2} = E_2  $ &
$ B_{\chi_{12}, 2} = E_2 $ \\
$ B_{\chi_{5}, 3} = E_3 + E_4 + E_6 + E_7 + E_8 + E_9 $ &
$ B_{\chi_{12}, 3} = E_3 $ \\
$ B_{\chi_{6}, 1} = E_1 $ & \\
$ B_{\chi_{6}, 2} = E_2 $ & \\
$ B_{\chi_{6}, 3} = E_3 + E_4 + E_6 $ & \\
\end{tabular}
\end{footnotesize}
}
\end{align}

\subsection{Skew-commutative cubes corresponding to $\Psi(\mathcal{O}_0 \otimes \chi)$}

Now, for each $\chi \in G^\vee$, we compute the skew-commutative cube
of line bundles $\hex(\chi^{-1})_{\tilde{\mathcal{M}}}$ (see Section 
\ref{section-associated-representation}). The total complex
of the cube gives then $\Psi(\mathcal{O}_0 \otimes \chi)$ in $D(Y)$.

Recall that the subquiver $\hex(\chi^{-1})$ 
was defined by sending vertex $v \subset \{1,2,3\}$ of the 
cube quiver to vertex $\kappa(v)^{-1} \chi^{-1}$ in $\mckquiv$. 
Therefore at vertex $v$ of the cube of line bundles corresponding 
to $\hex(\chi^{-1})_{\tilde{\mathcal{M}}}$ we get the summand of
$\tilde{\mathcal{M}}=\bigoplus \mathcal{L}(-M'_{\chi})$ which 
$G$ acts on by $\kappa(v)^{-1} \chi^{-1}$. So $\mathcal{L}_{v} = 
\mathcal{L}(-M'_{\kappa(v)\chi}) =
\mathcal{L}(M_{(\kappa(v)\chi)^{-1}})$. 
On each arrow in the cube we mark its divisor of the zeroes, using 
the data of \eqref{eqn-divisors-of-zeroes}: 
$D^i_{jk} = B_{\chi^{-1},i}$, $D^i_j = B_{\kappa(x_k)\chi^{-1}, j}$
and $D^i = B_{\kappa(x_i)^{-1}\chi^{-1}, i}$ for any 
$i \neq j \neq k$. We use a following shorthand: $E_{456} = E_4 + E_5 + E_6$ 
et cetera. 

\begin{small}
\begin{align*} 
\Psi(\mathcal{O}_0 \otimes \chi_0):\quad
\xymatrix{
& \mathcal{L}(M_{\chi_{12}}) \ar"2,5"^<<{E_3} \ar"3,5"_<<<{E_2}
& & & 
\mathcal{L}(M_{\chi_1}) \ar"2,6"^{E_1} 
& \\
\mathcal{L}(M_{\chi_0}) \ar"1,2"^{E_{1456789}} \ar"2,2"^{E_{2456789}}
\ar"3,2"_{E_{3456789}} 
& \mathcal{L}(M_{\chi_8}) \ar"1,5"^>>>>{E_3} \ar"3,5"_>>>>>{E_1}
& & & \mathcal{L}(M_{\chi_5}) \ar"2,6"^{E_2} 
& \mathcal{L}(M_{\chi_0}) \\
& \mathcal{L}(M_{\chi_6}) \ar"1,5"^<<<{E_2} \ar"2,5"_<<{E_1} 
& & & \mathcal{L}(M_{\chi_{7}}) \ar"2,6"_{E_3} 
&
}
\end{align*}
\begin{align*} 
\Psi(\mathcal{O}_0 \otimes \chi_1):\quad
\xymatrix{
& \mathcal{L}(M_{\chi_{11}}) \ar"2,5"^<<{E_{36}} \ar"3,5"_<<<{E_{25}}
& & & 
\mathcal{L}(M_{\chi_0}) \ar"2,6"^{E_{1456789}} 
& \\
\mathcal{L}(M_{\chi_{12}}) \ar"1,2"^{E_{1789}} \ar"2,2"^{E_{2}}
\ar"3,2"_{E_{3}} 
& \mathcal{L}(M_{\chi_7}) \ar"1,5"^>>>>{E_{3}} \ar"3,5"_>>>>>{E_{15789}}
& & & \mathcal{L}(M_{\chi_4}) \ar"2,6"^{E_{245}} 
& \mathcal{L}(M_{\chi_{12}}) \\
& \mathcal{L}(M_{\chi_5}) \ar"1,5"^<<<{E_{2}} \ar"2,5"_<<{E_{16789}} 
& & & \mathcal{L}(M_{\chi_{6}}) \ar"2,6"_{E_{346}} 
&
}
\end{align*}
\begin{align*} 
\Psi(\mathcal{O}_0 \otimes \chi_2):\quad
\xymatrix{
& \mathcal{L}(M_{\chi_{10}}) \ar"2,5"^<<{E_{3689}} \ar"3,5"_<<<{E_{25}}
& & & 
\mathcal{L}(M_{\chi_{12}}) \ar"2,6"^{E_{1789}} 
& \\
\mathcal{L}(M_{\chi_{11}}) \ar"1,2"^{E_{1}} \ar"2,2"^{E_{25}}
\ar"3,2"_{E_{36}} 
& \mathcal{L}(M_{\chi_6}) \ar"1,5"^>>>>{E_{346}} \ar"3,5"_>>>>>{E_{1}}
& & & \mathcal{L}(M_{\chi_3}) \ar"2,6"^{E_{2457}} 
& \mathcal{L}(M_{\chi_{11}}) \\
& \mathcal{L}(M_{\chi_4}) \ar"1,5"^<<<{E_{245}} \ar"2,5"_<<{E_{189}} 
& & & \mathcal{L}(M_{\chi_{5}}) \ar"2,6"_{E_{346789}} 
&
}
\end{align*}
\begin{align*} 
\Psi(\mathcal{O}_0 \otimes \chi_3):\quad
\xymatrix{
& \mathcal{L}(M_{\chi_{9}}) \ar"2,5"^<<{E_{3}} \ar"3,5"_<<<{E_{257}}
& & & 
\mathcal{L}(M_{\chi_{11}}) \ar"2,6"^{E_{1}} 
& \\
\mathcal{L}(M_{\chi_{10}}) \ar"1,2"^{E_{1689}} \ar"2,2"^{E_{25}}
\ar"3,2"_{E_{3689}} 
& \mathcal{L}(M_{\chi_5}) \ar"1,5"^>>>>{E_{346789}} \ar"3,5"_>>>>>{E_{16789}}
& & & \mathcal{L}(M_{\chi_2}) \ar"2,6"^{E_{2457}} 
& \mathcal{L}(M_{\chi_{10}}) \\
& \mathcal{L}(M_{\chi_3}) \ar"1,5"^<<<{E_{2457}} \ar"2,5"_<<{E_{1}} 
& & & \mathcal{L}(M_{\chi_{4}}) \ar"2,6"_{E_{34}} 
&
}
\end{align*}
\begin{align*} 
\Psi(\mathcal{O}_0 \otimes \chi_4):\quad
\xymatrix{
& \mathcal{L}(M_{\chi_{8}}) \ar"2,5"^<<{E_{3}} \ar"3,5"_<<<{E_{2578}}
& & & 
\mathcal{L}(M_{\chi_{10}}) \ar"2,6"^{E_{1689}} 
& \\
\mathcal{L}(M_{\chi_{9}}) \ar"1,2"^{E_{19}} \ar"2,2"^{E_{257}}
\ar"3,2"_{E_{3}} 
& \mathcal{L}(M_{\chi_4}) \ar"1,5"^>>>>{E_{34}} \ar"3,5"_>>>>>{E_{189}}
& & & \mathcal{L}(M_{\chi_1}) \ar"2,6"^{E_{245678}} 
& \mathcal{L}(M_{\chi_{9}}) \\
& \mathcal{L}(M_{\chi_2}) \ar"1,5"^<<<{E_{2457}} \ar"2,5"_<<{E_{19}} 
& & & \mathcal{L}(M_{\chi_{3}}) \ar"2,6"_{E_{346}} 
&
}
\end{align*}
\begin{align*} 
\Psi(\mathcal{O}_0 \otimes \chi_5):\quad
\xymatrix{
& \mathcal{L}(M_{\chi_{7}}) \ar"2,5"^<<{E_{3}} \ar"3,5"_<<<{E_{2578}}
& & & 
\mathcal{L}(M_{\chi_{9}}) \ar"2,6"^{E_{19}} 
& \\
\mathcal{L}(M_{\chi_{8}}) \ar"1,2"^{E_{1}} \ar"2,2"^{E_{2578}}
\ar"3,2"_{E_{3}} 
& \mathcal{L}(M_{\chi_3}) \ar"1,5"^>>>>{E_{346}} \ar"3,5"_>>>>>{E_{1}}
& & & \mathcal{L}(M_{\chi_0}) \ar"2,6"^{E_{2456789}} 
& \mathcal{L}(M_{\chi_{8}}) \\
& \mathcal{L}(M_{\chi_1}) \ar"1,5"^<<<{E_{245678}} \ar"2,5"_<<{E_{1}} 
& & & \mathcal{L}(M_{\chi_{2}}) \ar"2,6"_{E_{3469}} 
&
}
\end{align*}
\begin{align*} 
\Psi(\mathcal{O}_0 \otimes \chi_6):\quad
\xymatrix{
& \mathcal{L}(M_{\chi_{6}}) \ar"2,5"^<<{E_{346}} \ar"3,5"_<<<{E_{2}}
& & & 
\mathcal{L}(M_{\chi_{8}}) \ar"2,6"^{E_{1}} 
& \\
\mathcal{L}(M_{\chi_{7}}) \ar"1,2"^{E_{15789}} \ar"2,2"^{E_{2578}}
\ar"3,2"_{E_{3}} 
& \mathcal{L}(M_{\chi_2}) \ar"1,5"^>>>>{E_{3469}} \ar"3,5"_>>>>>{E_{19}}
& & & \mathcal{L}(M_{\chi_{12}}) \ar"2,6"^{E_{2}} 
& \mathcal{L}(M_{\chi_{7}}) \\
& \mathcal{L}(M_{\chi_0}) \ar"1,5"^<<<{E_{2456789}} \ar"2,5"_<<<{E_{1456789}} 
& & & \mathcal{L}(M_{\chi_{1}}) \ar"2,6"_{E_{346}} 
&
}
\end{align*}
\begin{align*} 
\Psi(\mathcal{O}_0 \otimes \chi_7):\quad
\xymatrix{
& \mathcal{L}(M_{\chi_{5}}) \ar"2,5"^<<{E_{346789}} \ar"3,5"_<<<{E_{2}}
& & & 
\mathcal{L}(M_{\chi_{7}}) \ar"2,6"^{E_{15789}} 
& \\
\mathcal{L}(M_{\chi_{6}}) \ar"1,2"^{E_{1}} \ar"2,2"^{E_{2}}
\ar"3,2"_{E_{346}} 
& \mathcal{L}(M_{\chi_1}) \ar"1,5"^>>>>{E_{346}} \ar"3,5"_>>>>>{E_{1}}
& & & \mathcal{L}(M_{\chi_{11}}) \ar"2,6"^{E_{25}} 
& \mathcal{L}(M_{\chi_{6}}) \\
& \mathcal{L}(M_{\chi_{12}}) \ar"1,5"^<<<{E_{2}} \ar"2,5"_<<{E_{1789}} 
& & & \mathcal{L}(M_{\chi_{0}}) \ar"2,6"_{E_{3456789}} 
&
}
\end{align*}
\begin{align*} 
\Psi(\mathcal{O}_0 \otimes \chi_8):\quad
\xymatrix{
& \mathcal{L}(M_{\chi_{4}}) \ar"2,5"^<<{E_{34}} \ar"3,5"_<<<{E_{245}}
& & & 
\mathcal{L}(M_{\chi_{6}}) \ar"2,6"^{E_{1}} 
& \\
\mathcal{L}(M_{\chi_{5}}) \ar"1,2"^{E_{16789}} \ar"2,2"^{E_{2}}
\ar"3,2"_{E_{346789}} 
& \mathcal{L}(M_{\chi_0}) \ar"1,5"^>>>>{E_{3456789}}
\ar"3,5"_>>>>>{E_{1456789}}
& & & \mathcal{L}(M_{\chi_{10}}) \ar"2,6"^{E_{25}} 
& \mathcal{L}(M_{\chi_{5}}) \\
& \mathcal{L}(M_{\chi_{11}}) \ar"1,5"^<<<{E_{25}} \ar"2,5"_<<{E_{1}} 
& & & \mathcal{L}(M_{\chi_{12}}) \ar"2,6"_{E_{3}} 
&
}
\end{align*}
\begin{align*} 
\Psi(\mathcal{O}_0 \otimes \chi_9):\quad
\xymatrix{
& \mathcal{L}(M_{\chi_{3}}) \ar"2,5"^<<{E_{346}} \ar"3,5"_<<<{E_{2457}}
& & & 
\mathcal{L}(M_{\chi_{5}}) \ar"2,6"^{E_{16789}} 
& \\
\mathcal{L}(M_{\chi_{4}}) \ar"1,2"^{E_{189}} \ar"2,2"^{E_{245}}
\ar"3,2"_{E_{34}} 
& \mathcal{L}(M_{\chi_{12}}) \ar"1,5"^>>>>{E_{3}}
\ar"3,5"_>>>>>{E_{1789}}
& & & \mathcal{L}(M_{\chi_{9}}) \ar"2,6"^{E_{257}} 
& \mathcal{L}(M_{\chi_{4}}) \\
& \mathcal{L}(M_{\chi_{10}}) \ar"1,5"^<<<{E_{25}} \ar"2,5"_<<{E_{1689}} 
& & & \mathcal{L}(M_{\chi_{11}}) \ar"2,6"_{E_{36}} 
&
}
\end{align*}
\begin{align*}
\Psi(\mathcal{O}_0 \otimes \chi_{10}):\quad
\xymatrix{
& \mathcal{L}(M_{\chi_{2}}) \ar"2,5"^<<{E_{3469}} \ar"3,5"_<<<{E_{2457}}
& & & 
\mathcal{L}(M_{\chi_{4}}) \ar"2,6"^{E_{189}} 
& \\
\mathcal{L}(M_{\chi_{3}}) \ar"1,2"^{E_{1}} \ar"2,2"^{E_{2457}}
\ar"3,2"_{E_{346}} 
& \mathcal{L}(M_{\chi_{11}}) \ar"1,5"^>>>>{E_{36}} \ar"3,5"_>>>>>{E_{1}}
& & & \mathcal{L}(M_{\chi_{8}}) \ar"2,6"^{E_{2578}} 
& \mathcal{L}(M_{\chi_{3}}) \\
& \mathcal{L}(M_{\chi_{9}}) \ar"1,5"^<<<{E_{257}} \ar"2,5"_<<{E_{19}} 
& & & \mathcal{L}(M_{\chi_{10}}) \ar"2,6"_{E_{3689}} 
&
}
\end{align*}
\begin{align*}
\Psi(\mathcal{O}_0 \otimes \chi_{11}):\quad
\xymatrix{
& \mathcal{L}(M_{\chi_{1}}) \ar"2,5"^<<{E_{346}} \ar"3,5"_<<<{E_{245678}}
& & & 
\mathcal{L}(M_{\chi_{3}}) \ar"2,6"^{E_{1}} 
& \\
\mathcal{L}(M_{\chi_{2}}) \ar"1,2"^{E_{19}} \ar"2,2"^{E_{2457}}
\ar"3,2"_{E_{3469}} 
& \mathcal{L}(M_{\chi_{10}}) \ar"1,5"^>>>>{E_{3689}} \ar"3,5"_>>>>>{E_{1689}}
& & & \mathcal{L}(M_{\chi_{7}}) \ar"2,6"^{E_{2578}} 
& \mathcal{L}(M_{\chi_{2}}) \\
& \mathcal{L}(M_{\chi_{8}}) \ar"1,5"^<<<{E_{2578}} \ar"2,5"_<<{E_{1}} 
& & & \mathcal{L}(M_{\chi_{9}}) \ar"2,6"_{E_{3}} 
&
}
\end{align*}
\begin{align*}
\Psi(\mathcal{O}_0 \otimes \chi_{12}):\quad
\xymatrix{
& \mathcal{L}(M_{\chi_{0}}) \ar"2,5"^<<<{E_{3456789}} \ar"3,5"_<<<{E_{2456789}}
& & & 
\mathcal{L}(M_{\chi_{2}}) \ar"2,6"^{E_{19}} 
& \\
\mathcal{L}(M_{\chi_{1}}) \ar"1,2"^{E_{1}} \ar"2,2"^{E_{24578}}
\ar"3,2"_{E_{346}} 
& \mathcal{L}(M_{\chi_{9}}) \ar"1,5"^>>>>{E_{3}} \ar"3,5"_>>>>>{E_{19}}
& & & \mathcal{L}(M_{\chi_{6}}) \ar"2,6"^{E_{2}} 
& \mathcal{L}(M_{\chi_{1}}) \\
& \mathcal{L}(M_{\chi_{7}}) \ar"1,5"^<<<{E_{2578}} \ar"2,5"_<<{E_{15789}} 
& & & \mathcal{L}(M_{\chi_{8}}) \ar"2,6"_{E_{3}} 
&
}
\end{align*}
\end{small}

\subsection{Conclusion}
\label{section-exmpl-conclusion}

We now analyze each of the skew-commutative cubes
computed in the previous section and use Lemma {\ref{lem:cube_cohomology}
to write down supports of the cohomology sheaves of their total
complexes, i.e. of $\Psi(\mathcal{O}_0 \otimes \chi)$. 

For example, consider the skew-commutative cube for $\Psi(\mathcal{O}_0
\otimes \chi_1)$:
\begin{align*} 
\Psi(\mathcal{O}_0 \otimes \chi_1):\quad
\xymatrix{
& \mathcal{L}(M_{\chi_{11}}) \ar"2,5"^<<{E_{36}} \ar"3,5"_<<<{E_{25}}
& & & 
\mathcal{L}(M_{\chi_0}) \ar"2,6"^{E_{1456789}} 
& \\
\mathcal{L}(M_{\chi_{12}}) \ar"1,2"^{E_{1789}} \ar"2,2"^{E_{2}}
\ar"3,2"_{E_{3}} 
& \mathcal{L}(M_{\chi_7}) \ar"1,5"^>>>>{E_{3}} \ar"3,5"_>>>>>{E_{15789}}
& & & \mathcal{L}(M_{\chi_4}) \ar"2,6"^{E_{245}} 
& \mathcal{L}(M_{\chi_{12}}) \\
& \mathcal{L}(M_{\chi_5}) \ar"1,5"^<<<{E_{2}} \ar"2,5"_<<{E_{16789}} 
& & & \mathcal{L}(M_{\chi_{6}}) \ar"2,6"_{E_{346}} 
&
}
\end{align*}

We see that $D^1 \cap D^2 \cap D^3 = E_4$, that $\gcd(D^1_2, D^2_1) 
= \gcd(D^2_3, D^3_2) = \gcd(D^1_3,D^3_1) = 0$ and that
$\gcd(D^1_{23},D^2_{13},D^3_{12}) = 0$. Therefore by Lemma
{\ref{lem:cube_cohomology} we have $\supp
H^0(\Psi(\mathcal{O}_0 \otimes \chi_1)) = E_4$ while $\supp
H^{-1}(\Psi(\mathcal{O}_0 \otimes \chi_1)) = \supp
H^{-2}(\Psi(\mathcal{O}_0 \otimes \chi_1)) = 0$. 

On the other hand, looking at the the skew-commutative cube for
$\Psi(\mathcal{O}_0 \otimes \chi_{11})$
\begin{align*}
\Psi(\mathcal{O}_0 \otimes \chi_{11}):\quad
\xymatrix{
& \mathcal{L}(M_{\chi_{1}}) \ar"2,5"^<<{E_{346}} \ar"3,5"_<<<{E_{245678}}
& & & 
\mathcal{L}(M_{\chi_{3}}) \ar"2,6"^{E_{1}} 
& \\
\mathcal{L}(M_{\chi_{2}}) \ar"1,2"^{E_{19}} \ar"2,2"^{E_{2457}}
\ar"3,2"_{E_{3469}} 
& \mathcal{L}(M_{\chi_{10}}) \ar"1,5"^>>>>{E_{3689}} \ar"3,5"_>>>>>{E_{1689}}
& & & \mathcal{L}(M_{\chi_{7}}) \ar"2,6"^{E_{2578}} 
& \mathcal{L}(M_{\chi_{2}}) \\
& \mathcal{L}(M_{\chi_{8}}) \ar"1,5"^<<<{E_{2578}} \ar"2,5"_<<{E_{1}} 
& & & \mathcal{L}(M_{\chi_{9}}) \ar"2,6"_{E_{3}} 
&
}
\end{align*}
we see that $D^1 \cap D^2 \cap D^3 = 0$ and
$\gcd(D^1_{23},D^2_{13},D^3_{12}) = 0$, while $\gcd(D^1_2,D^2_1) = 0$, 
$\gcd(D^2_3,D^3_2) = E_{46}$, $\gcd(D^1_3,D^3_1) = E_{689}$. Computing further 
$$ D^1 + \lcm(D^2_1,D^3_1) - \tilde{D}^3_2 - D^2 = 
E_1 + E_{2356789} - E_3 - E_{2578} \\ 
= E_{169} $$
and similarly  
$$ D^2 + \lcm(D^1_2, \tilde{D}^3_2) - \tilde{D}^1_3 - D^3 = 
E_{2578} + E_{13} - E_1 - E_3
= E_{2578}. $$
Consulting Figure \ref{figure-03} to compute intersections we see 
that $$\gcd(D^2_3,D^3_2) \cap (D^1 + \lcm(D^2_1,D^3_1) - \tilde{D}^3_2
- D^2) =  E_{46} \cap E_{169} = E_6$$ and 
$$ \gcd(D^1_3,D^3_1) \cap (D^2 + \lcm(D^1_2, \tilde{D}^3_2) -
\tilde{D}^1_3 - D^3) = E_{689} \cap E_{2578} = E_8.$$
Therefore by Lemma \ref{lem:cube_cohomology} we have
$\supp H^{-1}(\Psi(\mathcal{O}_0 \otimes \chi_{11})) = 
E_6 \cup E_8$ while
$\supp H^{0}(\Psi(\mathcal{O}_0 \otimes \chi_{11})) =
\supp H^{-2}(\Psi(\mathcal{O}_0 \otimes \chi_{11})) = 0$. 

Computing thus for each $\Psi(\mathcal{O}_0 \otimes \chi)$, we obtain
finally the data presented in 
Table \ref{table-13-1-5-7-supports-of-the-images} in the Introduction. 

\appendix

\section{Locally free resolutions of finite-length sheaves on
$\mathbb{C}^n$}

In this section we present a general method to resolve any
finite-length sheaf on $\mathbb{C}^n$ by locally free sheaves. It is
sufficiently intrinsic to lift to a relative setup of a flat family
of finite-length sheaves over a scheme of finite type over
$\mathbb{C}$. This result may be of independent interest. 

We then use this result to give an alternate proof of Proposition 
\ref{prps-Fdual} where we use this resolution to explicitly compute
the dual $\mathcal{F}^\vee$ in $D(Y)$.

\subsection{Resolving finite-length sheaves on $\mathbb{C}^n$ by
locally frees} 

The first step is the following technical lemma:

Let $A$ be a Hopf $\mathbb{C}$-algebra with a bijective antipode map. 
Denote by $\Delta$ its  comultiplication map, by $\epsilon$ the counit 
and by $S$ the antipode map. Let $W$ be a left $A$-module, which is
finitely generated as a $\mathbb{C}$-module. Denote by $M_{W}$ 
the left $A$-module structure on the $\mathbb{C}$-module $A \otimes W$ 
given by, in Sweedler's sigma notation 
(\cite{Sweedler-HopfAlgebras}, Section 1.2):
\begin{align*}
a \cdot (b \otimes w) = \sum bS(a_1) \otimes a_2 \cdot w
\end{align*}
Denote by $M'_{W}$ the structure of a free left $A$-module on $A
\otimes W$ where $A$ acts on the first factor of the tensor product 
by left multiplication.

Denote by $\alpha$ the surjective $\mathbb{C}$-module map 
$A \otimes W \twoheadrightarrow W$ given by $a \otimes w
\mapsto a \cdot w$. Denote by $\beta$ the surjective map  
$A \otimes W \twoheadrightarrow \mathbb{C} \otimes_A M_{W}$, 
where the $A$-module structure on $\mathbb{C}$ being given by the counit. 
We claim that $\alpha$ filters through $\beta$. 
Indeed, the kernel of $\beta$ is generated by the elements of form
$$ \sum (b S(a_1) \otimes a_2 \cdot w) - \epsilon(a) b \otimes w $$
for all $a,b \in A$ and $w \in W$. Rewriting as
\begin{align*}
&\sum (b S(a_1) \otimes a_2 \cdot w) - b \otimes (\sum
S(a_1)a_2) \cdot w = \\
&= \sum (b S(a_1) \otimes a_2 \cdot w - b \otimes (S(a_1) 
a_2) \cdot w)
\end{align*}
we see that each such element lies also in the kernel of $\alpha$. 
We thus have a commutative diagram:
\begin{align} \label{eqn-the-link-between-two-quotient-structures}
\xymatrix{
A \otimes_{\mathbb{C}} W 
\ar@{>>}[r]^{\beta} \ar@{>>}@/_1pc/[rr]_{\alpha} &
\mathbb{C} \otimes_A M_{W} \ar@{>>}[r]^{\gamma} 
& W
} .
\end{align}

\begin{lemma} \label{lemma-hopf-algebra-generalities}
\begin{enumerate}
\item \label{item-free-modules} $M_{W} \simeq M'_{W}$ as left $A$-modules.
\item \label{item-natural-transformation-is-isomoprhism} The map
$\gamma$ in \eqref{eqn-the-link-between-two-quotient-structures}
is a $\mathbb{C}$-module isomorphism.
\end{enumerate}
\end{lemma}
\begin{proof}
\ref{item-free-modules}). The isomorphism is given by the map $M'_{W}
\rightarrow M_{W}$ which sends $a \otimes w$ to $\sum S(a_1) \otimes a_2
\cdot w$. Its inverse is the map $M_{W} \rightarrow M'_{W}$ which
sends $a \otimes w$ to $\sum S^{-1}(a_2) \otimes a_1 \cdot w$. This 
can be verified using \cite{Sweedler-HopfAlgebras}, Prop. 4.0.1,
identities 3) and 4), and the identity 4) (iv) of the Exercises 
which immediately follow.

\ref{item-natural-transformation-is-isomoprhism}). From 
\ref{item-free-modules}) we have that $M_{W}$ is free. Hence 
$\mathbb{C} \otimes_A M_W$ is isomorphic to $W$ as a
$\mathbb{C}$-module. As $\mathbb{C} \otimes_A M_W \xrightarrow{\gamma} W$ 
is surjective it must also be an isomorphism.
\end{proof}

Let $\mathcal{W}$ be a finite length $G$-sheaf on $\mathbb{C}^n$. 
Then $\Gamma(\mathcal{W})$ is a finite dimensional vector space with a 
$G$-action. Denote by $\mathcal{L}$ the locally free 
$G$-sheaf $\mathcal{O}_{\mathbb{C}^n} \otimes
\Gamma(\mathcal{W})$. Its $\mathcal{O}_{\mathbb{C}^n}$-module 
structure is the natural action on the first factor of the tensor product,
while the action of $G$ is on both the factors. 

\begin{prps} \label{prps-lffr-resolution-of-finite-sheaf-over-c}
$\mathcal{W}$ has a $G$-equivariant locally free resolution of
given by 
\begin{align} \label{eqn-finite-sheaf-lffr-resolution-point}
0 \rightarrow \wedge^n \givrep \otimes \mathcal{L} \xrightarrow{\delta_n} \dots
\xrightarrow{\delta_{k+1}} \wedge^k \givrep \otimes \mathcal{L}
\xrightarrow{\delta_k} \dots \xrightarrow{\delta_1} \mathcal{L}
\xrightarrow{\alpha} \mathcal{W} \rightarrow 0
\end{align}
where the maps $\delta_k$ are defined for $m_i \in \givrep^\vee$, $f
\in \mathcal{O}_{\mathbb{C}^n}$ and $w \in \Gamma(\mathcal{W})$ by
\begin{align}
\label{eqn-differentials-of-lffr-resolution-of-finite-sheaves}
(m_1 \wedge \dots \wedge
m_k) \; f \otimes w \mapsto 
& \sum_{i = 1}^{k} (-1)^{i+1}  
(\dots m_{i-1} \wedge m_{i+1} \dots) 
(f \otimes m_i \cdot w - m_i f \otimes w)
\end{align}
and the map $\alpha$ is defined by $f \otimes w \mapsto f \cdot w$.
\end{prps}
\begin{proof}
Let $\mathcal{L}'$ be a $G$-sheaf on $\mathbb{C}^n$, 
which has the same underlying sheaf of $\mathbb{C}$-modules 
as $\mathcal{L}$, but whose $\mathcal{O}_{\mathbb{C}^n}$-module
structure is induced by the action 
\begin{align*}
m\cdot f\otimes w = f \otimes m \cdot w - m f \otimes w.
\end{align*}
of $\givrep^\vee$. The action of $G$ remains the same. 

It is easy to check that the differential maps $\delta_k$ of 
the complex $\eqref{eqn-finite-sheaf-lffr-resolution-point}$ respect this 
new $G$-sheaf structure. The resulting complex 
\begin{align*} 
0 \rightarrow \wedge^n \givrep \otimes \mathcal{L}'
\xrightarrow{\delta_n} \dots
\xrightarrow{\delta_{k+1}} \wedge^k \givrep \otimes \mathcal{L}'
\xrightarrow{\delta_k} \dots \xrightarrow{\delta_1} \mathcal{L}'
\rightarrow 0
\end{align*}
is precisely the Koszul complex of 
$\mathcal{L}'$ with respect to any basis of $\givrep^\vee$. 
Let $\mathfrak{m}$ be the ideal sheaf of $\mathcal{O}_{\mathbb{C}^n}$
generated by the elements of $\givrep^\vee$. For grossly general
reasons, as can be seen from Lemma
\ref{lemma-hopf-algebra-generalities},
\ref{item-free-modules}), $\mathcal{L}'$ is locally free. In
particular, then, it is a Cohen-Macaulay
$\mathcal{O}_{\mathbb{C}^n}$-module. Furthermore, the support of
$\mathcal{L'} / \mathfrak{m}\mathcal{L}'$ has codimension $n$, 
being the point $0 \in \mathbb{C}^n$. Any basis of $\givrep^\vee$, as
it consists of $n$ elements and generates $\mathfrak{m}$, must 
therefore be a regular sequence for $\mathcal{L}'$. Hence the 
$G$-sheaf complex 
\begin{align*} 
0 \rightarrow \wedge^n \givrep \otimes \mathcal{L}'
\xrightarrow{\delta_n} \dots
\xrightarrow{\delta_{k+1}} \wedge^k \givrep \otimes
\mathcal{L}'
\xrightarrow{\delta_k} \dots \xrightarrow{\delta_1} \mathcal{L}'
\xrightarrow{\beta} \mathcal{L}'/\mathfrak{m}\mathcal{L}'
\rightarrow 0
\end{align*}
is exact. It then follows from Lemma \ref{lemma-hopf-algebra-generalities},
\ref{item-natural-transformation-is-isomoprhism}), that the 
following is also exact as a complex of sheaves of
$\mathbb{C}$-modules:
\begin{align*} 
0 \rightarrow \wedge^n \givrep \otimes \mathcal{L}'
\xrightarrow{\delta_n} \dots
\xrightarrow{\delta_{k+1}} \wedge^k \givrep \otimes
\mathcal{L}'
\xrightarrow{\delta_k} \dots \xrightarrow{\delta_1} \mathcal{L}'
\xrightarrow{\alpha} \mathcal{W}
\rightarrow 0.
\end{align*}
As it is precisely the complex of sheaves of $\mathbb{C}$-modules 
underlying the $G$-sheaf complex
\eqref{eqn-finite-sheaf-lffr-resolution-point}, the complex
\eqref{eqn-finite-sheaf-lffr-resolution-point} must also be exact.
\end{proof}

Proposition \ref{prps-lffr-resolution-of-finite-sheaf-over-c} 
naturally generalises to flat families of finite sheaves. Let $S$ 
be a scheme of finite type over $\mathbb{C}$ endowed with a trivial 
$G$-action.  Let $\mathcal{F}$ be a coherent $G$-sheaf on 
$S \times \mathbb{C}^n$, flat over $S$, whose fiber $\mathcal{F}_{|p}$ 
over any closed point $p \in Y$ is a finite-length $G$-sheaf on 
$\mathbb{C}^n$. Imitating the construction above, let $\mathcal{L}$ be
the locally-free sheaf $\pi_S^* \pi_{S *} \mathcal{F}$ and let
$\mathcal{L} \overset{\alpha}{\twoheadrightarrow} \mathcal{F}$ be 
the surjection map defined by adjunction. 

There are two natural actions of $\regring$ on $\mathcal{L}$.  
The first action is via the inclusion $\regring \hookrightarrow
\Gamma(\mathcal{O}_{S\times\mathbb{C}^n})$ induced by the pullback 
from $\mathbb{C}^n$. The second action is via the the inclusion 
$\regring \hookrightarrow  \pi_{S *} \mathcal{O}_{S \times \mathbb{C}^n}$
which acts on $\pi_* \mathcal{F}$ and subsequently on 
$\mathcal{L} = \pi_S^* \pi_{S *} \mathcal{F}$.
For any $r \in \regring$ and $s \in \mathcal{L}$ we denote the
first action by $r \cdot s$ and the second action by $r \cdot_2 s$.

\begin{cor} \label{cor-lffr-resolution-of-finite-sheaf-family}
$\mathcal{F}$ has a $G$-equivariant locally-free resolution given by
\begin{align} \label{eqn-finite-sheaf-lffr-resolution-family}
0 \rightarrow \wedge^n \givrep \otimes \mathcal{L} 
\xrightarrow{\delta_n} \dots
\xrightarrow{\delta_{k+1}} \wedge^k \givrep \otimes
\mathcal{L} 
\xrightarrow{\delta_k} \dots \xrightarrow{\delta_1} \mathcal{L}
\xrightarrow{\alpha} \mathcal{F} \rightarrow 0
\end{align}
where the maps $\delta_k$ are defined for $m_i \in \givrep^\vee$ and $s
\in \mathcal{L}$ by
\begin{align}
\label{eqn-differentials-in-finite-sheaf-lffr-resolution-family}
(m_1 \wedge \dots \wedge
m_k) \otimes s \mapsto 
& \sum_{i = 1}^{k} (-1)^{i+1}  
(\dots m_{i-1} \wedge m_{i+1} \dots) \otimes
(m_i \cdot_2 s - m_i \cdot s)
\end{align}
\end{cor}

\subsection{Proof of Proposition \ref{prps-Fdual}}

As an application of the results in previous section, we use them to 
give an alternative proof of Prop. \ref{prps-Fdual}.

\begin{proof}[Proof of Prop. \ref{prps-Fdual}]
By Corollary \ref{cor-lffr-resolution-of-finite-sheaf-family}
we have a locally free resolution of $\mathcal{F}$. Taking 
its dual we find that $\mathcal{F}^\vee[n]$ is isomorphic in $D^G(Y
\times \mathbb{C}^n)$ to 
\begin{align*}
0 \rightarrow
\mathcal{L}^\vee \xrightarrow{\delta_1^\vee} 
\dots \xrightarrow{\delta_k^\vee} 
\wedge^k \givrep^\vee \otimes \mathcal{L}^\vee \xrightarrow{\delta_{k+1}^\vee} 
\dots \xrightarrow{\delta_n^\vee} 
{\wedge^n \givrep^\vee \otimes
\mathcal{L}^\vee}
\rightarrow 0
\end{align*}
where the rightmost term, ${\wedge^n \givrep^\vee \otimes
\mathcal{L}^\vee}$, lies in degree zero. 

On the other hand, by definition  $\pi_{Y*} \tilde{\mathcal{F}} = (\pi_{Y*} \mathcal{F})^\vee$, which is locally free. 
Hence $\pi^*_{Y} \pi_{Y*} \tilde{\mathcal{F}} = \mathcal{L}^\vee$.
Applying Corollary \ref{cor-lffr-resolution-of-finite-sheaf-family} 
again we find that $\tilde{\mathcal{F}}$ has a locally free resolution
\begin{align*}
0 \rightarrow 
\wedge^n \givrep \otimes \mathcal{L}^\vee 
\xrightarrow{\gamma_n} \dots
\xrightarrow{\gamma_{k+1}} \wedge^k \givrep \otimes
\mathcal{L}^\vee 
\xrightarrow{\gamma_k} \dots \xrightarrow{\gamma_1}
{\mathcal{L}^\vee}
\rightarrow 0
\end{align*}
where differential maps $\gamma_i$ are 
defined by \eqref{eqn-differentials-in-finite-sheaf-lffr-resolution-family}, 
with the action $\cdot_2$ of $\regring$ on $\mathcal{L}^\vee$ being
induced by the action \eqref{eqn-natural-twalg-module-structure-on-dual} 
of $\regring$ on $\shhomm_{\mathcal{O}_Y}(\pi_*\mathcal{F}, \mathcal{O}_Y)$. 

We now prove that these two complexes are isomorphic, thus showing
$\mathcal{F}^\vee[n] \simeq \tilde{\mathcal{F}}$. 
Let $\phi$ be the interior product 
isomorphism $\wedge^k \givrep^\vee \otimes \wedge^{n} 
\givrep \rightarrow \wedge^{n-k} \givrep$. It is $G$-equivariant.
Since $G \subseteq \gsl_n(\mathbb{C})$, $\wedge^n \givrep$ is
trivial as a representation of $G$. Thus we obtain $G$-equivariant 
isomorphisms $\wedge^k \givrep^\vee \xrightarrow{\sim}
\wedge^{n-k} \givrep$ which induce 
isomorphisms $\alpha_i:\; \wedge^{n+i} \givrep^\vee \otimes
\mathcal{L}^\vee \xrightarrow{\sim}  \wedge^{-i} \givrep \otimes \mathcal{L}^\vee$. 

It remains only to verify that $\alpha_i$ define a chain map 
between the two complexes, i.e. that the following diagram commutes
\begin{align*}
\xymatrix{
\wedge^k \givrep^\vee \otimes \mathcal{L}^\vee 
\ar[r]^{\delta_k^\vee} \ar[d]_{\alpha_{k-n}} &
\wedge^{k+1} \givrep^\vee \otimes \mathcal{L}^\vee
\ar[d]^{\alpha_{k+1-n}} \\
\wedge^{n-k} \givrep \otimes \mathcal{L}^\vee 
\ar[r]^{\gamma_{n-k}} &
\wedge^{n-k+1} \givrep \otimes \mathcal{L}^\vee 
}.
\end{align*}
It is a straightforward calculation which we leave to the reader.
\end{proof}

\bibliography{references}
\bibliographystyle{amsalpha}

\end{document}